\documentclass[leqno]{article}
\usepackage{amssymb,latexsym,amscd,pst-all,mathrsfs}
\newcommand{\se}[1]{{\section{#1}} {\setcounter{equation}{0}}}
\newtheorem{prop}{Proposition}[section]
\newtheorem{theorem}[prop]{Theorem}
\newtheorem{lm}[prop]{Lemma}

\newtheorem{co}[prop]{Corollary}
\newtheorem{re}[prop]{Remark}

\def\k{{K\"{a}hler}}

\input epsf
\begin{document}
\hbadness=10000
\title{{\bf The Fukaya category of symplectic neighborhood of a non-Hausdorff manifold}}
\author{Wei-Dong Ruan\\}
\maketitle

\se{Introduction}
In this paper, using similar idea as in \cite{FO}, we devise a method to compute the Fukaya category of certain exact symplectic manifold by reducing it to the corresponding Morse category of a non-Hausdorff manifold as perturbation of the Lagrangian skeleton of the exact symplectic manifold. The main theorem is theorem \ref{mk}.\\

\stepcounter{subsection}
{\bf \S \thesubsection\ Background:} Symplectic geometry (topology) is one of the fastest growing branch of geometry (topology) that exhibits both the characteristics of geometry and topology. Most of the new developments start with the well known fact that for any symplectic manifold $(M, \omega)$, there exists a compatible almost complex structure $J$ on $M$ (usually not integrable) that gives rise to the natural almost K\"{a}hler structure $(M, J, g, \omega)$. In the following we will discuss 3 such developments relevant to us that profoundly altered the landscape of symplectic geometry.\\

{\bf (1)} In a seminal work \cite{Gromov}, Gromov realized that although higher dimensional pseudo holomorphic submanifolds of an almost K\"{a}hler manifold are non-generic, in the one-dimensional case, pseudo holomorphic curves are generic and usually form finite dimensional moduli space that can be used to understand symplectic topology of the symplectic manifold. This important work, together with the work of Witten from physics perspective, through the great efforts of many others, grew into the highly successful Gromov-Witten theory of pseudo holomorphic curves that also in turn greatly impacted the complex and algebraic geometry. Gromov also considered the case of pseudo holomorphic disk with Lagrangian boundary condition that through works of Floer and many others was developed into the Floer theory, which is more relevant to us in this paper.\\

{\bf (2)} In a separate development \cite{D1}, Donaldson realized that although pseudo holomorphic hypersurfaces are non-generic, one can still construct smooth almost holomorphic hypersurfaces that give rise to a lot of smooth symplectic hypersurfaces, through similar technique as the construction of complex hyperplane sections in complex geometry. In \cite{D2}, Donaldson further constructed a family of almost holomorphic sections that gives rise to symplectic Lefschetz pencil, which through the works of Donaldson and many others have been applied with great success to the understanding of symplectic 4-manifolds, where coincidentally, Gromov's pseudo holomorphic curves and Donaldson's almost holomorphic symplectic hypersurfaces meet.\\

{\bf (3)} In yet another separate development \cite{EG}, Eliashberg and Gromov discussed the Lagrangian skeleton of convex symplectic open manifolds that is related to the Lagrangian skeleton of Stein manifolds studied by Eliashberg. P. Biran (\cite{Biran1}) and others studied the special case of complement of hyperplane sections in a projective algebraic manifold as convex symplectic open manifold, and yielded many interesting results. (It is worth mentioning that the almost holomorphic hyperplane sections constructed by Donaldson in \cite{D1} is ideal to be used to figure out Lagrangian skeleton for general symplectic manifolds.)\\

The greatest common strength of all these three very important developments is their general applicability to ALL symplectic manifolds, which makes them extremely powerful tools in symplectic geometry. While Gromov-Witten invariants of pseudo holomorphic curves (more precisely, the quantum cohomology) can be effectively computed quite generally owe to the powerful machinery of algebraic geometry, there are very few techniques and results available that can effectively compute Gromov-Witten-Floer invariants of pseudo holomorphic disks with Lagrangian boundary conditions (more precisely, the Fukaya categories) in general.\\

To compute Fukaya category effectively, certain structure theorem of symplectic manifolds is needed. Donaldson's construction of almost holomorphic symplectic hypersurfaces and Lefschetz pencils provide such kind of structure theorems. Based on his work of symplectic Lefschetz pencil, Donaldson devised a process to compute the Fukaya category inductively through symplectic Lefschetz pencil. Such process was carried out in P. Siedel's proof of homological mirror symmetry for quartic K3 surfaces (\cite{S}). However, due to its inductive nature, this process will become increasingly complicated as dimension increases. (The dimension 2 case is already quite involved.)\\

Our work starts from the view point that Donaldson's Lefschetz pencil and Lagrangian skeleton of Eliashberg and Gromov (both can be constructed via certain almost holomorphic sections) can be viewed as two sides of a coin. One should be able to uncover the symplectic geometry (topology) of the symplectic manifold through either one of the two complementary structures. On the other hand, the mechanisms of the two sides are quite different. Donaldson's Lefschetz pencil approach inductively reduces higher dimensional symplectic geometry (topology) to lower dimensional ones, while Lagrangian skeleton reduces the symplectic geometry (topology) of the symplectic manifold to the real geometry (topology) of the Lagrangian skeleton. The purpose of this paper is to devise an effective method to compute Fukaya category through Lagrangian skeleton. It turns out that the complexity of the computation through Lagrangian skeleton does not increase significantly with dimension, which makes it potentially valuable especially for computing the Fukaya category of higher dimensional symplectic manifolds.\\

In \cite{sum}, we give an alternative construction of Donaldson's symplectic hypersurface and Lefschetz pencil, which can be expressed rather explicitly according to a ball cover of the symplectic manifold. The Lagrangian skeleton will be constructed in a sequel of \cite{sum}, which can also be expressed quite explicitly. In principle, results in this paper should be helpful for computing the Fukaya category of a general symplectic manifold through such Lagrangian skeleton.\\

The homological mirror symmetry conjecture by Kontsevich \cite{KO} is an important mathematical conjecture grows out of the physical mirror symmetry conjecture, that identifies the Fukaya category in symplectic geometry with the derived category of coherent sheaves through mirror symmetry. The main difficult is in the symplectic side due to the lack of understanding and means of computing the Fukaya category. This work will be used in my joint work with A. Bondal \cite{BR} to prove Kontsevich's homological mirror symmetry conjecture for weighted projective space of general dimension generalizing works of Auroux, Kazarkov, Orlov \cite{AKO} and Seidel \cite{S1, S2}. This work can also be used to generalize Seidel's proof of the homological mirror symmetry for quartic K3 \cite{S} to higher dimension Calabi-Yau \cite{CY}.\\

\stepcounter{subsection}
{\bf \S \thesubsection\ Setting:} A symplectic manifold $(M, \omega)$ is called exact if $\omega = d\alpha$ for a 1-form on $M$. An exact symplectic manifold is necessarily open. Major examples of exact symplectic manifolds are cotangent bundle of a real manifold and K\"{a}hler manifolds with global K\"{a}hler potential, in particular, when $M$ is an affine algebraic manifold and $\omega$ is the K\"{a}hler form of a complete K\"{a}hler metric on $M$ with global K\"{a}hler potential $h$ so that $h^{-1}((-\infty,c])$ is compact in $M$ for any $c\in \mathbb{R}$.\\

Without loss of generality, we may assume that $h$ is a Morse function. Let $\mathbb{L}$ be the union of the stable manifolds of the gradient flow of $h$ with respect to the K\"{a}hler metric. It is straightforward to check that $\omega|_\mathbb{L}=0$. $\mathbb{L}$ is the isotropic skeleton of $(M, \omega)$. When $\mathbb{L}$ is of middle dimension, it is called the Lagrangian skeleton of $(M, \omega)$. We assume the existence of isotropic (Lagrangian) skeleton in this paper.\\

As pointed out in \cite{EG}, isotropic (Lagrangian) skeleton exists in more general situations, for example, when the exact symplectic manifold is of Weinstein type. Let $(M, \omega = d\alpha)$ be an exact symplectic manifold. Let $\phi_t$ be the flow generated by the vector field $v$ satisfying $\imath_v \omega = \alpha$. Notice that the flow is exactly the gradient flow of $h$ for the K\"{a}hler case. The purpose of the Weinstein type condition is to ensure that $\phi_t$ behave like a gradient flow. In general,

\[
{\cal L}_v \alpha = \imath_v d\alpha + d\langle \alpha, v\rangle =\alpha,\ \frac{d\phi_t^* \alpha}{dt} = \phi_t^* \alpha,\ \phi_t^* \alpha = e^t \alpha.
\]

Let $\mathbb{L} \subset M$ be the corresponding Lagrangian skeleton. Since $v$ is tangent to $\mathbb{L}$ and $\omega|_\mathbb{L}=0$, we have $\alpha|_\mathbb{L}=0$. In particular, $\mathbb{L}$ is exact. (In the simplest situation of $M = T^*L$ being the cotangent bundle of a real manifold $L$ with the canonical exact symplectic structure, the Lagrangian skeleton $\mathbb{L} = L$ is smooth.)\\

Let $L_0$ be an exact Lagrangian submanifold, namely, $\alpha|_{L_0} =df_0$. Then $\imath_v \omega|_{\phi_t (L_0)}$ $= d (e^tf_0)$. Hence $\phi_t (L_0)$ is a Hamiltonian deformation family starting from $L_0$, and approaches the Lagrangian skeleton $\mathbb{L}$ when $t\rightarrow +\infty$.\\

Heuristically, under such limiting process, the pseudo holomorphic polygons needed to compute the Fukaya category will degenerate to certain gradient tree on the Lagrangian skeleton in similar fashion as in the much simpler situation of $M$ being the cotangent bundle of a real manifold in \cite{FO}. In such way, the Fukaya category of the exact symplectic manifold $(M, \omega)$ can be identified with certain Morse category of the Lagrangian skeleton $\mathbb{L}$, therefore can be computed through counting of rigid gradient trees when the Lagrangian skeleton is reasonably structured.\\

It will be ideal and more canonical to make sense of such Morse category directly with objects being the piecewise smooth top dimensional submanifolds in $\mathbb{L}$. Since Fukaya category was only defined for smooth Lagrangian submanifolds, for technical reasons from analysis, such more canonical approach is very difficult to realize.\\

Instead, difficulty with singularities can be avoided by perturbing the singular Lagrangian skeleton $\mathbb{L}$ into a smooth Lagrangian non-Hausdorff submanifold ${\mathscr L}$ that is exact, which we will still call Lagrangian skeleton for simplicity. (See remark \ref{bi} for non-Hausdorff manifold.) The non-Hausdorff submanifold perturbation ${\mathscr L}$ is chosen in the way such that the exact Lagrangian submanifolds that is relevant to the specific $A^\infty$-product computation can all be realized as smooth Hausdorff submanifolds of ${\mathscr L}$. For most of our applications, there is a natural choice of the non-Hausdorff submanifold perturbation that is enough to compute all the $A^\infty$-products among the generators of the category. In more complicated cases, one may need to use different non-Hausdorff submanifold perturbations to compute different $A^\infty$-products. But The technical difficulty for all this cases are essentially the same and is discussed in this work. (We should mention here that \cite{KO} with very different purpose in mind does show some close resemblance and relation to a special case of our work.)\\

\stepcounter{subsection}
{\bf \S \thesubsection\ The relation to \cite{FO}:} Work of Fukaya and Oh (\cite{FO} that bears the influences of Gromov, Floer, Donaldson, Kontsevich, etc.) contributed greatly to our understanding of Fukaya category and its relation to classical Morse category, and was quite ahead of its time. Our work is greatly inspired by \cite{FO}, and is very much following the foot steps of \cite{FO} closely in general framework and many ideas. To better understand our work, it should be helpful to consult \cite{FO} for background, references and the origin of many of these ideas. It is our hope that our work could contribute in a small way to the improvement and generalization of this important work.\\

{\bf (1)} On the other hand, our method is quite different from the traditional approach starting from Floer \cite{floer2} and is subsequently followed in Fukaya-Oh \cite{FO}. First of all, we use H\"{o}lder norms (More precisely, $C^{1,\alpha}$-norm) for our estimates, which is quite different from the traditional approach using various Sobolev norms that owe to the great success of psedo-differential operators via the advanced technology of Fourier transformation. Secondly, rescaling method and Cheeger-Gromov convergence is used extensively and systematically in our work. Thirdly, elementary inverse function theorem is used to prove the existence results of pseudo holomorphic disks, instead of more sophisticated method of parametrix used in \cite{FO}.\\

{\bf (1a)} Being pointwise norms, H\"{o}lder norms are more classical than Sobolev norms and are in a sense more natural geometrically. One obvious disadvantage of H\"{o}lder norms is that they can not take advantage of the powerful pseudo-differential operator estimates via Fourier transformation, especially for higher derivative estimates, Sobolev norms are more flexible to use. One important motivation for us to use the H\"{o}lder norms is the realization that the purpose of studying pseudo holomorphic disks is to count them in order to understand the symplectic topology. Consequently, higher derivative estimates of pseudo holomorphic disks are not absolutely necessary for such topological purpose. It turns out, quite surprisingly to the author, $C^{1,\alpha}$-smoothness of pseudo holomorphic disks is all one need to count the pseudo holomorphic disks correctly and establish the Fukaya category in the context of our work. It also turns out, quite surprisingly to the author, that all these $C^{1,\alpha}$-estimates can be derived in quite elementary fashion using nothing more than Cauchy integration formula and its generalizations, in another word, via the potential theory.\\

Pseudo-holomorphic disk can be viewed as generalization of holomorphic disk or special case of minimal surface. It is interesting to see how much of pseudo holomorphic disk theory can be reduced to the theory of one variable holomorphic functions, which is based on simple things like Cauchy integration formula and is widely known, and how much of pseudo holomorphic disk theory has to resort to minimal surface theory, which in general involve quite sophisticated and difficult non-linear PDE methods that is only accessible to experts. Or as termed in Oh \cite{Oh} (with not necessarily the same exact intended meaning), holomorphicity versus harmonicity. Our work indicates that estimates of pseudo holomorphic disks in our context can be derived via holomorphicity instead of harmonicity.\\

{\bf (1b)} The method of rescaling and Cheeger-Gromov convergence was pioneered by Gromov and widely used since then. It is a very simple idea that turns out to be extremely powerful when used in combination with suitable estimates. It is also used for some of the arguments in \cite{FO}. The extensive and systematical use of this method in this work is crucial in reducing our arguments to the basic $C^{1,\alpha}$-estimates. It is interesting to notice that the Cheeger-Gromov convergence results needed for such arguments are of the most classical and standard type (see remark \ref{bk}).\\

{\bf (1c)} Implicit function theorem is among the simplest and most direct methods in proving existence results, provided suitable global norms can be constructed that enable one to prove the linearized operator is invertible or more precisely a quasi-isometry between Banach spaces. This method enable us to avoid the delicate gluing construction of the right inverse operator in \cite{FO}. In this paper, gluing argument is only used to construct the approximate solutions.\\

{\bf (2)} \cite{FO} mainly consists of 2 parts. The first part shows the convergence of pseudo holomorphic disks to gradient tree. The second part shows the construction and uniqueness of a holomorphic disk near a gradient tree. Our work concerns symplectic neighborhood of non-Hausdorff Lagrangian submanifold generalizing \cite{FO}. The generalization from \cite{FO} mainly occurs in the first part. For the second part, since gradient tree is contractible, it is always inside a Hausdorff Lagrangian submanifold $L$ in the non-Hausdorff Lagrangian submanifold ${\mathscr L}$. In principle, result in the second part of \cite{FO} (suitably modified) can be applied to construct the pseudo holomorphic polygon corresponding to the gradient tree.\\

Since the norm used in the first part of our work is different from \cite{FO}, to apply the second part of \cite{FO} to our situation, additional arguments would be needed to address the transition between 2 types of norms. Instead, in our second part, we provide an alternative proof of the existence and uniqueness result for the pseudo holomorphic polygon using $C^{1,\alpha}$-norm and inverse function theorem. Beside satisfying our curiosity on testing the potential and limitation of our elementary approach via $C^{1,\alpha}$-norm, the more important reason is that we need to treat a special case of gradient tree (called exceptional gradient tree here) that was ignored in \cite{FO} and is very important for our application later to homological mirror symmetry. A less substantial difference, which is also necessary to be addressed is that \cite{FO} was proved for the canonical almost complex structure of $M= T^*L$ determined by a metric on $L$, while we only assume the Lagrangian fibres to be orthogonal to $L$.\\
\begin{re}
In our proof, we also need the general property of pseudo holomorphic polygon that $\displaystyle \sup_{t\in \Theta} |Dw(t)|$ is finite, which, for example, is a consequence of the more precise estimate of lemma 9.3 in \cite{FO}. It comes down to the $C^{0,1}$-estimate according to the area bound, which is the $W^{1,2}$-bound of the pseudo holomorphic disk. Sobolev norm clearly has advantage for such result. This is the only estimate we use that we have to resort to the traditional method with Sobolev norm. The only result in this paper not fully proved is the index result of the linearized operator necessary to prove the surjectivity of the linearized operator. We only showed how the index result from \cite{FO} (not including the case of exceptional gradient tree) implies the index result in our case. A full proof of the index result (including the case of exceptional gradient tree) using a geometric $C^{1,\alpha}$ method will be provided in a forthcoming paper \cite{disk}, because the proof has its independent flavor and fit naturally into a different context.
\end{re}

\stepcounter{subsection}
{\bf \S \thesubsection\ Organization:} The layout of this paper is similar to \cite{FO}. In section 2, we explain in detail and prove the main result (theorem \ref{mk}) based on results from later sections.\\

Sections 3, 4 and 5 discuss basic results and estimates necessary mainly for section 6 and also for later sections. Section 3 contains basic results for Cheeger-Gromov convergence. Section 4 contains elementary discussion of conformal models of the disk and some basic facts of holomorphic functions of one complex variable that we need. Section 5 contains the basic interior and boundary $C^{1,\alpha}$ estimates of pseudo holomorphic disks and corresponding convergence results.\\

Rigid pseudo holomorphic polygons are related to rigid gradient trees in section 6 (theorem \ref{mo}), which effectively reduces our discussion from the non-Hausdorff situation to the Hausdorff situation (see remark \ref{moa}).\\

In section 7, the approximate pseudo holomorphic polygon necessary for the existence result is constructed from a rigid gradient tree with estimates derived.\\

Sections 8 and 9 establish the framework and basic estimates for proving the existence of pseudo holomorphic polygon through inverse function theorem.\\

The invertibility of the linearized operator is proved in section 10 (theorem \ref{pf}), therefore finishes the proof of the existence of pseudo holomorphic polygon. The uniqueness of pseudo holomorphic polygon near a rigid gradient tree is proved in section 11 (theorem \ref{qa}), therefore finishes the proof of our main theorem.\\

{\bf Convention of notations:} $z = z^\Re + i z^\Im$ denotes the decomposition of a complex number into its real and imaginary parts. ``$\lim$" means ``$\displaystyle\lim_{k \rightarrow +\infty}$" unless otherwise specified. $A = O(B)$ if there is universal constant $C>0$ such that $|A| \leq C|B|$. $A\sim B$ if $A = O(B)$ and $B = O(A)$. $A_k = o(B_k)$ if $\lim (|A_k|/|B_k|) =0$ uniformly. $C$ and its variations are used as positive constants for estimates that may differ from expression to expression. $\nabla f$ (resp. ${\cal H}f$) is used to denote the gradient (resp. the Hessian) of a function $f$. Each of $\alpha$, $i$, $C$, etc. has been used to represent different meanings that are unlikely to be confused from the contexts.\\

{\bf Acknowledgement:} This work grows out of the need of generalizing Fukaya-Oh's work \cite{FO} for application to my joint work \cite{BR} with A. Bondal on proving the homological mirror symmetry for weighted projective spaces. I am very grateful to Alexei for many very inspiring discussions. I would also like to thank P. Seidel, Yau and K. Zhu for helpful discussions.\\

\se{Fukaya category and Morse category}
Being an active part of the youthful symplectic geometry, Fukaya category is still in its infancy. Despite great effort by many people, the precise definition in general is still urgently needing to be pined down and worked out. Very few systematical methods are available to compute Fukaya category effectively. This is not to say we know nothing about Fukaya category. In fact, a great deal of general structure and properties of Fukaya category is already well understood. We also have quite complete understanding of some important subcategories of the Fukaya category in various cases. A rather effective method is to use our general understanding of the Fukaya category to show that certain well understood subcategory in effect generate the full Fukaya category. (In a remarkable success story, P. Seidel was able to use such method to determine the Fukaya category for the quartic K3 surface in \cite{S}.)\\

Fukaya category of a general symplectic manifold in its simplest version should contain Lagrangian submanifolds as objects. The morphism group between such Lagrangian objects in general position should be generated by the Lagrangian intersection points. Fukaya category is a $A^\infty$-category, namely, there are a sequence of $A^\infty$ products among the morphism groups satisfying $A^\infty$ compatibility conditions that can be computed via weighted counting of holomorphic polygons (disks) with boundary in the Lagrangian objects. The $A^\infty$ structure satisfies certain definite homotopy properties under Hamiltonian deformation of the Lagrangian objects, which more or less reduce the computation of the $A^\infty$ products to among the Lagrangian objects in general positions. The main invariants of the Fukaya category are encoded in the $A^\infty$ products. (We should remark here that although the general definition of Fukaya category is still under construction, such symplectic invariants for Lagrangian submanifolds modulo Hamiltonian deformations including the Floer homology as special case can already be rigorously defined.) Whatever the precise context and precise definition of the Fukaya category might be, the most important ingredient in our understanding of Fukaya category is the computation of these $A^\infty$ products. \cite{FO} did just such computation for exact Lagrangian objects in the special case of cotangent bundle of a real manifold as a symplectic manifold by reducing it to the $A^\infty$ products in the corresponding Morse category of the real manifold, which can be computed by classical means.\\

Therefore, in this paper, instead of worrying about the general definition of Fukaya category, we will concentrate on the computation of such $A^\infty$ products, which is rigorously defined for general positioned Lagrangian objects. More precisely, we will compute the $A^\infty$ products for general positioned exact Lagrangian objects in convex symplectic manifold (as termed by Eliashberg and Gromov) that is wrapped around so-called Lagrangian skeleton of Eliashberg and Gromov. Analogue to the approach of \cite{FO}, we identify the $A^\infty$ product of the Fukaya category with the corresponding $A^\infty$ products in certain Morse category of the Lagrangian skeleton. In some cases, when Fukaya-type category can be defined rigorously, and the full category is generated by such Lagrangian objects, (for example, the Fukaya-Seidel category in \cite{S1,S2}), our computation will be able to determine the full structure of the Fukaya-type category. \\

For technical reason, it is convenient for us to define the Morse category on the Hamiltonian perturbation of the singular Lagrangian skeleton into a smooth non-Hausdorff Lagrangian submanifold (skeleton). Such perturbation is usually straightforward to be done in a case by case basis. In this paper, we will start with such smooth non-Hausdorff perturbation of the Lagrangian skeleton, which is necessarily exact as indicated in the introduction, and consider Lagrangian submanfolds as small Hamiltonian deformation of top dimensional Hausdorff submanifolds of the smooth non-Hausdorff Lagrangian skeleton.\\

{\samepage
\stepcounter{subsection}
{\bf \S \thesubsection\ Symplectic neighborhood of a non-Hausdorff manifold}\\
\nopagebreak

We start with} the following version of Darboux-Weinstein neighborhood theorem.\\
\begin{prop}
\label{mf}
For a $n$-manifold $L$ and a Lagrangian immersion $\imath_L: L \rightarrow (M,\omega)$,and a Lagrangian subbundle $N$ of $\imath_L^* TM$ that is transverse to $TL$ as subbundle of $\imath_L^* TM$. There exists a tubular neighborhood $U_L$ of $L$ in $(T^*L, \omega_L)$ together with the Lagrangian fibration $\pi_L: U_L \rightarrow L$, so that $\imath_L$ can be extended to an open symplectic immersion $e_L: (U_L, \omega_L) \rightarrow (M, \omega)$, under which, $T^*L$ as subbundle of $TU_L|_L$ is identified with $N$ as subbundle of $(e_L^* TM)|_L = \imath_L^* TM$.
\end{prop}
{\bf Proof:} The proof of this result is a straightforward generalization of the proof of the standard Darboux-Weinstein neighborhood theorem.\\

Since the proof is local in nature, without loss of generality, we may assume that $\imath_L$ is an embedding. $\omega$ induce an identification $N \cong T^*L$. It is straightforward to extend this identification to a smooth open embedding $i_L: U_L \hookrightarrow M$ for a tubular neighborhood $U_L$ of $L$ in $T^*L$, so that $T^*L$ as subbundle of $TU_L|_L$ is identified with $N$ as subbundle of $(i_L^* TM)|_L = \imath_L^* TM$.

\[
\omega_t = (1-t) \omega_L + t i_L^* \omega
\]

defines a family of symplectic forms on $U_L$. There exists a smooth 1-form $\alpha$ on $U_L$ such that $\omega_1 - \omega_0 = d\alpha$.\\

Through adjustment of exact 1-form on $U_L$, one can ensure that $\alpha = O(|p|^2)$ as a 1-form on $U_L$, where $L = \{p=0\}$.\\

Let $\phi_t$ be the flow generated by the family of vector fields $\{v_t\}$ such that $\imath_{v_t} \omega_t = - \alpha$. Then $\phi_t^* \omega_t = \omega_0$ for $t \in [0,1]$. In particular, $e_L = i_L \circ \phi_1: (U_L, \omega) \rightarrow (M, \omega)$ gives the desired open symplectic embedding.
\hfill$\Box$\\

For any subset $S \subset L$, define $U_S := U_L \cap T^*L|_{S}$. We have the following refinement of proposition \ref{mf}.

\begin{prop}
\label{mh}
In proposition \ref{mf}, assume that $N_{x_1} = N_{x_2}$ for any $x_1,x_2 \in L$ such that $\imath_L (x_1) = \imath_L (x_2)$. $e_L$ can be constructed in such a way so that $e_L (U_{x_1}) = e_L (U_{x_2})$ for any $x_1,x_2 \in L$ satisfying $\imath_L (x_1) = \imath_L (x_2)$. Furthermore, if $(\imath_L)_* T_{x_1} L$ $= (\imath_L)_* T_{x_2} L$, then the affine structures on the Lagrangian fibre $e_L (U_{x_1}) = e_L (U_{x_2})$ determined through $U_{x_1}$ and $U_{x_2}$ coincide.
\end{prop}
{\bf Proof:} Let $L_1$ (resp. $L_2$) be the image of a neighborhood of $x_1$ (resp. $x_2$) in $L$ under the map $\imath_L: L \rightarrow M$. Let $L_{12} = L_1 \cap L_2$. We may first construct $e_{L_1}: (U_{L_1}, \omega_{L_1}) \rightarrow (M, \omega)$ according to proposition \ref{mf}. Then extend $e_{L_1} |_{U_{L_{12}}}$ to $i_{L_2}$ on $U_{L_2}$. Then the $\alpha = O(|p|^2)$ can be adjusted by exact 1-form so that $\alpha|_{U_{L_{12}}}=0$. Consequently, the Lagrangian fibres under $e_{L_1}$ and $e_{L_2}$ coincide over $L_{12}$. When $\imath_L$ is constructed in such a way, $e_L (U_{x_1}) = e_L (U_{x_2})$ for any $x_1,x_2 \in L$ satisfying $\imath_L (x_1) = \imath_L (x_2)$. If $L_1$ is tangent to $L_2$ at $\imath_L (x_1) = \imath_L (x_2)$, we may further ensure that $\alpha$ viewed as a section of $TU_{L_2}$ is vanishing on $L_{12}$. Consequently, the affine structures on the Lagrangian fibre $e_L (U_{x_1}) = e_L (U_{x_2})$ determined through $U_{x_1}$ and $U_{x_2}$ coincide.
\hfill$\Box$

\begin{re}
\label{bi}
{\bf [\textit{non-Hausdorff manifold}]:} Let $\tilde{\mathscr L}$ be a compact non-Hausdorff $n$-manifold, by identifying any 2 points in $\tilde{\mathscr L}$ that is inseparable, we get a separable topological space ${\mathscr L}$ and the natural projection $\tilde{\pi}: \tilde{\mathscr L} \rightarrow {\mathscr L}$. ${\mathscr L}$ still carry a system of (not necessarily open) charts inherited from $\tilde{\mathscr L}$. ${\mathscr L}$ is very much like a manifold, tangent bundle and submanifolds etc. all make perfect sense. Use the system of charts, the original non-Hausdorff manifold $\tilde{\mathscr L}$ can be recovered from ${\mathscr L}$. To avoid confusion, we will call such ${\mathscr L}$ together with the charts a NH manifold. We should note here that the system of charts on ${\mathscr L}$ that is compatible with the tangent bundle of ${\mathscr L}$ is not unique and give rise to different non-Hausdorff manifolds $\tilde{\mathscr L}$. In our discussion, ${\mathscr L}$ is a smooth Lagrangian NH skeleton of $(M, \omega)$, which is Hamiltonian equivalent to a usual Lagrangian skeleton. We will also consider $n$-manifold $L_i$ together with a map $\imath_{L_i}: L_i \rightarrow {\mathscr L}$, where $\imath_{L_i} = \tilde{\pi} \circ \tilde{\imath}_{L_i}$ and $\tilde{\imath}_{L_i}: L_i \rightarrow \tilde{\mathscr L}$ is an immersion. One can observe that only charts of ${\mathscr L}$ that come from charts of some $L_i$ are relevant in our discussion. In such situation, we may also consider $\imath_{L_i}$ as a Lagrangian immersion from $L_i$ to $(M, \omega)$. We should point out here that for all practical purpose, it is possible (although quite inconvenient) to avoid mentioning the NH manifold ${\mathscr L}$ and only refer to $\imath_{L_i}: L_i \rightarrow M$ in all our discussion.
\end{re}

Let ${\mathscr L}$ be a compact NH Lagrangian submanifold in a symplectic manifold $(M, \omega)$. Let $L_i$ be a $n$-manifold together with an immersion $\imath_{L_i}: L_i \rightarrow {\mathscr L}$. Due to the structure of ${\mathscr L}$, when $L_i$ self-intersect, the tangent spaces always coincide. Apply propositions \ref{mf} and \ref{mh} to $N$ being the normal bundle of $L_i$ with respect to the metric, we have the open symplectic immersion $e_{L_i}: (U_{L_i}, \omega_{L_i}) \rightarrow (M, \omega)$ together with the Lagrangian fibration $\pi_i = \pi_{L_i}: U_{L_i} \rightarrow L_i$. Define the Lagrangian immersion $\imath_{\epsilon f_i} = e_{L_i} (\epsilon df_i): L_i \rightarrow (M, \omega)$, where $f_i$ is a function on $L_i$. Let $\Lambda^\epsilon_{f_i} = \Lambda_{\epsilon f_i} := \imath_{\epsilon f_i} (L_i)$.\\

For $i\not=j$, it is reasonable to assume that the fibre product $L_{ij} := L_i \times_{\mathscr L} L_j$ is a union of closed regions with smooth boundaries. $f_{ij} (x_1,x_2) := f_j(x_2) - f_i(x_1)$ for $(x_1,x_2) \in L_{ij} = L_i \times_{\mathscr L} L_j$ defines a function $f_{ij}$ on $L_{ij}$, which we will assume to be a Morse function with only critical points in the interior of $L_{ij}$.\\

Via projections of $L_{ij}$ to $L_i$ and $L_j$, locally, $L_{ij}$ can be viewed as subset in $L_i$ and $L_j$. Via $L_i$ (resp. $L_j$), $L_{ij}$ can be extended to an open manifold $\hat{L}_i$ (resp. $\hat{L}_j$) as a small open neighborhood of $L_{ij}$, so that $L_{ij} = \hat{L}_i \cap \hat{L}_j$, and $\hat{L}_{ij} = \hat{L}_i \cup \hat{L}_j$ can be viewed as a smooth NH manifold with Lagrangian immersion $\hat{\imath}_{\hat{L}_i}: \hat{L}_i \rightarrow L_i \rightarrow (M, \omega)$ (resp. $\hat{\imath}_{\hat{L}_j}: \hat{L}_j \rightarrow L_j \rightarrow (M, \omega)$).\\

Locally near a point in $\partial L_{ij}$, we may identify $L_i$ (resp. $L_j$) and $U_{L_i}$ (resp. $U_{L_j}$) with their embedded images in $M$ and view $L_{ij}$ as a closed region in $L_i$ (resp. $L_j$). Let $\hat{L}_i$ (resp. $\hat{L}_j$) be a small neighborhood of $L_{ij}$ in $L_i$ (resp. $L_j$). Then $L_j$ coincides with $\Lambda_{h_{ij}}$ in $U_{\hat{L}_i}$ for a smooth function $h_{ij}$ defined on $\hat{L}_i$ such that $L_{ij} = \{h_{ij}=0\}$. Similarly, there is a smooth function $\hat{f}_j$ on $\hat{L}_i$ such that $\Lambda_{\epsilon f_j}$ coincides with $\Lambda_{h_{ij} + \epsilon \hat{f}_j}$ in $U_{\hat{L}_i}$. By proposition \ref{mh}, one can ensure $\hat{f}_j = f_j$ and $\hat{f}_{ij} = f_{ij}$ on $L_{ij} \subset \hat{L}_i$, where $\hat{f}_{ij} := \hat{f}_j -f_i$.\\

We will assume that $\nabla f_{ij}$ is pointing outward (resp. inward) on $\partial L_{ij}$ if $h_{ij}\geq 0$ (resp. $h_{ij}\leq 0$) near $\partial L_{ij}$. Under our previous assumption that $df_{ij} \not=0$ on $\partial L_{ij}$, this condition is equivalent to $(\nabla \hat{f}_{ij}, \nabla h_{ij}) >0$ outside of $L_{ij}$. This condition will prevent any intersection of $\Lambda^\epsilon_{f_i}$ and $\Lambda^\epsilon_{f_j}$ near $\partial L_{ij}$. Consequently, under our assumptions, $\textsc{x}^\epsilon \in \Lambda^\epsilon_{f_i} \cap \Lambda^\epsilon_{f_j}$ if and only if $\textsc{x} = \pi_i(\textsc{x}^\epsilon) = \pi_j(\textsc{x}^\epsilon)$ is in the interior of $L_{ij}$ and is a critical point of $f_{ij}$.

\begin{re}
Since immersion is locally embedding. When we perform local discussion, we may regard $\imath_{L_i}$ (resp. $e_{L_i}$) as embedding, then it is convenient to identify $L_i$ (resp. $U_{L_i}$) with its embedding image in $M$. We will keep such practice through out the paper without mentioning it each time, whenever we are performing such local discussion.
\end{re}

\begin{re}
{\bf [\textit{embedded versus immersed Lagrangian submanifolds}]:}\\
The usual (embedded) Lagrangian submanifold in a symplectic manifold $(M,\omega)$ can be viewed as the image of a Lagrangian embedding from a manifold $L$ to $(M, \omega)$. We can generalize this concept to immersed Lagrangian submanifold as the image of a generic Lagrangian immersion from a manifold $L$ to $(M, \omega)$. By generic immersion, we mean that the image of the manifold only self intersects transversely at isolated points.\\

The simplest generic condition on individual $\Lambda^\epsilon_{f_i}$ is to require it to be an embedded Lagrangian submanifold. This condition is not as restrictive as it might sound. In usual Fukaya category, one only consider embedded Lagrangian submanifold as objects. Under the flow discussed in the introduction, an embedded Lagrangian submanifold will remain embedded before it is collapsed into the Lagrangian skeleton at the limit. If we are considering the usual Fukaya category, $L_i$ is meant to be such collapsed remain of an embedded Lagrangian submanifold under the flow, and should be Hamiltonian equivalent to an embedded Lagrangian submanifold. Consequently, it would be quite reasonable to assume $\Lambda^\epsilon_{f_i}$ to be an embedded Lagrangian submanifold.\\

On the other hand, it make sense to generalize Fukaya category to include immersed Lagrangian submanifolds as objects (see remark \ref{bh}). Our computation actually work in this more general context. Namely, we only need to assume that $\Lambda^\epsilon_{f_i}$ is an immersed Lagrangian submanifold. Such added flexibility should be beneficial even for computation in the usual Fukaya category.\\
\end{re}

{\samepage
\stepcounter{subsection}
{\bf \S \thesubsection\ $A^\infty$ products in Fukaya category}\\
\nopagebreak

Let $\mathbb{D}_{[\scriptsize{\textsc{n}}]}$} denote the unit disk $\mathbb{D} = \{t\in \mathbb{C}: |t|\leq 1\}$ with $\textsc{n}$ marked points $\{\textsc{p}_1, \cdots, \textsc{p}_{\scriptsize{\textsc{n}}}\}$ located in $S^1 \cong \partial \mathbb{D}$ according to the counter clockwise cyclic ordering. There is the natural decomposition $\partial \mathbb{D}_{[\scriptsize{\textsc{n}}]} := \partial \mathbb{D} \setminus \{\textsc{p}_1, \cdots, \textsc{p}_{\scriptsize{\textsc{n}}}\} = \partial_1 \mathbb{D}_{[\scriptsize{\textsc{n}}]} \cup \cdots \cup \partial_{\scriptsize{\textsc{n}}} \mathbb{D}_{[\scriptsize{\textsc{n}}]}$, where $\partial_i \mathbb{D}_{[\scriptsize{\textsc{n}}]}$ is between $\textsc{p}_{i-1}$ and $\textsc{p}_i$ with index $i$ viewed as an integer modulo $\textsc{n}$. Let ${\cal I}_{\scriptsize{\textsc{n}}}$ denote the moduli space of such $\mathbb{D}_{[\scriptsize{\textsc{n}}]}$.\\

For a compatible almost complex structure $J$ on $(M, \omega)$, functions $\vec{f} = (f_1, \cdots, f_{\scriptsize{\textsc{n}}})$ and Lagrangian submanifolds $\vec{\Lambda}^\epsilon = \{\Lambda^\epsilon_{f_i}\}$ with Lagrangian intersection points $\vec{\textsc{x}}^\epsilon = (\textsc{x}^\epsilon_1, \cdots, \textsc{x}^\epsilon_{\scriptsize{\textsc{n}}})$, where $\textsc{x}^\epsilon_i \in \Lambda^\epsilon_{f_i} \cap \Lambda^\epsilon_{f_{i+1}}$, we may define the moduli space ${\cal M}_J(M, \vec{\Lambda}^\epsilon, \vec{\textsc{x}}^\epsilon)$ that consists of pairs $\mathbb{D}_{[\scriptsize{\textsc{n}}]} = (\mathbb{D}, \{\textsc{p}_i\}) \in {\cal I}_{\scriptsize{\textsc{n}}}$ and pseudo holomorphic map $w: (\mathbb{D}, \partial_i \mathbb{D}_{[\scriptsize{\textsc{n}}]}, \textsc{p}_i) \rightarrow (M, \Lambda^\epsilon_{f_i}, \textsc{x}^\epsilon_i)$.

\begin{re}
\label{bh}
If $\Lambda^\epsilon_{f_i}$ is an immersed Lagrangian submanifold, we need to require that $w(\partial_i \mathbb{D}_{[\scriptsize{\textsc{n}}]})$ to be a smooth curve in $\Lambda^\epsilon_{f_i}$. Namely, $w(\partial_i \mathbb{D}_{[\scriptsize{\textsc{n}}]})$ does not connect 2 local components of $\Lambda^\epsilon_{f_i}$ through a self-intersection point.
\end{re}

Clearly, $\textsc{x}_i = \pi_i(\textsc{x}^\epsilon_i) \in L_{i,i+1}$ is a non-degenerate critical point of $f_{i,i+1}$, where $\pi_i: U_{L_i} \rightarrow L_i$ is the natural projection. Let $\mu(\textsc{x}_i)$ be the corresponding Morse index. Any version of the definition of the Fukaya category necessarily should have shown that under the corresponding norm, the virtual dimension of ${\cal M}_J(M, \vec{\Lambda}^\epsilon, \vec{\textsc{x}}^\epsilon)$ is

\[
{\rm Ind} (\vec{f}, \vec{\textsc{x}}) = \sum_{i=1}^{\scriptsize{\textsc{n}}} \mu(\textsc{x}_i) - (\textsc{n}-1)n + (\textsc{n}-3).
\]

To define the $A^\infty$ products, it is convenient to assume $\vec{f} = (f_0, \cdots, f_{\scriptsize{\textsc{n}}})$. Then

\[
{\rm Ind} (\vec{f}, \vec{\textsc{x}}) = \mu(\textsc{x}_{\scriptsize{\textsc{n}}}) + \sum_{i=0}^{{\scriptsize{\textsc{n}}}-1} (\mu(\textsc{x}_i) - n) + (\textsc{n}-2).
\]

Define ${\rm Hom} (\Lambda^\epsilon_{f_i}, \Lambda^\epsilon_{f_j}) = {\rm CF}^* (\Lambda^\epsilon_{f_i}, \Lambda^\epsilon_{f_j}; \mathbb{R}) = {\rm Span}_{\mathbb{R}} (\Lambda^\epsilon_{f_i} \cap \Lambda^\epsilon_{f_j})$, then the $A^\infty$ product $m_{\scriptsize{\textsc{n}}}: {\rm Hom} (\Lambda^\epsilon_{f_0}, \Lambda^\epsilon_{f_1}) \otimes \cdots \otimes {\rm Hom} (\Lambda^\epsilon_{f_{{\scriptsize{\textsc{n}}}-1}}, \Lambda^\epsilon_{f_{\scriptsize{\textsc{n}}}}) \rightarrow {\rm Hom} (\Lambda^\epsilon_{f_0}, \Lambda^\epsilon_{f_{\scriptsize{\textsc{n}}}}) [2-\textsc{n}]$ every version of the definition of the Fukaya category can agree on is

\begin{equation}
\label{mi}
m^{\scriptsize{\textsc{n}}}_{\rm Fukaya} ([\textsc{x}^\epsilon_0], \cdots, [\textsc{x}^\epsilon_{{\scriptsize{\textsc{n}}}-1}]) = \sum_{\tiny{\begin{array}{c}\textsc{x}^\epsilon_{\tiny{\textsc{n}}} \in \Lambda^\epsilon_{f_0} \cap \Lambda^\epsilon_{f_{\tiny{\textsc{n}}}}\\ {\rm Ind} (\vec{f}, \vec{\textsc{x}}) = 0 \end{array}}} \# {\cal M}_J(M, \vec{\Lambda}^\epsilon, \vec{\textsc{x}}^\epsilon) [\textsc{x}^\epsilon_{\scriptsize{\textsc{n}}}].
\end{equation}

\begin{re}
\label{bg}
Notice that ${\cal M}_J(M, \vec{\Lambda}^\epsilon, \vec{\textsc{x}}^\epsilon)$ is the same under whatever norm one uses to discuss it. Therefore, to actually compute the product, one may use whatever norm one please, as long as one can show that for generic $f_i$ and small $\epsilon$, when ${\rm Ind} (\vec{f}, \vec{\textsc{x}}) = 0$, the linearized operator at $w \in {\cal M}_J(M, \vec{\Lambda}^\epsilon, \vec{\textsc{x}}^\epsilon)$ has zero kernel, which implies $\dim {\cal M}_J(M, \vec{\Lambda}^\epsilon, \vec{\textsc{x}}^\epsilon) =0$ in particular.\\
\end{re}

{\samepage
\stepcounter{subsection}
{\bf \S \thesubsection\ $A^\infty$ products in Morse category}\\
\nopagebreak

Recall that ${\mathscr L}$} is a smooth n-dimensional NH manifold. An object of the Morse category of ${\mathscr L}$ is a pair $(L_i, f_i)$ together with an open immersion $\imath_{L_i}: L_i \rightarrow {\mathscr L}$, where $L_i$ is a n-dimensional manifold and $f_i$ is a smooth function on $L_i$. The genericity assumptions for $(L_i, f_i)$ are stated in page \pageref{mp}. (For brevity, we will simply refer to the object $(L_i, f_i)$ by $f_i$ in this section.)\\

A graph $\Gamma$ has a set $C^1(\Gamma)$ of legs and a set $C^0(\Gamma)$ of vertices. 1-valent vertices form the set of external vertices $C^0_{\rm ext}(\Gamma)$. The rest in $C^0(\Gamma)$ form the set of internal vertices $C^0_{\rm int}(\Gamma)$. Legs reaching 1-valent vertices form the set of external legs $C^1_{\rm ext}(\Gamma)$. The rest in $C^1(\Gamma)$ form the set of internal legs $C^1_{\rm int}(\Gamma)$.\\

A graph $\Gamma$ can be viewed as a metric graph or a topological one. As a metric graph, $\Gamma$ comes with leg length $\gamma_\nu >0$ for each $\nu\in C^1(\Gamma)$. We assume that the metric graph is complete, namely, $\gamma_\nu = +\infty$ if and only if $\nu \in C^1_{\rm ext}(\Gamma)$. The topology of $\Gamma$ is captured by the set of strata $C^*(\Gamma) := C^0_{\rm int}(\Gamma) \cup C^1(\Gamma)$ together with the inclusion relations. A topological map from $\Gamma_1$ to $\Gamma_2$ is a map $C^*(\Gamma_1) \rightarrow C^*(\Gamma_2)$ that preserve the inclusion relations. In particular, $C^0_{\rm int}(\Gamma_1) \rightarrow C^0_{\rm int}(\Gamma_2)$. \\

We are interested in tree graph with cyclically ordered external vertices $C^0_{\rm ext}(\Gamma) = \{\textsc{p}_1, \cdots, \textsc{p}_{\scriptsize{\textsc{n}}}\}$ and no 2-valent vertices, which is called a ribbon tree in \cite{FO}. It is convenient to denote $C^1_{\rm ext}(\Gamma) = \{\nu_1, \cdots, \nu_{\scriptsize{\textsc{n}}}\}$, where $\nu_i$ leads to $\textsc{p}_i$. In \cite{FO}, $G_{\scriptsize{\textsc{n}}}$ (resp. $Gr_{\scriptsize{\textsc{n}}}$) denotes the moduli space of ribbon trees $\Gamma$ with $C^0_{\rm ext}(\Gamma) = \{\textsc{p}_1, \cdots, \textsc{p}_{\scriptsize{\textsc{n}}}\}$ as topological graphs (resp. complete metric graphs). For any ribbon tree $\Gamma$, there is a unique embedding $i_\Gamma: \Gamma \rightarrow \mathbb{D}$ up to isotopy so that only $C^0_{\rm ext}(\Gamma) = \{\textsc{p}_1, \cdots, \textsc{p}_{\scriptsize{\textsc{n}}}\}$ is mapped to $\partial \mathbb{D}$ with counter clockwise order. $i_\Gamma (\Gamma)$ separates $\mathbb{D}$ into $k$ regions $\mathbb{D}_1, \cdots, \mathbb{D}_{\scriptsize{\textsc{n}}}$, where $\mathbb{D}_i$ contains $i_\Gamma(\textsc{p}_{i-1})$ and $i_\Gamma(\textsc{p}_i)$. For $\nu \in C^* (\Gamma)$, let $I(\nu) = \{i: i_\Gamma (\nu) \subset \mathbb{D}_i\}$.\\

A gradient tree is a map $\Upsilon = \Upsilon_\Gamma: \Gamma \rightarrow {\mathscr L}$ such that $\textsc{x}_i = \Upsilon (\textsc{p}_i)$ is a critical point of $f_{i,i+1}$, and for $\nu \in C^1 (\Gamma)$ with a direction, $\Upsilon (\nu)$ is a gradient flow line determined by $\nabla f_{ij}$, where $\mathbb{D}_i$ (resp. $\mathbb{D}_j$) is in the right (resp. left) side of $i_\Gamma (\nu)$. In particular, for $\nu_i \in C^1_{\rm ext}(\Gamma)$ pointing toward $\textsc{p}_i$, $\Upsilon (\nu_i)$ is a gradient flow line determined by $\nabla f_{i,i+1}$. With slight abuse of notation, $\Upsilon$ is also used to denote the image $\Upsilon (\Gamma)$ of $\Upsilon$. $\Upsilon (\nu)$ for $\nu \in C^0 (\Gamma)$ (resp. $\nu \in C^1 (\Gamma)$) is often referred to as a vertex (resp. a leg) of $\Upsilon$ (resp. if $\dim \Upsilon (\nu) =1$).\\

For $\vec{f} = (f_1, \cdots, f_{\scriptsize{\textsc{n}}})$ and $\vec{\textsc{x}} = (\textsc{x}_1, \cdots, \textsc{x}_{\scriptsize{\textsc{n}}})$, where $\textsc{x}_i$ is a critical point of $f_{i,i+1}$, we may define the moduli space ${\cal M}_g ({\mathscr L}, \vec{f}, \vec{\textsc{x}})$ that consists of pairs $\Gamma \in Gr_{\scriptsize{\textsc{n}}}$ and gradient tree $\Upsilon: \Gamma \rightarrow {\mathscr L}$ defined according to $(\vec{f}, \vec{\textsc{x}})$.

\begin{re}
{\bf [\textit{exceptional gradient tree}]:} There is one special type of gradient trees that we call exceptional, was not paid explicit attention in \cite{FO}, that turns out to be the only gradient trees encountered in our application (\cite{BR}) to the homological mirror symmetry for weighted projective spaces.\\

$\textsc{x}_i = \Upsilon (\textsc{p}_i)$ is called an {\sf exceptional vertex} of $\Upsilon$ if $\textsc{x}_i$ is a local minimal of $f_{i,i+1}$, or $\mu (\textsc{x}_i)=0$. Then $\nabla f_{i,i+1}$ is pointing outward, and for $\nu_i \in C^1_{\rm ext}(\Gamma)$ pointing toward $\textsc{p}_i$, $\Upsilon (\nu_i)$ has to be a single point $\textsc{x}_i$. Generically, $\textsc{x}_i$ is a smooth point of $\Upsilon (\Gamma)$. (In figure 1, $\Upsilon_2$ is an exceptional gradient tree with the exceptional vertex $\textsc{x}_3$.) In such situation, it is more appropriate to consider a smaller {\sf reduced} graph $\hat{\Gamma}$ (through removing $\nu$ and related vertices) together with injective map $\hat{\Upsilon}: \hat{\Gamma} \rightarrow {\mathscr L}$ and surjective map $\hat{\pi}: \Gamma \rightarrow \hat{\Gamma}$ such that $\Upsilon = \hat{\Upsilon} \circ \hat{\pi}$. $\Gamma = \hat{\Gamma}$ if and only if $\mu(\textsc{x}_i) \not=0$ for all $i$, or in another word, $\Upsilon$ has no exceptional vertex. We will call $\Upsilon$ {\sf exceptional} if $\Gamma \not= \hat{\Gamma}$. On each directed leg of $\hat{\Gamma}$, the gradient flow determines a natural coordinate up to shift. In another word, $\Upsilon$ determines a complete metric graph structure on the reduced graph $\hat{\Gamma}$. (In such context, with slight abuse of notation, we will also call $\textsc{p}_i$ an exceptional vertex of $\Gamma$ and $\nu_i \in C^1_{\rm ext} (\Gamma)$ an exceptional leg of $\Gamma$.)
\end{re}
\pspicture(-6,-2.5)(6,2)
\psset{linewidth=.5pt}

\rput(-2.5,0){
\pscustom{\pscurve{-<}(-2,1.5)(-2,1.5)(-1.5,.8)\pscurve(-1.5,.8)(-1,0)}
\pscustom{\pscurve{-<}(-2,-1.5)(-2,-1.5)(-1.5,-.8)\pscurve(-1.5,-.8)(-1,0)}
\pscustom{\pscurve{->}(-1,0)(-1,0)(0,.05)\pscurve(0,.05)(1,0)}
\pscustom{\pscurve{-<}(2,1.5)(2,1.5)(1.5,.8)\pscurve(1.5,.8)(1,0)}
\pscustom{\pscurve{-<}(2,-1.5)(2,-1.5)(1.5,-.8)\pscurve(1.5,-.8)(1,0)}

\rput(0,.3){$\nabla f_{24}$} \rput(-1.9,.6){$\nabla f_{12}$} \rput(-1.9,-.6){$\nabla f_{41}$} \rput(1.9,.5){$\nabla f_{23}$} \rput(1.9,-.6){$\nabla f_{34}$}
\rput(-2.3,1.4){$\textsc{x}_1$} \rput(-2.3,-1.5){$\textsc{x}_4$} \rput(2.3,1.4){$\textsc{x}_2$} \rput(2.3,-1.5){$\textsc{x}_3$}
\rput(0,-1){$\Upsilon_1$} \rput(0,-1.4){(non-exceptional)}
\qdisk(-2,1.5){1.5pt} \qdisk(-2,-1.5){1.5pt} \qdisk(2,1.5){1.5pt} \qdisk(2,-1.5){1.5pt}}

\rput(3,0){
\pscustom{\pscurve{-<}(-1.5,-1.4)(-1.5,-1.4)(-.8,-.7) \pscurve(-.8,-.7)(0,0) \pscurve{->}(0,0)(.8,.62) \pscurve(.8,.62)(1.5,1.1)}

\rput(-.5,-1){$\nabla f_{12}$} \rput(1.1,.3){$\nabla f_{23}$}
\rput(-1.8,-1.4){$\textsc{x}_1$} \rput(1.9,1.2){$\textsc{x}_2$} \rput(-.2,.2){$\textsc{x}_3$}
\rput(1,-1){$\Upsilon_2$} \rput(1,-1.4){(exceptional)}
\qdisk(-1.5,-1.4){1.5pt} \qdisk(0,0){1.5pt} \qdisk(1.5,1.1){1.5pt}}

\stepcounter{figure}
\uput{2}[d](0,0){Figure \thefigure: Examples of gradient trees}
\endpspicture

For a gradient tree $\Upsilon: \Gamma \rightarrow {\mathscr L}$, since tree $\Gamma$ is contractible, there exist an open manifold $L_\Upsilon$ with open immersion $\imath_{L_\Upsilon}: L_\Upsilon \rightarrow {\mathscr L}$ and a lift $\tilde{\Upsilon}: \Gamma \rightarrow L_\Upsilon$ such that $\Upsilon = \imath_{L_\Upsilon} \circ \tilde{\Upsilon}$. (We want to think of $L_\Upsilon$ as a tubular neighborhood of $\tilde{\Upsilon} (\Gamma)$.) $f_i$ defined on $L_i$ can also be lifted to $\tilde{f_i}$ defined on part of $L_\Upsilon$ that can be suitably extended to all of $L_\Upsilon$ so that $\tilde{f}_{ij} = \tilde{f}_j - \tilde{f}_i$ is a Morse function on $L_\Upsilon$ for $i\not=j$. $\tilde{\Upsilon}: \Gamma \rightarrow L_\Upsilon$ is a gradient tree in the exact sense of \cite{FO} with respect to the functions $\tilde{f}_i$ on $L_\Upsilon$. For this reason, results in \cite{FO} concerning gradient trees can essentially be carried over to our case directly. In particular, according to theorem 1.4 in \cite{FO}, the moduli space ${\cal M}_g ({\mathscr L}, \vec{f}, \vec{\textsc{x}})$ is a $C^\infty$ manifold of dimension ${\rm Ind} (\vec{f}, \vec{\textsc{x}})$. We are now in a position to summarize all our genericity assumptions (that we assume through out the paper) concerning $(L_i, f_i)$ in the following:\\

{\bf Genericity assumptions:}\label{mp} $\Lambda^\epsilon_{f_i}$ is an immersed Lagrangian submanifold when $\epsilon>0$. For $i\not=j$, $\partial L_{ij}$ is smooth, $f_{ij}$ is a Morse function on $L_{ij}$ with only critical points in the interior of $L_{ij}$ and $\nabla f_{ij}$ is pointing outward (resp. inward) on $\partial L_{ij}$ if $h_{ij}\geq 0$ (resp. $h_{ij}\leq 0$) near $\partial L_{ij}$. None of the critical points of $f_{ij}$ for all $i,j$ coincide. There is no gradient tree with ${\rm Ind} (\vec{f}, \vec{\textsc{x}}) <0$ and rigid gradient tree (${\rm Ind} (\vec{f}, \vec{\textsc{x}})=0$) is isolated with only 3-valent interior vertices.\\

With the immersion $\imath_{L_\Upsilon}: L_\Upsilon \rightarrow {\mathscr L} \rightarrow M$, we can also get an open manifold $M_\Upsilon$ from $M$ as a tubular neighborhood of $\tilde{\Upsilon} (\Gamma)$ and $L_\Upsilon$. In $M_\Upsilon$, it is convenient to piece together $\pi_i: U_{L_i} \rightarrow L_i$ to form the Lagrangian fibration (projection) $\label{bj}\pi_\Upsilon: M_\Upsilon \rightarrow L_\Upsilon$. More precisely, $\pi_\Upsilon := \pi_i$ near $\tilde{\Upsilon} (\nu)$ for $\nu \in C^*(\Gamma)$ and $i \in I(\nu)$. Proposition \ref{mh} implies that the definitions coincide in the common regions.\\

To define the $A^\infty$ products, it is convenient to assume $\vec{f} = (f_0, \cdots, f_{\scriptsize{\textsc{n}}})$. Define ${\rm Hom} (f_i, f_j) = {\rm C}^* (f_i, f_j; \mathbb{R}) = {\rm Span}_{\mathbb{R}} ({\rm Crit} (f_{ij}))$. We may define the $A^\infty$ product $m_{\scriptsize{\textsc{n}}}: {\rm Hom} (f_0, f_1) \otimes \cdots \otimes {\rm Hom} (f_{{\scriptsize{\textsc{n}}}-1}, f_{\scriptsize{\textsc{n}}}) \rightarrow {\rm Hom} (f_0, f_{\scriptsize{\textsc{n}}}) [2-\textsc{n}]$

\begin{equation}
\label{mj}
m^{\scriptsize{\textsc{n}}}_{\rm Morse} ([\textsc{x}_0], \cdots, [\textsc{x}_{{\scriptsize{\textsc{n}}}-1}]) = \sum_{\tiny{\begin{array}{c}\textsc{x}_{\tiny{\textsc{n}}} \in {\rm Crit} (f_{0\tiny{\textsc{n}}})\\ {\rm Ind} (\vec{f}, \vec{\textsc{x}}) = 0 \end{array}}} \# {\cal M}_g ({\mathscr L}, \vec{f}, \vec{\textsc{x}}) [\textsc{x}_{\scriptsize{\textsc{n}}}].
\end{equation}

\begin{theorem}
\label{mk}
For $\vec{f} = (f_0, \cdots, f_{\scriptsize{\textsc{n}}})$ satisfying genericity conditions (p.~\pageref{mp}) and $\epsilon>0$ small enough, we have $m^{\scriptsize{\textsc{n}}}_{\rm Fukaya} = m^{\scriptsize{\textsc{n}}}_{\rm Morse}$, under the natural identifications ${\rm Hom}_{\rm Fukaya} (\Lambda^\epsilon_{f_i}, \Lambda^\epsilon_{f_j}) \cong {\rm Hom}_{\rm Morse} (f_i, f_j)$.
\end{theorem}
{\bf Proof:} With considerations in remark \ref{bg}, compare (\ref{mi}) and (\ref{mj}), it is clear that one only need to show that ${\cal M}_J(M, \vec{\Lambda}^\epsilon, \vec{\textsc{x}}^\epsilon) \cong {\cal M}_g ({\mathscr L}, \vec{f}, \vec{\textsc{x}})$ with signs preserved, when ${\rm Ind} (\vec{f}, \vec{\textsc{x}})=0$. The identification of the 2 moduli spaces for $\epsilon$ small is a consequence of theorems \ref{pf} and \ref{qa}.
\hfill$\Box$

\begin{re}
{\bf [\textit{Orientation}]:} Orientations of the moduli spaces are not addressed in this paper. (Coincidentally, orientation was also not addressed in \cite{FO}.) Therefore the identification we proved in this theorem should be understood as over $\mathbb{Z}_2$ for the moment. The orientation of ${\cal M}_g ({\mathscr L}, \vec{f}, \vec{\textsc{x}})$ is a rather straightforward generalization of the case of classical Morse theory. The orientation of ${\cal M}_J(M, \vec{\Lambda}^\epsilon, \vec{\textsc{x}}^\epsilon)$ and the identification of the orientations will be discussed in \cite{disk}.
\end{re}

\begin{re}
{\bf [\textit{Generic conditions}]:} Genericity perturbation arguments for gradient trees are much easier than genericity perturbation arguments for pseudo holomorphic polygons, because one is achieving genericity for finite dimension submanifolds (stable, unstable submanifolds, etc.) through infinite dimension perturbations (of $f_{ij}$ near its critical points). Most of such difficulty was already encountered in classical Morse theory. One important consequence of \cite{FO} and subsequently this paper is that not only pseudo holomorphic polygons can be counted through counting gradient trees, but also the much more difficult genericity perturbation arguments for pseudo holomorphic polygons can in a way be reduced to much easier genericity perturbation arguments for gradient trees. One reflection is that instead of perturbing almost complex structure, one perturb the Morse functions $f_{ij}$.
\end{re}

\se{Cheeger-Gromov convergence}
Cheeger-Gromov convergence is a very convenient tool for discussing limits concerning metric spaces, Riemannian (\k) manifolds with possibly various tensor on the manifolds or submanifolds, etc. One of the great strength of the Cheeger-Gromov convergence is its general applicability. One start with Gromov-Hausdorff convergence to get the limit, which applies to general compact metric spaces. (The concept can be generalized to pointed Gromov-Hausdorff convergence that applies to complete metric space with marked point.) When the objects are manifolds, the convergence can usually be strengthened to $C^l$-convergence over the smooth part of the limit in the sense of Cheeger-Gromov. The existence of the Gromov-Hausdorff limit is usually achieved after possibly passing to subsequence through certain compactness theorems (for example, Gromov compactness theorem) based on Ascoli type arguments.\\
\begin{re}
\label{bk}
The Cheeger-Gromov convergence results we use in this paper are among the most classical types (that were developed by Cheeger, Fukaya, Gromov, etc. in the Riemannian case) requiring the uniform bound of local geometry. Such Cheeger-Gromov convergence results enable us to consider the limit of a sequence of submanifolds or maps as if the domain and the target space are fixed, therefore reduces to the classical consideration of limit of functions. Such simple yet powerful application was indeed the original motivation of the Cheeger-Gromov convergence.
\end{re}

Consider an almost K\"{a}hler structure $(M,g,J,\omega)$, where the triple $(g,J,\omega)$ of Riemannian metric, almost complex structure and symplectic form on $M$ are mutually compatible. $(M,g,J,\omega)$ is said to have uniformly bounded local geometry if: 1. The full curvature of $g$ and all its covariant multi-derivatives are uniformly bounded. 2. The injective radius of $g$ is uniformly bounded from below. 3. $J$ and all its covariant multi-derivatives are uniformly bounded. This condition is equivalent to the fact that there exists $r_0>0$ so that $(M,g,J,\omega)$ can be covered by coordinate charts $U \rightarrow B(r_0) \subset (\mathbb{C}^n, g_0)$ such that $g$ and $g_0$ are uniformly quasi-isometric on $U$, $g$ and $J$ and all their multi-derivatives with respect to the coordinate on $U$ are uniformly bounded. We will call such charts preferred charts. It is straightforward to check that the transition functions among preferred charts have uniformly bounded multi-derivatives.\\

A submanifold $L \subset (M,g,J,\omega)$ is of bounded local geometry if $(M,g,J,\omega)$ is of bounded local geometry and for $r_0>0$ small enough $L$ in each corresponding preferred chart can be expressed as a graph with uniformly bounded multi-derivatives. The generalization of this concept to non-Hausdorff submanifolds is rather straightforward.\\

When $M$ has smooth boundary $\partial M$, $(M, \partial M,g,J,\omega)$ is said to have uniformly bounded local geometry if in addition, neighborhood of $\partial M$ can be covered by coordinate charts $U \rightarrow B_+(r_0) \subset (\mathbb{C}_+^n, g_0)$ such that $g$ and $g_0$ are uniformly quasi-isometric on $U$, $g$ and $J$ and all their multi-derivatives with respect to the coordinate on $U$ are uniformly bounded, where $B_+(r_0) = B(r_0) \cap \mathbb{C}_+^n$ and $\mathbb{C}_+^n = \{z \in \mathbb{C}^n: z_n^\Re \geq0\}$. The generalization of this concept to manifolds with corners (boundary with real normal crossing singularity) is rather straightforward.\\

$M$ is called a complex if $M$ is stratified with each strata a manifold with corners and a local neighborhood near singularity of $M$ can be identified with a local neighborhood of a simplicial complex. A complex $M$ is said to have uniformly bounded local geometry if each of its strata is a manifold with corners that is of uniformly bounded local geometry and the number of strata near each point of $M$ is uniformly bounded. We are only going to deal with graph (proposition \ref{ac}), which is the simplest complex. In the graph case, the condition of uniformly bounded local geometry reduces to that the number of legs intersecting a ball of fixed radius $r_0>0$  is uniformly bounded.\\

In the following, we state the convergence results that we need. The case of Riemannian manifold is standard that goes back to Cheeger, Fukaya, Gromov, etc. The case of K\"{a}hler manifold have been discussed by Tian, the author, etc. in various situations. (Please consult \cite{Ruan1} for the Riemannian convergence result, its precise references and some K\"{a}hler convergence results.) The proofs of all these cases are straightforward generalization of the Riemannian case without real additional technical difficulties. But the proofs get increasingly more tedious and notations get increasingly more complicated for the more and more complex cases. Such proofs are best left as exercises to the readers.\\
\begin{prop}
\label{ba}
Let $(M_k,g_k,J_k,\omega_k)$ be a sequence of complete K\"{a}hler manifolds with uniformly bounded local geometry and marked points $z_k \in M_k$. By possibly taking subsequence, $(M_k,z_k, g_k,J_k,\omega_k)$ pointed converge in the sense of Cheeger-Gromov to a complete K\"{a}hler manifold $(M_o,z_o, g_o,J_o,\omega_o)$ with bounded local geometry and the marked point $z_o \in M_o$.
\hfill$\Box$
\end{prop}

\begin{co}
\label{be}
Let $(M,g,J,\omega)$ be a complete K\"{a}hler manifold with uniformly bounded local geometry. For a sequence of points $p_k \in M$, by possibly taking subsequence, $(M,p_k, g,J,\omega)$ pointed converge in the sense of Cheeger-Gromov to a complete K\"{a}hler manifold $(M_o,p_o, g_o,J_o,\omega_o)$ with bounded local geometry and the marked point $p_o \in M_o$.
\hfill$\Box$
\end{co}

\begin{co}
\label{bd}
Let $(M,g,J,\omega)$ be a complete K\"{a}hler manifold. For $p = \lim p_k \in M$ and $\lim \epsilon_k = 0$, $(M,p_k, g/\epsilon_k,J,\omega/\epsilon_k)$ pointed converge in the sense of Cheeger-Gromov to $(T_p M,0, g_p,J_p,\omega_p)$.
\hfill$\Box$
\end{co}

\begin{prop}
\label{bb}
In proposition \ref{ba}, let $L_k$ be a Lagrangian submanifold in $(M_k,g_k,J_k,\omega_k)$ with uniformly bounded local geometry. By possibly taking subsequence, $L_k$ converge to a Lagrangian submanifold $L_o$ in $(M_o,g_o,J_o,\omega_o)$ with bounded local geometry, while $(M_k,z_k, g_k,J_k,\omega_k)$ pointed converge in the sense of Cheeger-Gromov to $(M_o, z_o, g_o, J_o, \omega_o)$.
\hfill$\Box$
\end{prop}

\begin{co}
\label{bf}
In corollary \ref{be}, let $L$ be a Lagrangian submanifold in $(M,g,J,\omega)$ with uniformly bounded local geometry. For a sequence of points $p_k \in M$, by possibly taking subsequence, $L$ converge to a Lagrangian submanifold $L_o$ in $(M_o, g_o, J_o, \omega_o)$ with bounded local geometry, while $(M,p_k, g,J,\omega)$ pointed converge in the sense of Cheeger-Gromov to $(M_o, p_o, g_o, J_o, \omega_o)$.
\hfill$\Box$
\end{co}

\begin{prop}
\label{bc}
In proposition \ref{ba}, assume $M_k$ has smooth boundary $\partial M_k$ and $(M_k,\partial M_k, g_k,J_k,\omega_k)$ is complete with uniformly bounded local geometry. By possibly taking subsequence, $(M_k, \partial M_k, z_k, g_k,J_k,\omega_k)$ pointed converge to $(M_o, \partial M_o, z_o, g_o, J_o, \omega_o)$, which is complete with bounded local geometry.
\hfill$\Box$
\end{prop}

\begin{re}
An important fact concerning pointed convergence in the sense of Cheeger-Gromov (stated in the context of proposition \ref{ba} for simplicity) is that the limit $M_o$ is determined by $B_{R_k}^{g_k} (M_k, z_k) \subset M_k$ as long as $\displaystyle \lim R_k = +\infty$. More detail of pointed convergence is illustrated in remark \ref{ae} through a concrete example.
\end{re}

{\bf Global Gromov-Hausdorff limit:} Pointed convergent limit $(M_o, g_o, z_o)$ of $(M_k, g_k, z_k)$ is local in nature, and can be viewed as a local component of a ``global limit" of $(M_k, g_k)$. It is straightforward to show that for a different set of marked points $z'_k \in M_k$ satisfying ${\rm Dist}_{g_k} (z_k, z'_k) \leq C$ for certain $C>0$, the corresponding pointed limit $(M'_o, g'_o, z'_o)$ can be canonically identified with $(M_o, g_o, z_o)$ up to isometry. Accordingly, we may define the global limit $(M_o, g_o)$ to be a collection of pointed limits $(M^i_o, g^i_o, z^i_o)$ (as local components of the global limit) satisfying certain condition (e.g. with maximal dimension, etc.) determined by the marked point sets $\{z_k^i\}$ indexed by $i \in I$, such that $\lim {\rm Dist}_{g_k} (z^i_k, z^j_k) = +\infty$ for $i\not=j \in I$, furthermore, we assume $I$ is maximal in the sense that no more pointed limit can be added to the collection with all the conditions satisfied. Such global limit $(M_o, g_o)$ is most useful when it contains only finite many local components, which in particular implies that the global limit is essentially unique up to isometry. The finiteness of $I$ can quite often be shown through certain variations of the volume method, the index method or a combination of the two. The volume method start with a uniform finite upper bound of the volume of $(M_k, g_k)$, then achieve finiteness of $I$ by proving the existence of a positive lower bound of the volume of the pointed limits. The index method start with a uniform finite upper bound of certain subadditive index of $(M_k, g_k)$, then achieve finiteness of $I$ by proving the index of any pointed limit is a positive integer.\\

\se{Conformal models of a disk}
Conformal models constructed in this section are essentially the same as the ones used in \cite{FO}. These conformal models are formulated in a form that is more convenient for our discussion and the construction is more elementary in nature.\\

Recall that $\mathbb{D}_{[\scriptsize{\textsc{n}}]}$ is the closed unit disk removing $n$ marked points $\{\textsc{p}_k\}_{k=1}^{\scriptsize{\textsc{n}}}$ circularly located in $\partial \mathbb{D}$. We are looking for graph like conformal model $(\Theta, g_\Theta)$ of $\mathbb{D}_{[\scriptsize{\textsc{n}}]}$, by which we mean that there is a homeomorphism from $\mathbb{D}_{[\scriptsize{\textsc{n}}]}$ to $\Theta$ that is biholomorphic in the interior of $\mathbb{D}_{[\scriptsize{\textsc{n}}]}$.\\
\begin{prop}
\label{aa}
For any $\mathbb{D}_{[\scriptsize{\textsc{n}}]} \in {\cal I}_{\scriptsize{\textsc{n}}}$, there exists a conformal model $(\Theta \cong \mathbb{D}_{[\scriptsize{\textsc{n}}]}, g_\Theta)$, a complete metric tree graph $\Gamma \in Gr_{\scriptsize{\textsc{n}}}$ and a decomposition

\[
\Theta = \bigcup_{\nu \in C^0_{\rm int}(\Gamma) \cup  C^1(\Gamma)} \Theta (\nu),
\]

where for $\nu \in C^0_{\rm int}(\Gamma)$, $(\Theta (\nu), g_\Theta)$ has bounded geometry (namely, diameter and curvature are bounded and the injective radius is bounded from below), with bounds only depend on $\textsc{n}$; for $\nu \in C^1(\Gamma)$, $(\Theta (\nu), g_\Theta) \cong ((0, \gamma_\nu) \times [0, 1], g_0)$, where $g_0$ is the standard flat metric. Furthermore, there is a continuous map $\phi: \Theta \rightarrow \Gamma$, where $\Theta (\nu)$ is surjectively mapped to $\nu$ for $\nu \in C^*(\Gamma) = C^0_{\rm int}(\Gamma) \cup  C^1(\Gamma)$ in the obvious way.
\end{prop}
{\bf Proof:} For the convenience of doubling the disk into a Riemann surface, we like to view $\mathbb{D}$ as the upper half of $(\mathbb{CP}^1, \mathbb{RP}^1)$ with the marked points $\{\textsc{p}_i\}_{i=1}^{\scriptsize{\textsc{n}}}$ in $\partial \mathbb{D} \cong \mathbb{RP}^1$. We will prove the proposition by induction for the conformal model $\breve{\Theta}$ of $\mathbb{CP}^1_{[\scriptsize{\textsc{n}}]} = \mathbb{CP}^1 \setminus \{\textsc{p}_i\}_{i=1}^{\scriptsize{\textsc{n}}}$ so that $\Theta$ is the upper half of $\breve{\Theta}$.\\

We will start with a multiple of the Fubini-Study metric on $\mathbb{CP}^1$. Assume that there is already a conformal model $(\breve{\Theta}', g')$ for $\mathbb{CP}^1_{[\scriptsize{\textsc{n}}-1]}$ with $\Gamma' \in Gr_{\scriptsize{\textsc{n}} -1}$ and $\textsc{n}-1$ ends corresponding to $\{\textsc{p}_i\}_{i=1}^{\scriptsize{\textsc{n}}-1}$ such that $g'$ is invariant under complex conjugation. Without loss of generality, we may assume that $\gamma_\nu>3$ for $\nu \in C^1_{\rm int} (\Gamma')$. (If $\gamma_\nu\leq 3$ for $\nu \in C^1_{\rm int} (\Gamma')$, we may take the union of $\breve{\Theta}' (\nu)$ and the 2 adjacent components to form a new component that effectively collapses the leg $\nu$ into a vertex.)\\

For the position of $\textsc{p}_{\scriptsize{\textsc{n}}}$ in $(\breve{\Theta}', g')$, there is either a $\nu \in C^1 (\Gamma')$ such that $\textsc{p}_{\scriptsize{\textsc{n}}} \in \breve{\Theta}' (\nu)$ and ${\rm Dist} (\textsc{p}_{\scriptsize{\textsc{n}}}, \partial \breve{\Theta}' (\nu)) > 1$, or a $\nu \in C^0_{\rm int} (\Gamma')$ such that ${\rm Dist} (\textsc{p}_{\scriptsize{\textsc{n}}}, \breve{\Theta}' (\nu)) \leq 1$, in which case, by suitably expanding $\breve{\Theta}'(\nu)$, one can ensure $\textsc{p}_{\scriptsize{\textsc{n}}} \in \breve{\Theta}' (\nu)$ and ${\rm Dist} (\textsc{p}_{\scriptsize{\textsc{n}}}, \partial \breve{\Theta}' (\nu)) > 1$. In sum, we have $\nu \in C^* (\Gamma')$ such that $\textsc{p}_{\scriptsize{\textsc{n}}} \in \breve{\Theta}' (\nu)$ and $B_1 (\textsc{p}_{\scriptsize{\textsc{n}}}) \subset \breve{\Theta}' (\nu)$. \\

Let $z$ be an affine coordinate on $\mathbb{CP}^1$ that is real on $\mathbb{RP}^1$, $z(\textsc{p}_{\scriptsize{\textsc{n}}})=0$ and $\{|z|<3\}\subset B_1 (\textsc{p}_{\scriptsize{\textsc{n}}})$. Let $g_{\scriptsize{\textsc{n}}}$ be the cylindrical metric on the punctured disk $\{0<|z|<3\}$ determined by the cylindrical coordinate $t = \frac{\log z}{\pi}$. $g = (1-\rho) g' + \rho g_{\scriptsize{\textsc{n}}}$ gives us the desired conformal model $(\breve{\Theta}, g)$ of $\mathbb{CP}^1_{[\scriptsize{\textsc{n}}]}$ with $\textsc{n}$ ends corresponding to $\{\textsc{p}_i\}_{i=1}^{\scriptsize{\textsc{n}}}$ that is invariant under complex conjugation, where $0\leq \rho \leq 1$, $\rho|_{|z| \geq 2}=0$ and $\rho|_{|z| \leq 1}=1$.\\

When $\nu \in C^0_{\rm int} (\Gamma')$, we may get $\Gamma$ from $\Gamma'$ by replacing $\nu \in C^0_{\rm int} (\Gamma')$ with $\nu \in C^0_{\rm int} (\Gamma)$ and $\nu_{\scriptsize{\textsc{n}}} \in C^1_{\rm ext} (\Gamma)$ connecting $\nu$ and $\textsc{p}_{\scriptsize{\textsc{n}}}$ such that $\breve{\Theta}(\nu_{\scriptsize{\textsc{n}}}) = \{0< |z| <1\}$ with the cylindrical coordinate $t$ and $\breve{\Theta}(\nu) = \breve{\Theta}'(\nu) \setminus \breve{\Theta}(\nu_{\scriptsize{\textsc{n}}})$.\\

When $\nu \in C^1 (\Gamma')$, it is straightforward to get the decomposition $\breve{\Theta}'(\nu) = \breve{\Theta}(\nu_0) \cup \breve{\Theta}(\nu_1) \cup \breve{\Theta}(\nu_2) \cup \breve{\Theta}(\nu_{\scriptsize{\textsc{n}}})$, so that $\breve{\Theta}(\nu_{\scriptsize{\textsc{n}}}) = \{0< |z| <1\}$ with the cylindrical coordinate $t$, $\{1\leq |z| \leq 2\} \subset \breve{\Theta}(\nu_0)$ and $\breve{\Theta}(\nu_1)$, $\breve{\Theta}(\nu_2)$ are cylindrical under the metric $g$. We may get $\Gamma$ from $\Gamma'$ by replacing $\nu \in C^1 (\Gamma')$ with $\nu_1, \nu_2 \in C^1 (\Gamma)$, $\nu_0 \in C^0_{\rm int} (\Gamma)$ and $\nu_{\scriptsize{\textsc{n}}} \in C^1_{\rm ext} (\Gamma)$ connecting $\nu_0$ and $\textsc{p}_{\scriptsize{\textsc{n}}}$.
\hfill$\Box$\\

For $\nu \in C^1 (\Gamma)$, it is convenient to identify $(0, \gamma_\nu) \times [0, 1]$ with a subset in $\mathbb{C} \cong \mathbb{R}^2$ and to think of the cylindrical coordinate as a complex coordinate $t_\nu: \Theta (\nu) \rightarrow \mathbb{C}$ (as we already did in the proof of proposition \ref{aa}). For $\nu \in C^0_{\rm int}(\Gamma)$, one may choose $t_\nu$ to be the cylindrical coordinate corresponding to any of the leg in $C^1 (\Gamma)$ that is attached to $\nu$. Since each of $t_\nu$ as defined in the proof of proposition \ref{aa} is the logarithm of an affine coordinate on $\mathbb{CP}^1$ that is real on $\mathbb{RP}^1$, the proof of proposition \ref{aa} in fact implies the following:\\
\begin{prop}
\label{ad}
For any $\nu_1,\nu_2 \in C^0_{\rm int}(\Gamma) \cup  C^1(\Gamma)$, there exists $a,b,c,d \in \mathbb{R}$ such that\\

$\hspace{1.5in} \displaystyle t_{\nu_2} = - \frac{1}{\pi} \log \frac{a e^{-\pi t_{\nu_1}} +b}{c e^{-\pi t_{\nu_1}} +d}$.
\hfill$\Box$\\
\end{prop}

\begin{re}
A more geometric and canonical proof of proposition \ref{aa} starts with a multiple of the hyperbolic Poincare metric on $\mathbb{CP}^1_{[\scriptsize{\textsc{n}}]}$, then modifies the hyperbolic cusps (resp. cylinders) explicitly into cylindrical ends (resp. cylinders). One purpose of our proof is to illustrate the elementary nature of proposition \ref{aa} without referring to the powerful hyperbolic Poincare metric. Of course it is also nice to have the coordinate transformation formula in proposition \ref{ad} as a byproduct.
\end{re}

In this paper, a conformal model $(\Theta, g_\Theta)$ is always referred to one associated with a metric tree graph $\Gamma$, the map $\phi: \Theta \rightarrow \Gamma$ (as defined in proposition \ref{aa}) and coordinates (as in proposition \ref{ad}).\\
\begin{prop}
\label{ab}
For a sequence of conformal model $(\Theta_k, g_{\Theta_k})$ with marked point $t^\circ_k \in \Theta_k$ and at most $\textsc{n}$ ends, by possibly taking subsequence, $(\Theta_k, t^\circ_k, g_{\Theta_k})$ is pointed convergent to $(\Theta_o, t^\circ_o, g_{\Theta_o})$ in the sense of Cheeger-Gromov, where $(\Theta_o, g_{\Theta_o})$ is a conformal model with marked point $t^\circ_o \in \Theta_o$ and at most $\textsc{n}$ ends.
\hfill$\Box$
\end{prop}

With map $\phi: \Theta \rightarrow \Gamma$, it is easy to see that the Gromov-Hausdorff distance between the conformal model $(\Theta, g_\Theta)$ and the complete metric tree graph $(\Gamma, \gamma)$ as metric spaces is bounded with bound only depending on $\textsc{n}$. Consequently, the Gromov-Hausdorff distance between $(\Theta, \epsilon g_\Theta)$ and the complete metric tree graph $(\Gamma, \epsilon \gamma)$ is $O(\epsilon)$ for any $\epsilon >0$.\\
\begin{prop}
\label{ac}
For $\lim \epsilon_k =0$ and a sequence of conformal model $(\Theta_k, g_{\Theta_k})$ with marked point $t^\circ_k \in \Theta_k$ and at most $n$ ends, by possibly taking subsequence, $(\Theta_k, t^\circ_k, \epsilon_k^2 g_{\Theta_k})$ (resp. $(\Gamma_k, \tau^\circ_k, \epsilon_k \gamma_k)$) is pointed Gromov-Hausdorff convergent to a complete metric tree graph $(\Gamma^\circ, \tau^\circ, \gamma^\circ)$ with marked point and at most $n$ ends, where $\tau^\circ_k = \phi (t^\circ_k)$.
\hfill$\Box$
\end{prop}

\begin{re}
\label{ae}
{\bf [\textit{pointed convergence}]:} The Gromov-Hausdorff distance was originally defined for compact metric spaces. More rigorously, the statement: ``$(\Theta_k, t^\circ_k, \epsilon_k^2 g_{\Theta_k})$ (resp. $(\Gamma_k, \tau^\circ_k, \epsilon_k \gamma_k)$) is pointed Gromov-Hausdorff convergent to $(\Gamma^\circ, \tau^\circ, \gamma^\circ)$" involve 2 limiting steps. In the first step, for $R>0$, we have the Gromov-Hausdorff convergence for compact metric spaces:

\[
(\Gamma^\circ_R, \tau^\circ, \gamma^\circ_R) := \lim_{k \rightarrow +\infty} (B_R(\Theta_k, t^\circ_k), t^\circ_k, \epsilon_k^2 g_{\Theta_k}) = \lim_{k \rightarrow +\infty} (B_R(\Gamma_k, \tau^\circ_k), \tau^\circ_k, \epsilon_k \gamma_k).
\]

In the second step, for any $R_2 > R_1>0$, we have the isometric embedding $(\Gamma^\circ_{R_1}, \tau^\circ, \gamma^\circ_{R_1}) \rightarrow (\Gamma^\circ_{R_2}, \tau^\circ, \gamma^\circ_{R_2})$. Then we have the following limit as infinite union:

\[
(\Gamma^\circ, \tau^\circ, \gamma^\circ) = \lim_{R \rightarrow +\infty} (\Gamma^\circ_R, \tau^\circ, \gamma^\circ_R) = \bigcup_{R >0} (\Gamma^\circ_R, \tau^\circ, \gamma^\circ_R).
\]

In this paper, whenever pointed convergence is considered, such process is always implicitly implied without explicitly mentioned each time. Similar to the case of smooth manifolds, the convergence can be further strengthened to $C^k$-convergence over the smooth part of $\Gamma^\circ$ in the sense of Cheeger-Gromov.\\
\end{re}

{\samepage
\stepcounter{subsection}
{\bf \S \thesubsection\ Holomorphic functions of one complex variable}\\
\nopagebreak

In this section,} we will summarize some well known results for holomorphic functions of one complex variable.

\begin{prop}
\label{ea}
Let $w: \mathbb{C} \rightarrow \mathbb{C}^n$ be a holomorphic function with $|Dw| \leq C$, then $w$ is linear.
\hfill$\Box$
\end{prop}

The following 2 propositions are consequences of proposition \ref{ea} through extension by reflection.

\begin{prop}
\label{eb}
Let $w: (\mathbb{C}_+, \mathbb{R}) \rightarrow (\mathbb{C}^n, \mathbb{R}^n)$ be a holomorphic function with $|Dw| \leq C$, then $w$ is linear.
\hfill$\Box$
\end{prop}

\begin{prop}
\label{ec}
Let $\Theta = \{t\in \mathbb{C}: 0 \leq t^\Im \leq 1\}$ and $w: \Theta \rightarrow \mathbb{C}^n$ be a holomorphic function such that $w (\partial \Theta)$ is in shifts of $\mathbb{R}^n$ and $|Dw| \leq C$, then $w$ is linear.
\hfill$\Box$
\end{prop}

\begin{prop}
\label{ed}
Let $\Theta$ be a conformal model of the disk with cylindrical coordinates $\{t_i\}_{i=1}^{\scriptsize{\textsc{n}}}$ at ends $\{\textsc{p}_i \}_{i=1}^{\scriptsize{\textsc{n}}}$, and $w: (\Theta, \partial \Theta) \rightarrow (\mathbb{C}^n, \mathbb{R}^n)$ be a holomorphic function such that for certain $0 \leq c<1$, $|Dw (t_i)| \leq C e^{c\pi t^\Re_i}$ near the end $\textsc{p}_i$, then $w$ is a constant map.
\end{prop}
{\bf Proof:} Since $\Theta$ has $\textsc{n}$ ends, its holomorphic double $\mathbb{CP}_{[\scriptsize{\textsc{n}}]}^1$ is biholomorphic to $\mathbb{CP}^1$ removing $\textsc{n}$ punctures $\{\textsc{p}_i \}_{i=1}^{\scriptsize{\textsc{n}}}$ on $\mathbb{RP}^1$. By reflection, $w$ can be extended to $w: (\mathbb{CP}_{[\scriptsize{\textsc{n}}]}^1, \mathbb{RP}_{[\scriptsize{\textsc{n}}]}^1) \rightarrow (\mathbb{C}^n, \mathbb{R}^n)$.\\

$\hat{t}_i = e^{-\pi t_i}$ is the natural coordinate of $\mathbb{CP}_{[\scriptsize{\textsc{n}}]}^1$ near the puncture $\textsc{p}_i$. $|\frac{\partial w}{\partial \hat{t}_i}| = |\frac{1}{\pi \hat{t}_i} \frac{\partial w}{\partial t_i}| \leq \frac{1}{\pi|\hat{t}_i|^{1+c}}$ implies that $|w| \leq C|\hat{t}_i|^{-c}$ (resp. $|w| \leq C|\log |\hat{t}_i||$ when $c=0$). Consequently, $w$ is holomorphic at the puncture $\hat{t}_i =0$. Namely $w$ can be extended to $w: \mathbb{CP}^1 \rightarrow \mathbb{C}^n$, which has to be a constant map.
\hfill$\Box$

\se{Estimates of pseudo holomorphic disks and convergence}
In this section, we discuss local interior and boundary $C^{1,\alpha}$-estimates of pseudo holomorphic disks with Lagrangian boundary condition. The proofs are elementary generalizations of the standard method in one complex variable using Cauchy integration formula. As application, we establish convergence results of pseudo holomorphic disks with Lagrangian boundary condition.\\

Let $\mathbb{D}'$ be a smaller disk in the unit disk $\mathbb{D} \subset \mathbb{C}$. Let $(\mathbb{R}^{2n}, J)$ be an almost complex manifold, with $J_w$ smooth on $w\in\mathbb{R}^{2n}$ be the complex structure on $T_w \mathbb{R}^{2n} \cong \mathbb{R}^{2n}$. It is straightforward to construct a family of linear identifications $A_w: (\mathbb{R}^{2n}, J_w) \cong (\mathbb{C}^n, i)$ that is smoothly depending on $w\in \mathbb{R}^{2n}$. We start with the local interior $C^{1,\alpha}$-estimate.\\
\begin{prop}
For a pseudo holomorphic map $f: \mathbb{D} \rightarrow (\mathbb{R}^{2n}, J)$,
\begin{equation}
\label{nf}
|Df|_{C^{\alpha} (\mathbb{D}')} \leq C_\alpha \sup_{z\in \mathbb{D}} (|Df(z)| + |Df(z)|^2).
\end{equation}
\end{prop}
{\bf Proof:} It is straightforward to derive:

\[
f(z) = \frac{1}{2\pi i}\int_{\partial \mathbb{D}} \frac{f(\zeta)}{(\zeta -z)} d\zeta - \frac{1}{2\pi i}\int_{\mathbb{D}} \frac{df(\zeta)}{(\zeta -z)} d\zeta.
\]

The pseudo holomorphic condition can be expressed as $A_{f(\zeta)} df(\zeta) \wedge d\zeta =0$. Hence

\[
A_w f(z) = \frac{1}{2\pi i}\int_{\partial \mathbb{D}} \frac{A_w f(\zeta)}{(\zeta -z)} d\zeta - \frac{1}{2\pi i}\int_{\mathbb{D}} \frac{(A_w - A_{f(\zeta)}) df(\zeta)}{(\zeta -z)} d\zeta.
\]

Assume $w = f(z)$, then

\[
A_w Df (z) = \frac{Dz}{2\pi i} \left[ \Xi_1(z) - \Xi_2(z) \right],
\]
\[
\Xi_1(z) = \int_{\partial \mathbb{D}} \frac{A_w(f(\zeta) - f(z_0))}{(\zeta -z)^2} d\zeta,\ \Xi_2(z) = \int_{\mathbb{D}} \frac{(A_w - A_{f(\zeta)}) df(\zeta)}{(\zeta -z)^2} d\zeta.
\]

For the estimate of $\Xi_2(z)$, assume $z_1,z_2 \in \mathbb{D}'$ and $\zeta \in \mathbb{D}_1 := \{ \zeta \in \mathbb{D}: |\zeta - z_1| \leq |\zeta - z_2|\}$, then $|z_2 - z_1| \leq 2|\zeta - z_2|$. We have

\[
\left| \frac{A_{f(z_2)} - A_{f(\zeta)}}{(\zeta -z_2)^2} - \frac{A_{f(z_1)} - A_{f(\zeta)}}{(\zeta -z_1)^2} \right| \leq |A_{f(z_1)} - A_{f(\zeta)}|\left| \frac{1}{(\zeta -z_2)^2} - \frac{1}{(\zeta -z_1)^2} \right|
\]
\[
+ \left| \frac{A_{f(z_1)} - A_{f(z_2)}}{(\zeta -z_2)^2} \right| \leq C|Df|_\mathbb{D} \left( \frac{1}{|\zeta -z_2|^{1+\alpha}} + \frac{1}{|\zeta -z_1|^{1+\alpha}} \right) |z_1 -z_2|^\alpha.
\]

Here we used

\[
\left| \frac{1}{(\zeta -z_1)^l} - \frac{1}{(\zeta -z_2)^l} \right| \leq \frac{C|z_1 -z_2|}{|\zeta -z_1|^l|\zeta -z_2|}, \mbox{ if } |\zeta -z_1| \leq |\zeta -z_2|.
\]

The estimate can be symmetrically done for $\zeta \in \mathbb{D}\setminus \mathbb{D}_1$. Hence

\[
|\Xi_2 (z_2) - \Xi_2 (z_1)| \leq \int_\mathbb{D} \left| \frac{A_{f(z_2)} - A_{f(\zeta)}}{(\zeta -z_2)^2} - \frac{A_{f(z_1)} - A_{f(\zeta)}}{(\zeta -z_1)^2} \right| \left| \frac{\partial f}{\partial \bar{\zeta}} \right| d\zeta d\bar{\zeta}
\]
\[
\leq C|z_1 -z_2|^\alpha |Df|_\mathbb{D}^2 \int_\mathbb{D} \left( \frac{1}{|\zeta -z_2|^{1+\alpha}} + \frac{1}{|\zeta -z_1|^{1+\alpha}} \right)  d\zeta d\bar{\zeta}.
\]

Consequently,

\begin{equation}
\label{nfa}
\ \ \frac{|\Xi_2 (z_2) - \Xi_2 (z_1)|}{|z_2-z_1|^\alpha} \leq C|Df|_\mathbb{D}^2, \mbox{ for } z_1,z_2 \in \mathbb{D}'.
\end{equation}

For the estimate of $\Xi_1(z)$, assume $z_1,z_2 \in \mathbb{D}'$, $\zeta \in \partial \mathbb{D}$. Then

\[
\left| \frac{A_{f(z_2)}}{(\zeta -z_2)^2} - \frac{A_{f(z_1)}}{(\zeta -z_1)^2} \right| \leq C (1 + |Df|_\mathbb{D}) |z_1 -z_2|.
\]

Hence

\[
|\Xi_1 (z_2) - \Xi_1 (z_1)| \leq \int_{\partial \mathbb{D}} \left| \frac{A_{f(z_2)}}{(\zeta -z_2)^2} - \frac{A_{f(z_1)}}{(\zeta -z_1)^2} \right| \left| f(\zeta) - f(z_0)\right| |d\zeta|
\]
\[
\leq C|z_1 -z_2| (1 + |Df|_\mathbb{D}) |Df|_\mathbb{D} .
\]

Consequently,

\begin{equation}
\label{nfe}
\frac{|\Xi_1 (z_2) - \Xi_1 (z_1)|}{|z_2-z_1|^\alpha} \leq C(1 + |Df|_\mathbb{D}) |Df|_\mathbb{D}, \mbox{ for } z_1,z_2 \in \mathbb{D}'.
\end{equation}

Estimates (\ref{nfa}) and (\ref{nfe}) together imply that for $z_1,z_2 \in \mathbb{D}'$

\begin{equation}
\label{nff}
\frac{|A_{f(z_2)} Df (z_2) - A_{f(z_1)} Df (z_1)|}{|z_2-z_1|^\alpha} \leq C|Df|_\mathbb{D} (1 + |Df|_\mathbb{D}).
\end{equation}

\[
\frac{A_{f(z_2)} Df (z_2) - A_{f(z_1)} Df (z_1)}{|z_2-z_1|^\alpha} = \frac{A_{f(z_2)} - A_{f(z_1)}}{|z_2-z_1|^\alpha} Df (z_1) + A_{f(z_2)} \frac{Df (z_2) - Df (z_1)}{|z_2-z_1|^\alpha}
\]

implies

\[
\frac{|Df (z_2) - Df (z_1)|}{|z_2-z_1|^\alpha} \leq C\frac{|A_{f(z_2)} Df (z_2) - A_{f(z_1)} Df (z_1)|}{|z_2-z_1|^\alpha} + C|Df|_\mathbb{D} \frac{|f (z_2) - f (z_1)|}{|z_2-z_1|^\alpha}.
\]

Apply the estimate (\ref{nff}), we have

\[
\frac{|Df (z_2) - Df (z_1)|}{|z_2-z_1|^\alpha} \leq C|Df|_\mathbb{D} (1 + |Df|_\mathbb{D}), \mbox{ for } z_1,z_2 \in \mathbb{D}',
\]

which implies the estimate (\ref{nf}).
\hfill$\Box$\\

This estimate is responsible for the following useful convergence result.

\begin{prop}
\label{nh}
For a sequence of pseudo holomorphic maps $f_k: (\Theta_k, t_k, g_{\Theta_k}) \rightarrow (M_k, z_k, g_k)$ with uniform $C^1$-bound $\displaystyle \sup_{t\in \Theta_k} |D f_k (t)| \leq C$, where $(M_k,g_k)$ is an almost K\"{a}hler manifold, assume that
\[
(\Theta_o, t_o, g_{\Theta_o}) = \displaystyle \lim (\Theta_k, t_k, g_{\Theta_k}) \mbox{ and } (M_o, z_o, g_o) = \displaystyle \lim (M_k, z_k, g_k)
\]
in the sense of Cheeger-Gromov. Then there exists a pseudo holomorphic map $f_o: (\Theta_o, t_o, g_{\Theta_o}) \rightarrow (M_o, z_o, g_o)$ such that $f_o = \displaystyle \lim f_k$. The convergence is in $C^{1,\alpha}$ for $0< \alpha <1$ on any compact subset in the interior of $\Theta_o$. In particular, if $t_o$ is in the interior of $\Theta_o$, then $Df_o (t_o) = \displaystyle \lim Df_k (t_k)$.
\hfill$\Box$
\end{prop}

Let $\mathbb{D}_+$ (resp. $\mathbb{D}'_+$) be the half disk as intersection of $\mathbb{D}$ (resp. $\mathbb{D}'$) and the upper half plane. Let $\partial_0 \mathbb{D}_+$ be the intersection of $\mathbb{D}$ and the real axis. Let $L$ be a smooth middle dimensional totally real submanifold of $\mathbb{C}^n$ passing through 0. For the boundary estimate, without loss of generality, we may assume $L = \mathbb{R}^n \times \{0\} \subset \mathbb{R}^{2n} \cong \mathbb{R}^n \times \mathbb{R}^n$. Then we have

\begin{prop}
For a pseudo holomorphic map $f: (\mathbb{D}_+, \partial_0 \mathbb{D}_+) \rightarrow (\mathbb{R}^{2n}, L, J)$,

\begin{equation}
\label{nk}
|Df|_{C^{\alpha} (\mathbb{D}'_+)} \leq C_\alpha \sup_{z\in \mathbb{D}_+} (|Df(z)| + |Df(z)|^2).
\end{equation}
\end{prop}
{\bf Proof:} Without loss of generality, we may assume $A_w (L) = \mathbb{R}^n \subset \mathbb{C}^n$ for $w \in L$. This condition implies that

\begin{equation}
\label{nka}
A^\Im_{w_1} (w_2 - w_1) = 0, \mbox{ for } w_1,w_2 \in L.
\end{equation}

For $z \in \mathbb{D}_+$, we have

\[
A_w f(z) = \frac{1}{2\pi i}\int_{\partial \mathbb{D}_+} \frac{A_w f(\zeta)}{(\zeta -z)} d\zeta - \frac{1}{2\pi i}\int_{\mathbb{D}_+} \frac{(A_w - A_{f(\zeta)}) df(\zeta)}{(\zeta -z)} d\zeta.
\]
\[
0 = \frac{1}{2\pi i}\int_{\partial \mathbb{D}_-} \frac{\bar{A}_w f(\bar{\zeta})}{(\zeta -z)} d\zeta - \frac{1}{2\pi i}\int_{\mathbb{D}_-} \frac{(\bar{A}_w - \bar{A}_{f(\bar{\zeta})}) df(\bar{\zeta})}{(\zeta -z)} d\zeta.
\]

Add the 2 equations, we have

\[
A_w f(z) = \frac{1}{\pi}\int_{\partial_0 \mathbb{D}_+} \frac{A^\Im_w f(\zeta)}{(\zeta -z)} d\zeta + \frac{1}{2\pi i} \left[ \int_{\partial_1 \mathbb{D}_+} \frac{A_w f(\zeta)}{(\zeta -z)} d\zeta + \int_{\partial_1 \mathbb{D}_-} \frac{\bar{A}_w f(\bar{\zeta})}{(\zeta -z)} d\zeta \right]
\]
\[
- \frac{1}{2\pi i} \left[ \int_{\mathbb{D}_+} \frac{(A_w - A_{f(\zeta)}) df(\zeta)}{(\zeta -z)} d\zeta + \int_{\mathbb{D}_-} \frac{(\bar{A}_w - \bar{A}_{f(\bar{\zeta})}) df(\bar{\zeta})}{(\zeta -z)} d\zeta \right].
\]

Assume $w = f(z)$, then

\[
A_w Df (z) = \frac{Dz}{2\pi i}\int_{\partial \mathbb{D}_+} \frac{A_w(f(\zeta) - f(z^\Re))}{(\zeta -z)^2} d\zeta - \frac{Dz}{2\pi i}\int_{\mathbb{D}_+} \frac{(A_w - A_{f(\zeta)}) df(\zeta)}{(\zeta -z)^2} d\zeta.
\]
\[
= \frac{Dz}{\pi} \Xi_0 (z) + \frac{Dz}{2\pi i} [\Xi_1 (z) - \Xi_2 (z)],
\]

where

\[
\Xi_0 (z) = \int_{\partial_0 \mathbb{D}_+} \frac{A^\Im_w (f(\zeta) - f(z^\Re))}{(\zeta -z)^2} d\zeta,
\]
\[
\Xi_2 (z) = \int_{\mathbb{D}_+} \frac{(A_w - A_{f(\zeta)}) df(\zeta)}{(\zeta -z)^2} d\zeta + \int_{\mathbb{D}_-} \frac{(\bar{A}_w - \bar{A}_{f(\bar{\zeta})}) df(\bar{\zeta})}{(\zeta -z)^2} d\zeta,
\]
\[
\Xi_1 (z) = \int_{\partial_1 \mathbb{D}_+} \frac{A_w(f(\zeta) - f(z^\Re))}{(\zeta -z)^2} d\zeta + \int_{\partial_1 \mathbb{D}_-} \frac{\bar{A}_w(f(\bar{\zeta}) - f(z^\Re))}{(\zeta -z)^2} d\zeta .
\]

Since $f(\zeta), f(z^\Re)$ are in $L$, (\ref{nka}) implies that $A^\Im_{f(z^\Re)} (f(\zeta) - f(z^\Re)) =0$. Consequently,

\[
\Xi_0 (z) = \int_{\partial_0 \mathbb{D}_+} \frac{(A^\Im_w - A^\Im_{f(z^\Re)}) (f(\zeta) - f(z^\Re))}{(\zeta -z)^2} d\zeta.
\]

Assume $z_1,z_2 \in \mathbb{D}'_+$, $\zeta \in \partial_0 \mathbb{D}_+$ and $|\zeta - z_1| \leq |\zeta - z_2|$, then $|z_2 - z_1| \leq 2|\zeta - z_2|$. We have

\[
\left| \frac{(A^\Im_{f(z_2)} - A^\Im_{f(z_2^\Re)}) (f(\zeta) - f(z_2^\Re))}{(\zeta -z_2)^2} - \frac{(A^\Im_{f(z_1)} - A^\Im_{f(z_1^\Re)}) (f(\zeta) - f(z_1^\Re))}{(\zeta -z_1)^2} \right|
\]
\[
\leq |(A^\Im_{f(z_1)} - A^\Im_{f(z_1^\Re)}) (f(\zeta) - f(z_1^\Re))|\left| \frac{1}{(\zeta -z_2)^2} - \frac{1}{(\zeta -z_1)^2} \right|
\]
\[
+ \left| \frac{(A^\Im_{f(z_2)} - A^\Im_{f(z_2^\Re)}) (f(\zeta) - f(z_2^\Re)) - (A^\Im_{f(z_1)} - A^\Im_{f(z_1^\Re)}) (f(\zeta) - f(z_1^\Re))}{(\zeta -z_2)^2} \right|
\]
\[
\leq C |Df|^2_{\mathbb{D}_+} |\zeta -z_2|^{-\alpha} |z_1 -z_2|^\alpha.
\]

Hence

\[
|\Xi_0 (z_2) - \Xi_0 (z_1)| \leq C |Df|^2_{\mathbb{D}_+} |z_1 -z_2|^\alpha \int_{\partial_0 \mathbb{D}_+} |\zeta -z_2|^{-\alpha} d\zeta\leq C |Df|^2_{\mathbb{D}_+} |z_1 -z_2|^\alpha.
\]

Consequently,

\begin{equation}
\label{nkb}
\frac{|\Xi_0 (z_2) - \Xi_0 (z_1)|}{|z_2-z_1|^\alpha} \leq C|Df|^2_{\mathbb{D}_+}, \mbox{ for } z_1,z_2 \in \mathbb{D}'_+.
\end{equation}

The terms $\Xi_1(z)$ and $\Xi_2(z)$ in the expression of $A_w Df (z)$ can be estimated in the same way as in the proof of the estimate (\ref{nf}). Consequently,

\[
\ \ \frac{|A_{f(z_2)} Df (z_2) - A_{f(z_1)} Df (z_1)|}{|z_2-z_1|^\alpha} \leq C(|Df|_{\mathbb{D}_+} + |Df|_{\mathbb{D}_+}^2), \mbox{ for } z_1,z_2 \in \mathbb{D}'_+.
\]

Following rest of the proof of the estimate (\ref{nf}), we get the estimate (\ref{nk}).
\hfill$\Box$\\

This boundary estimate extends our convergence to the boundary.

\begin{prop}
\label{ni}
In proposition \ref{nh}, assume that $\partial_0 \Theta_o = \displaystyle \lim \partial_0 \Theta_k$ and $L_o = \displaystyle \lim L_k$ in the context of Cheeger-Gromov limits in proposition \ref{nh}, where $L_k$ (resp. $L_o$) is a Lagrangian submanifold of $(M_k, \omega_k)$ (resp. $(M_o, \omega_o)$) and $\partial_0 \Theta_k$ (resp. $\partial_0 \Theta_o$) is a boundary component of $\Theta_k$ (resp. $\Theta_o$) such that $f_k(\partial_0 \Theta_k) \subset L_k$. Then $f_k$ converges to $f_o$ in $C^{1,\alpha}$ for $0< \alpha <1$ on any compact subset of $\Theta_o \cup \partial_0 \Theta_o$ such that $f_o(\partial_0 \Theta_o) \subset L_o$. In particular, if $t_o \in \Theta_o \cup \partial_0 \Theta_o$, then $Df_o (t_o) = \displaystyle \lim Df_k (t_k)$.
\hfill$\Box$
\end{prop}

\se{Convergence of pseudo holomorphic polygons to gradient tree}
Assume $\omega = d\alpha$ on $M$, $\alpha|_{{\mathscr L}} = dh$ and $\omega|_{U_{L_i}} = d\alpha_i$, where $\alpha_i$ is the canonical 1-form on the symplectic neighborhood $U_{L_i}$ of $L_i$. Let $\partial_i \Theta$ denote the boundary component of $\Theta \cong \mathbb{D}_{[\scriptsize{\textsc{n}}]}$ between $\textsc{p}_{i-1}$ and $\textsc{p}_i$.\\
\begin{prop}
\label{md}
For a pseudo holomorphic disk $w: \Theta \rightarrow M$ such that $w(\partial_i \Theta) \subset \Lambda^\epsilon_{f_i}$, we have $\int_\Theta w^* \omega = O(\epsilon)$.
\end{prop}
{\bf Proof:} Let $\textsc{x}^\epsilon_i = w (\textsc{p}_i) \in \Lambda^\epsilon_{f_i} \cap \Lambda^\epsilon_{f_{i+1}}$ be the Lagrangian intersection point, where $w(\partial_i \Theta)$ and $w(\partial_{i+1} \Theta)$ meet. $\textsc{x}_i =\pi_i (\textsc{x}^\epsilon_i) = \pi_{i+1} (\textsc{x}^\epsilon_i)$ is a critical point of $f_{i, i+1}$ in $L_{i, i+1}$. Let $l'_i$ be the path connecting $\textsc{x}^\epsilon_i$ and $\textsc{x}_i$ such that $\pi_i (l'_i) = \textsc{x}_i$. $D_i = \bigcup_{t\in [0,1]} tw(\partial_i \Theta)$ is a disk with boundary $\partial D_i = w(\partial_i \Theta) \cup l'_i \cup (\pi_i \circ w(\partial_i \Theta)) \cup l'_{i-1}$. Since $\alpha - \alpha_i$ is closed 1-form and $\alpha_i|_{l'_i} =0$, $\alpha_i|_{\pi_i \circ w(\partial_i \Theta)} =0$. We have

\[
\int_{\partial_i \Theta} w^* (\alpha -\alpha_i) = \int_{\partial_i \Theta} w^* \pi_i^* \alpha + \int_{l'_i} \alpha - \int_{l'_{i-1}} \alpha,
\]
\[
= h(\textsc{x}_i) - h(\textsc{x}_{i-1}) + \int_{l'_i} \alpha - \int_{l'_{i-1}} \alpha.
\]
\[
\int_\Theta w^* \omega = \sum_{i=1}^{\scriptsize{\textsc{n}}} \int_{\partial_i \Theta} w^* \alpha = \sum_{i=1}^{\scriptsize{\textsc{n}}} \left(\int_{\partial_i \Theta} w^* \alpha_i +\int_{\partial_i \Theta} w^* (\alpha -\alpha_i)\right)
\]
$\hspace{.9in} \displaystyle = \epsilon\sum_{i=1}^{\scriptsize{\textsc{n}}} (f_i(\textsc{x}_i) - f_i(\textsc{x}_{i-1})) = O(\epsilon).$
\hfill$\Box$\\
\begin{prop}
\label{mc}
Assume that $(M,g,J,\omega)$ is of bounded local geometry. Then $\displaystyle \sup_{t\in \Theta} |Dw(t)| = O(\epsilon)$ for the pseudo holomorphic disk in proposition \ref{md}.
\end{prop}
{\bf Proof:} Suppose the proposition is not true. Then there exists a sequence $\epsilon_k \rightarrow 0$, $w_k: \Theta_k \rightarrow M$ and $t_k \in \Theta_k$ such that

\[
\lim_{k\rightarrow +\infty} \frac{\epsilon_k}{\epsilon'_k} =0,\mbox{ where } \epsilon'_k = |Dw_k(t_k)| = \sup_{t\in \Theta_k} |Dw_k(t)|.
\]

Assume $\displaystyle \lim \epsilon'_k =0$. Since ${\mathscr L}$ is compact, without loss of generality, we may assume $w_k(t_k) \rightarrow \textsc{x} \in {\mathscr L} \subset M$. Corollary \ref{bd} implies that the sequence $(M, w_k(t_k), g/(\epsilon'_k)^2)$ (pointed) converges in the sense of Cheeger-Gromov to $(T_{\scriptsize{\textsc{x}}} M, 0, g_{\scriptsize{\textsc{x}}})$. Under the identification $T_{\scriptsize{\textsc{x}}} M \cong T_{\scriptsize{\textsc{x}}} {\mathscr L} \oplus T^*_{\scriptsize{\textsc{x}}} {\mathscr L}$, the corresponding limit of ${\mathscr L}$ is $p + T_{\scriptsize{\textsc{x}}} {\mathscr L}$ for certain $p\in T^*_{\scriptsize{\textsc{x}}} {\mathscr L}$. Since $\displaystyle \lim \frac{\epsilon_k}{\epsilon'_k} =0$, $\Lambda^{\epsilon_k}_{f_i}$ converge to $p + T_{\scriptsize{\textsc{x}}} L_i \cong p + T_{\scriptsize{\textsc{x}}} {\mathscr L} \subset T_{\scriptsize{\textsc{x}}} M$ if $\textsc{x} \in L_i$. Clearly, when $\textsc{x} \not\in L_i$, $\Lambda^{\epsilon_k}_{f_i}$ will disappear at the limit. Hence, the limit of $\{\Lambda^{\epsilon_k}_{f_i}\}$ collectively lies in $p + T_{\scriptsize{\textsc{x}}} {\mathscr L} \subset T_{\scriptsize{\textsc{x}}} M$.\\

Proposition \ref{ab} implies that by possibly taking subsequence, the sequence $(\Theta_k, t_k, g_{\Theta_k})$ (pointed) converge in the sense of Cheeger-Gromov to $(\Theta_o, t_o, g_{\Theta_o})$, and $\partial_i \Theta_o = \displaystyle \lim \partial_i \Theta_k$. Since $|D w_k| \leq 1$ for $w_k: (\Theta_k, t_k, g_{\Theta_k}) \rightarrow (M, w_k(t_k), g/(\epsilon'_k)^2)$, according to proposition \ref{ni}, by passing to subsequence, $w_k$ converge to $w_o: (\Theta_o, t_o, g_{\Theta_o}) \rightarrow (T_{\scriptsize{\textsc{x}}} M, 0, g_{\scriptsize{\textsc{x}}})$ such that $|D w_o| \leq |D w_o(t_o)| =1$ and $w_o (\partial \Theta_o) \subset p + T_{\scriptsize{\textsc{x}}} {\mathscr L}$.\\

It is convenient to identify $(T_{\scriptsize{\textsc{x}}} M \cong (p +T_{\scriptsize{\textsc{x}}} {\mathscr L}) \times T^*_{\scriptsize{\textsc{x}}} {\mathscr L}, p)$ with $(\mathbb{C}^n \cong \mathbb{R}^n \times i\mathbb{R}^n, 0)$. Then $w_o$ can be understood as a holomorphic map $(\Theta_o, \partial \Theta_o, t_o) \rightarrow (\mathbb{C}^n, \mathbb{R}^n, 0)$. According to proposition \ref{ed}, $w_o$ has to be a constant map, which contradict with the fact $|D w_o(t_o)| =1$.\\

Assume $\displaystyle \lim \epsilon'_k \geq 2c > 0$. Proposition \ref{bc} implies that by possibly taking subsequence, the sequence $(\Theta_k, t_k, (\epsilon'_k)^2 g_{\Theta_k})$ (pointed) converge in the sense of Cheeger-Gromov to $(\Theta_o, t_o, g_{\Theta_o})$. Since $|D w_k| \leq |D w_k (t_k)| = 1$ for $w_k: (\Theta_k, t_k, (\epsilon'_k)^2 g_{\Theta_k}) \rightarrow (M, w_k(t_k), g)$, if $z_o = \displaystyle \lim w_k(t_k)$ exists in $M$, according to proposition \ref{ni}, by possibly passing to subsequence, $w_k$ converge to $w_o: (\Theta_o, t_o, g_{\Theta_o}) \rightarrow (M, z_o, g)$ satisfying $|D w_o (t_o)| = 1$. In particular $\int_{\Theta_o} w_o^* \omega >0$. This contradicts with $\displaystyle \int_{\Theta_o} w_o^* \omega = \lim \int_{\Theta_k} w_k^* \omega = \lim O(\epsilon_k) = 0$.\\

If $\displaystyle \lim w_k(t_k)$ is not in the finite part of $M$, since $(M, g)$ is of uniformly bounded local geometry according to our assumption, by corollary \ref{be}, $(M, w_k(t_k), g)$ pointed converges to $(M_o, z_o, g_o)$. Similarly apply proposition \ref{ni}, by passing to subsequence, $w_k$ converge to $w_o: (\Theta_o, t_o, g_{\Theta_o}) \rightarrow (M_o, z_o, g_o)$ satisfying $|D w_o (t_o)| = 1$. In particular $\int_{\Theta_o} (w_o)^* \omega >0$. This contradicts with $\displaystyle \int_{\Theta_o} w_o^* \omega_o = \lim \int_{\Theta_k} w_k^* \omega = \lim O(\epsilon_k) = 0$.
\hfill$\Box$

\begin{re}
The technical condition ``$(M,g,J,\omega)$ is of uniformly bounded local geometry" can be omitted if one can show that $\displaystyle \lim w_k(t_k) = x_o$ is in the finite part of $M$ using maximal principle.
\end{re}

\begin{re}
The estimate $|Dw| = O(\epsilon)$ was also carried out in \cite{FO} for the special case of $M = T^*L$, where ${\mathscr L} = L$ is a Hausdorff manifold. The proof can be separated into 4 steps:\\

{\bf Step 1:} $\int_\Theta w^* \omega = O(\epsilon)$. (Lemma 9.2 in \cite{FO})\\

{\bf Step 2:} The $C^1$-estimate $\displaystyle \sup_{t\in \Theta} |Dw(t)| = O(\sqrt{\epsilon})$, consequently, the $C^0$-estimate $d(w(\Theta), L) = O(\sqrt{\epsilon})$. (Lemma 9.3 in \cite{FO})\\

{\bf Step 3:} The $C^0$-estimate $d(w(\Theta), L) = O(\epsilon)$ via maximum principle. (Proposition 9.4 in \cite{FO})\\

{\bf Step 4:} The $C^1$-estimate $\displaystyle \sup_{t\in \Theta} |Dw(t)| = O(\epsilon)$. (Proposition 9.7 in \cite{FO})\\

In \cite{FO}, the step 1 is not used anywhere except in the proof of the step 2. The step 2 in its precise form is not used anywhere. Only a rather tangential consequence of the step 2 (namely, $\displaystyle \sup_{t\in \Theta} |Dw(t)|$ is finite) is necessary for the proof of the step 4. In more general case when ${\mathscr L} =L$ and $M \not= T^*L$, if the method of \cite{FO} is to be used, the precise form of the step 2 would be needed to first show the disk is in effect inside the symplectic neighborhood of $L$ before the steps 3 and 4 can be performed. Of course such important function of the step 2 is unnecessary in \cite{FO}, where $M = T^*L$.\\

In our case, where ${\mathscr L}$ is non-Hausdorff, it is not clear wether the analogue of the $C^0$-estimate via maximum principle in step 3, which is of global nature, is possible. An alternative argument without maximal principle is desirable. The main point of our argument here is the realization that the much weaker result in step 1 is essentially enough to prove the $C^1$-estimate $|Dw| = O(\epsilon)$ already without any additional assumptions. Therefore eliminate the need for using maximum principle, and apply very nicely to our case when ${\mathscr L}$ is non-Hausdorff. (In our proof, we also need the general property of holomorphic polygon that $\displaystyle \sup_{t\in \Theta} |Dw(t)|$ is finite, which, for example, is a consequence of the more precise estimate of lemma 9.3 in \cite{FO}.)
\end{re}

Consider a sequence of pseudo holomorphic disks $w_k: \Theta_k \rightarrow (M, J)$ such that $w_k \in {\cal M}_J(M, \vec{\Lambda}^{\epsilon_k}, \vec{\textsc{x}}^{\epsilon_k})$ with $\lim \epsilon_k =0$. By possibly taking subsequence, we may assume that $\Gamma_k$ are topologically the same as a fixed graph $\Gamma$. Then for $\nu \in C^*(\Gamma)$, $\{i: \Theta_k (\nu) \cap \partial_i \Theta_k \not= \emptyset\}$ can be identified with $I(\nu)$ defined in section 2.3 and is independent of $k$.\\

If $\displaystyle \lim {\rm Diam}_g w_k(\Theta_k) =0$, then by possibly taking subsequence, $\displaystyle \lim w_k(\Theta_k)$ is a single point, which implies coincidence of critical points of certain $f_{ij}$ that would not occur in our generic situation. (Even when it occurs, it is a local situation that is very easy to understand.)\\

Assume $\displaystyle \lim {\rm Diam}_g w_k(\Theta_k) \not=0$. By possibly taking subsequence, one can find a reference point $\textsc{x}^\circ= \displaystyle \lim w_k(t^\circ_k) \in L$ that is not a critical point of any $f_{ij}$, where $t^\circ_k \in \Theta_k (\nu^\circ)$ for certain $\nu^\circ \in C^1 (\Gamma)$. According to proposition \ref{ac}, by possibly taking subsequence, we have the (pointed) Gromov-Hausdorff limit

\[
(\Gamma^\circ, \tau^\circ, \gamma^\circ) = \lim_{k \rightarrow +\infty} (\Theta_k, t^\circ_k, \epsilon_k^2 g_{\Theta_k}) = \lim_{k \rightarrow +\infty} (\Gamma_k, \tau^\circ_k, \epsilon_k \gamma_k).
\]

In this context, by possibly taking subsequence, one may assume that for $\nu\in C^*(\Gamma)$, the (pointed) Gromov-Hausdorff limit $\Gamma^\circ (\nu) = \displaystyle \lim_{k \rightarrow +\infty} (\Theta_k (\nu), t^\circ_k, \epsilon_k^2 g_{\Theta_k})$ exists (maybe empty) as a subset of $\Gamma^\circ$. Let $\Gamma_\circ$ be a subgraph of $\Gamma$ containing $\nu \in C^* (\Gamma)$ such that $\Gamma^\circ (\nu)$ is non-empty. Clearly, when $\Gamma^\circ (\nu)$ is non-empty, it is a strata of $\Gamma^\circ$. $\psi_\circ (\nu) := \Gamma^\circ (\nu)$ defines a surjective topological map $\psi_\circ: C^*(\Gamma_\circ) \rightarrow C^*(\Gamma^\circ)$. For $\nu \in C^*(\Gamma^\circ)$, we define $I(\nu)$ to be the union of all $I(\tilde{\nu})$ such that $\tilde{\nu} \in C^*(\Gamma_\circ)$ and $\psi_\circ (\tilde{\nu}) = \nu$.\\

Proposition \ref{mc} implies that $w_k: (\Theta_k, t^\circ_k, \epsilon_k^2 g_{\Theta_k}) \rightarrow (M, w_k(t^\circ_k), g)$ is $C^1$-bounded (uniformly with respect to $k$). By possibly taking subsequence, we may assume that $w_k$ restricted to $B_R^{\epsilon_k^2 g_{\Theta_k}} (\Theta_k, t^\circ_k)$ converges to a Lefschetz continuous map $\Upsilon^\circ_R: (\Gamma^\circ_R, \tau^\circ) \rightarrow (M, \textsc{x}^\circ)$ that can be extended to a Lefschetz continuous map $\Upsilon^\circ: (\Gamma^\circ, \tau^\circ) \rightarrow (M, \textsc{x}^\circ)$. We will need the following technical lemma\\
\begin{lm}
\label{mg}
Assume $\mathbb{C}^n$-valued (resp. $\mathbb{R}$-valued) function $u_o (x) = \lim u_k(x)$ for $x$ in a neighborhood of $x_o \in \mathbb{R}$ and $u'_o(x_o) \not= A$ (resp. ${\displaystyle \liminf_{x \rightarrow x_o}} \frac{u_o(x) - u_o(x_o)}{x -x_o} < A$). Then by possibly taking subsequence of $\{u_k\}$, there exists $c_1>0$ and $x_k \rightarrow x_o$ such that $|u'_k(x_k) - A| \geq c_1$ (resp. $u'_k(x_k) \leq A-c_1$).
\end{lm}
{\bf Proof:} $u'_o(x_o) \not= A$ (resp. ${\displaystyle \liminf_{x \rightarrow x_o}} \frac{u_o(x) - u_o(x_o)}{x -x_o} < A$) implies that there exists $c_1>0$ and $\tilde{x}_k \rightarrow x_o$ such that $|\frac{u_o(\tilde{x}_k) - u_o(x_o)}{\tilde{x}_k -x_o} -A| \geq 2c_1$ (resp. $\frac{u_o(\tilde{x}_k) - u_o(x_o)}{\tilde{x}_k -x_o} \leq A - 2c_1$). For any $k$, there exists $n_k$ such that $|\frac{u_{n_k}(\tilde{x}_k) - u_{n_k}(x_o)}{\tilde{x}_k -x_o} -A| \geq c_1$ (resp. $\frac{u_{n_k}(\tilde{x}_k) - u_{n_k}(x_o)}{\tilde{x}_k -x_o} \leq A - c_1$). Consequently, there exists $x_k$ between $x_o$ and $\tilde{x}_k$ such that $|u'_{n_k}(x_k) - A| \geq c_1$ (resp. $u'_{n_k}(x_k) \leq A - c_1$).
\hfill$\Box$

\begin{re}
\label{mga}
In applications of lemma \ref{mg}, $u_k$ are usually $C^\infty$-functions. For suitable choice of $A_{x_o}$, we usually prove that the assertion ``for suitable subsequence of $\{u_k\}$, there exists $c_1>0$ and $x_k \rightarrow x_o$ such that $|u'_k(x_k) - A_{x_o}| \geq c_1$" leads to contradiction. By lemma \ref{mg}, this contradiction imply $u'_o(x_o) = A_{x_o}$. It is straightforward to see that this contradiction also implies that $u'_k$ converges to $u'_o$ uniformly. Consequently $u_o$ is a $C^1$-function.
\end{re}

For $\nu \in C^1(\Gamma^\circ)$, there is a unique $\tilde{\nu} \in C^1(\Gamma)$ such that $\psi_\circ (\tilde{\nu}) = \nu$. Since $|I(\tilde{\nu})| =2$ for $\tilde{\nu} \in C^1(\Gamma)$, we may assume $I(\nu) = I(\tilde{\nu}) =\{i,j\}$. Then $\Upsilon^\circ (\nu)$ is a Lefschetz curve in $L_{ij}$.\\
\begin{prop}
\label{mb}
For $\tau_o \in \nu \in C^1(\Gamma^\circ)$ with $I(\nu) = I(\tilde{\nu}) =\{i,j\}$, $\displaystyle \left|\frac{d \Upsilon^\circ}{d\tau}\right|(\tau_o) := \liminf_{\tau \rightarrow \tau_o} \frac{|\Upsilon^\circ (\tau) - \Upsilon^\circ (\tau_o)|}{|\tau - \tau_o|} \geq |\nabla f_{ij} (\textsc{x}_o)|$, where $\textsc{x}_o = \Upsilon^\circ (\tau_o)$.
\end{prop}
{\bf Proof:} Without loss of generality, we may assume that $\nabla f_{ij} (\textsc{x}_o) \not=0$. Take $\tilde{t}_k \in \Theta_k (\tilde{\nu})$ such that $\tau_o = \lim \tilde{t}_k$ under the limit $(\Gamma^\circ, \tau^\circ, \gamma^\circ) = \lim (\Theta_k, t^\circ_k, \epsilon_k^2 g_{\Theta_k})$. Then $\textsc{x}_o = \Upsilon^\circ (\tau_o) = \lim w_k(\tilde{t}_k)$. The definition of $\Upsilon^\circ$ implies that

\[
\left\langle \Upsilon^\circ (\tau_o +x), \frac{\nabla f_{ij} (\textsc{x}_o)}{|\nabla f_{ij} (\textsc{x}_o)|} \right\rangle = \lim_{k\rightarrow +\infty} \left\langle w_k \left(\tilde{t}_k + \frac{x}{\epsilon_k}\right), \frac{\nabla f_{ij} (\textsc{x}_o)}{|\nabla f_{ij} (\textsc{x}_o)|} \right\rangle.
\]

Lemma \ref{mg} can be applied to these $\mathbb{R}$-valued functions of $x$ if the estimate of the proposition is not true. Then by possibly taking subsequence, there exists $c_1>0$ and $x_k \rightarrow 0$ that give rise to $t_k = \tilde{t}_k + \frac{x_k}{\epsilon_k} \in \Theta_k (\tilde{\nu})$ such that $\textsc{x}_o = \displaystyle \lim_{k \rightarrow +\infty} w_k(t_k)$ and for $k$ large,

\begin{equation}
\label{mba}
\frac{1}{\epsilon_k} \left\langle \frac{d w_k}{dt} (t_k), \nabla f_{ij} (\textsc{x}_o)\right\rangle \leq |\nabla f_{ij} (\textsc{x}_o)|^2 - c_1|\nabla f_{ij} (\textsc{x}_o)|.
\end{equation}

Proposition \ref{ab} implies that by possibly taking subsequence, the sequence $(\Theta_k, t_k, g_{\Theta_k})$ (pointed) converge in the sense of Cheeger-Gromov to $(\Theta_o, t_o, g_{\Theta_o})$. Since $\tau_o \in \nu$, there exists $c_2>0$ such that $[\tau_o - c_2, \tau_o + c_2] \subset \nu$. For $k$ large enough, $[t_k^\Re - \frac{c_2}{\epsilon_k}, t_k^\Re + \frac{c_2}{\epsilon_k}] \times [0,1] \subset \Theta_k (\tilde{\nu})$. Hence $(\Theta_o, t_o) \cong (\mathbb{R} \times [0,1], 0)$.\\

Corollary \ref{bd} implies that the sequence $(M, w_k(t_k), g/\epsilon_k^2)$ (pointed) converge in the sense of Cheeger-Gromov to $(T_{\scriptsize{\textsc{x}}_o} M, 0, g_{\scriptsize{\textsc{x}}_o})$. Under the identification $T_{\scriptsize{\textsc{x}}_o} M \cong T_{\scriptsize{\textsc{x}}_o} {\mathscr L} \oplus T^*_{\scriptsize{\textsc{x}}_o} {\mathscr L}$, by possibly taking subsequence, the corresponding limit of $\Lambda^{\epsilon_k}_{f_i}$ (resp. $\Lambda^{\epsilon_k}_{f_j}$) is $p_i + T_{\scriptsize{\textsc{x}}_o} {\mathscr L}$ (resp. $p_j + T_{\scriptsize{\textsc{x}}_o} {\mathscr L}$) for certain $p_i \in T^*_{\scriptsize{\textsc{x}}_o} {\mathscr L}$ (resp. $p_j \in T^*_{\scriptsize{\textsc{x}}_o} {\mathscr L}$). (Proposition \ref{mc} is implicitly used here to guarantee the finiteness of $p_i,p_j$.) Furthermore, $p_i,p_j \in T^*_{\scriptsize{\textsc{x}}_o} {\mathscr L}$ satisfy

\[
p_j -p_i = df_{ij}(\textsc{x}_o) + \lim_{k\rightarrow +\infty} \frac{d h_{ij}(w_k(t_k))}{\epsilon_k}.
\]

Proposition \ref{mc} implies that $|D w_k| \leq C$ for $w_k: (\Theta_k, t_k, g_{\Theta_k}) \rightarrow (M, w_k(t_k), g/\epsilon_k^2)$. According to proposition \ref{ni}, by possibly passing to subsequence, $w_k$ converge to $w_o: (\Theta_o, t_o, g_{\Theta_o}) \rightarrow (T_{\scriptsize{\textsc{x}}_o} M, 0, g_{\scriptsize{\textsc{x}}_o})$ such that $|D w_o| \leq C$ and $w_o (\partial_i \Theta_o) \subset p_i + T_{\scriptsize{\textsc{x}}_o} {\mathscr L}$ (resp. $w_o (\partial_j \Theta_o) \subset p_j + T_{\scriptsize{\textsc{x}}_o} {\mathscr L}$). Under suitable local coordinates, the limiting process can be expressed as

\[
w_o (t) = \lim_{k \rightarrow +\infty} \frac{1}{\epsilon_k} (w_k(t_k + t) - w_k(t_k)).
\]

Hence (\ref{mba}) implies that

\begin{equation}
\label{mbb}
\left\langle \frac{d w_o}{dt} (t_o), \nabla f_{ij} (\textsc{x}_o)\right\rangle \leq |\nabla f_{ij} (\textsc{x}_o)|^2 - c_1|\nabla f_{ij} (\textsc{x}_o)|.
\end{equation}

$w_o$ can be understood as a holomorphic map with bounded derivative from the strip $\{t\in \mathbb{C}: 0\leq t^\Im \leq 1\}$ to $\mathbb{C}^n$, with the 2 boundaries of the strip mapped to shifts of $\mathbb{R}^n \subset \mathbb{C}^n$. According to proposition \ref{ec}, $w_o$ must be linear, hence

\[
\frac{d w_o}{dt} (t_o) = \tilde{v}_o = \nabla f_{ij}(\textsc{x}_o) + \lim_{k\rightarrow +\infty} \frac{\nabla h_{ij}(w_k(t_k))}{\epsilon_k}.
\]

Since $\langle \nabla h_{ij}(w_k(t_k)), \nabla f_{ij}(w_k(t_k))\rangle \geq 0$, we have

\[
\left\langle \frac{d w_o}{dt} (t_o), \nabla f_{ij} (\textsc{x}_o)\right\rangle \geq |\nabla f_{ij} (\textsc{x}_o)|^2,
\]

which is a contradiction to (\ref{mbb}).
\hfill$\Box$\\

Without loss of generality, we may assume $h_{ij} \geq 0$. Let $B_{ij,c} = \{h_{ij}=c\}$ for $c>0$ and $B_{ij,0} = \displaystyle \lim_{c \rightarrow 0} B_{ij,c} = \partial L_{ij}$. Then $\nabla f_{ij}$ is pointing outward on $B_{ij,0}$ and $(\nabla \hat{f}_{ij}, \nabla h_{ij}) > 0$ outside of $L_{ij}$.

\begin{prop}
\label{ma}
For any $\tau^-, \tau^+ \in \nu \in C^1(\Gamma^\circ)$ with $I(\nu) = I(\tilde{\nu}) =\{i,j\}$, $\Upsilon^\circ ([\tau^-, \tau^+]) \setminus \Upsilon^\circ (\{\tau^-, \tau^+\})$ is in the interior of $L_{ij}$.
\end{prop}
{\bf Proof:} Assume $\textsc{x}_o = \Upsilon^\circ (\tau_o) \in \Upsilon^\circ ([\tau^-, \tau^+]) \setminus \Upsilon^\circ (\{\tau^-, \tau^+\})$ for some $\tau_o \in (\tau^-, \tau^+)$. One can find  $t_k \in \Theta_k (\tilde{\nu})$ such that $\textsc{x}_o = \displaystyle \lim w_k(t_k)$ and $\frac{d u_k(t^\Re)}{dt^\Re} \not=0$ for $t^\Re \in [0, \gamma_{\tilde{\nu}}]$, where $u_k(t^\Re) =w_k (t^\Re + it_k^\Im)$. (Recall that $\Theta_k (\tilde{\nu}) \cong \{t\in \mathbb{C}: (t^\Re, t^\Im) \in [0, \gamma_{\tilde{\nu}}] \times [0,1]\}$, and $\{t\in \Theta_k (\tilde{\nu}): \frac{d w_k}{dt}=0\}$ is discrete in $\Theta_k (\tilde{\nu})$.)\\

Assume that $\textsc{x}_o$ is also in $B_{ij,0} = \partial L_{ij}$. By suitably modifying $t_k^\Re \in [0, \gamma_{\tilde{\nu}}]$, one can ensure that $\textsc{x}_o = \displaystyle \lim w_k(t_k)$ and $v_o = \displaystyle \lim v_k$ is tangent to $B_{ij,0}$, where $v_k$ is the unit tangent of $u_k([0, \gamma_{\tilde{\nu}}])$ at $w_k (t_k) = u_k (t_k^\Re)$.\\

If not so, then there exist $c_1,c_2 >0$ such that the connected component of $[0, \gamma_{\tilde{\nu}}] \cap u_k^{-1} (B_{\scriptsize{\textsc{x}}_o}(c_1))$ containing $t_k^\Re$ is $[t_k^-, t_k^+]$, and $|\frac{d\tilde{h}(u_k(t^\Re))}{dt^\Re}|\geq c_2|\frac{du_k(t^\Re)}{dt^\Re}|$ for $t^\Re \in [t_k^-, t_k^+]$, where $\tilde{h}$ is a defining function of $L_{ij} = \{ \tilde{h} \leq 0\}$ such that $\tilde{h}$ grow linearly along the normal direction of $B_{ij,0}$. Since $\frac{du_k (t^\Re)}{dt^\Re} \not=0$ for $t^\Re \in [0, \gamma_{\tilde{\nu}}]$, without loss of generality, we may assume $\frac{d\tilde{h}(u_k(t^\Re))}{dt^\Re}>0$ for $t^\Re \in [0, \gamma_{\tilde{\nu}}]$. Consequently

\[
\tilde{h}(u_k(t_k^+)) - \tilde{h}(w_k(t_k)) = \int^{t_k^+}_{t_k^\Re} \frac{d\tilde{h}(u_k(t^\Re))}{dt^\Re} dt^\Re
\]
\[
\geq c_2\int^{t_k^+}_{t_k^\Re} \left|\frac{du_k(t^\Re)}{dt^\Re}\right|dt^\Re \geq c_2 {\rm Dist} (w_k(t_k), u_k(t_k^+)) = c_2 c_1.
\]

Since $\displaystyle \lim \tilde{h}(w_k(t_k)) =0$, $\displaystyle \lim \tilde{h}(u_k(t_k^+))$ has to be positive, which contradicts with the fact $\Upsilon^\circ ([\tau^-, \tau^+]) \subset \Upsilon^\circ (\nu) \subset L_{ij}$. Hence $v_o$ is tangent to $B_{ij,0}$.\\

Notice that $v_o = \tilde{v}_o/|\tilde{v}_o|$ for $\tilde{v}_o$ from the proof of proposition \ref{mb}. Since $\nabla h_{ij}(z_k)$ is normal to $B_{ij,c}$ for $c = h_{ij} (z_k)$, $\nabla f_{ij}$ is not tangent to $B_{ij,0}$ and is in the increasing direction of $h_{ij}$, $\tilde{v}_o$ is not tangent to $B_{ij,0}$, which is a contradiction. Therefore, $\Upsilon^\circ ([\tau^-, \tau^+]) \setminus \Upsilon^\circ (\{\tau^-, \tau^+\})$ is in the interior of $L_{ij}$.
\hfill$\Box$\\
\begin{prop}
\label{me}
$\Upsilon^\circ (\nu)$ is a gradient flow line for $f_{ij}$ on $\nu \in C^1(\Gamma^\circ)$ with $I(\nu) = \{i,j\}$. Furthermore, $w_k: (\Theta_k, t^\circ_k, \epsilon_k^2 g_{\Theta_k}) \rightarrow (M, w_k(t^\circ_k), g)$ converge to $\Upsilon^\circ: (\Gamma^\circ, \tau^\circ) \rightarrow (M, \textsc{x}^\circ)$ uniformly in $C^1$ over any compact subset in the smooth part of $\Gamma^\circ$ in the sense of Cheeger-Gromov.
\end{prop}
{\bf Proof:} Take $\tilde{t}_k \in \Theta_k (\tilde{\nu})$ such that $\tau_o = \lim \tilde{t}_k$ under the limit $(\Gamma^\circ, \tau^\circ, \gamma^\circ) = \lim (\Theta_k, t^\circ_k, \epsilon_k^2 g_{\Theta_k})$. Then $\textsc{x}_o = \Upsilon^\circ (\tau_o) = \lim w_k(\tilde{t}_k)$. The definition of $\Upsilon^\circ$ implies that

\[
\Upsilon^\circ (\tau_o +x) = \lim_{k\rightarrow +\infty} w_k \left(\tilde{t}_k + \frac{x}{\epsilon_k}\right).
\]

Lemma \ref{mg} can be applied to these functions of $x$ if the proposition is not true. Then by possibly taking subsequence, there exists $c_1>0$ and $x_k \rightarrow 0$ that give rise to $t_k = \tilde{t}_k + \frac{x_k}{\epsilon_k} \in \Theta_k (\tilde{\nu})$ such that $\textsc{x}_o = \displaystyle \lim w_k(t_k)$ and for $k$ large,

\begin{equation}
\label{mea}
\left| \frac{1}{\epsilon_k} \frac{d w_k}{dt} (t_k) - \nabla f_{ij} (\textsc{x}_o)\right| \geq c_1.
\end{equation}

Propositions \ref{mb} and \ref{ma} together imply that $\textsc{x}_o$ is in the interior of $L_{ij}$. Through the same arguments as in the proof of proposition \ref{mb} using Cheeger-Gromov convergence, we arrive at $w_o: (\Theta_o, t_o, g_{\Theta_o}) \rightarrow (T_{\scriptsize{\textsc{x}}_o} M, 0, g_{\scriptsize{\textsc{x}}_o})$ such that $|D w_o| \leq C$ and $w_o (\partial_i \Theta_o) \subset p_i + T_{\scriptsize{\textsc{x}}_o} {\mathscr L}$ (resp. $w_o (\partial_j \Theta_o) \subset p_j + T_{\scriptsize{\textsc{x}}_o} {\mathscr L}$), where $p_j -p_i = df_{ij}(\textsc{x}_o)$. Under suitable local coordinates, the limiting process can be expressed as

\[
w_o (t) = \lim_{k \rightarrow +\infty} \frac{1}{\epsilon_k} (w_k(t_k + t) - w_k(t_k)).
\]

Hence (\ref{mea}) implies that

\begin{equation}
\label{meb}
\left| \frac{d w_o}{dt} (0) - \nabla f_{ij} (\textsc{x}_o)\right| \geq c_1.
\end{equation}

$w_o$ can be understood as a holomorphic map with bounded derivative from the strip $\{t\in \mathbb{C}: 0\leq t^\Im \leq 1\}$ to $\mathbb{C}^n$, with the 2 boundaries of the strip mapped to shifts of $\mathbb{R}^n \subset \mathbb{C}^n$. According to proposition \ref{ec}, $w_o$ must be linear, hence

\[
\frac{d w_o}{dt} = \nabla f_{ij}(\textsc{x}_o),
\]

which is a contradiction to (\ref{meb}) that implies the first part of the proposition. According to remark \ref{mga}, such contradiction also implies the second part of the proposition.
\hfill$\Box$

\begin{re}
Proposition \ref{mb} is a weaker version of proposition \ref{me} and the proofs of the 2 propositions are similar. But proposition \ref{mb} is necessary in combination with proposition \ref{ma} to show that $\Upsilon^\circ (\nu)$ is in the interior of $L_{ij}$ that is needed for the proof of proposition \ref{me}. The proof of proposition \ref{ma} can be simplified if proposition \ref{mb} is used. We provide a proof of proposition \ref{ma} without using proposition \ref{mb} to show their independence.
\end{re}

\begin{co}
$\Gamma^\circ$ is a complete metric graph. $\Upsilon^\circ (\nu)$ contains no critical point of $f_{ij}$ for $i,j \in I(\nu)$, $\nu \in C^* (\Gamma^\circ)$. For an exterior leg $\nu \in C^1_{\rm ext}$, $\Upsilon^\circ (\tau)$ approaches a critical point of $f_{ij}$ when $\tau \in \nu$ approaches an exterior end of $\nu$.
\end{co}
{\bf Proof:} This is a simple corollary of proposition \ref{me} and the fact that a non-critical point take infinite time to reach a critical point along gradient flow. Since $\tau^\circ \in\Gamma^\circ$ is not a critical point for any $f_{ij}$, any point in $\Gamma^\circ$ has finite $\tau$-distance from $\tau^\circ$, $\Upsilon^\circ (\nu)$ contains no critical point of $f_{ij}$ for $i,j \in I(\nu)$, $\nu \in C^* (\Gamma^\circ)$.\\

Since $\Gamma_k$ is complete, exterior ends of $\Gamma^\circ_R$ have distance $R$ from $\tau^\circ$. Hence $\Gamma^\circ$ is complete. A gradient flow line always have to end somewhere, the only way for it to end in infinite time is to end at a critical point. Hence the last claim of the corollary.
\hfill$\Box$\\

If $\Gamma_\circ = \Gamma$ and $\lim w_k (t_k) = \textsc{x}_i$ for any $t_k \in \Theta_k (\nu_i)$ satisfying $\lim {\rm Dist}_{\epsilon_k^2 g_{\Theta_k}} (t_k, t_k^\circ) =+\infty$, where $\nu_i \in C_{\rm ext}^1 (\Gamma)$ leads to $\textsc{p}_i$, then the gradient tree $\Upsilon^\circ: \Gamma^\circ \rightarrow M$ is called the global limit of the pseudo holomorphic polygons $\{w_k\}$. Otherwise, $\Upsilon^\circ$ is called a local limit of $\{w_k\}$.\\

When $\Upsilon^\circ$ is not the global limit of $\{w_k\}$, by possibly taking subsequence, one can find $\textsc{x}^* = \displaystyle \lim w_k(t^*_k) \in {\mathscr L}$ that is not a critical point of any $f_{ij}$, where $t^*_k \in \Theta_k (\nu^*)$ for certain $\nu^* \in C^1 (\Gamma)$, such that $\lim {\rm Dist}_{\epsilon_k^2 g_{\Theta_k}} (t^*_k, t_k^\circ) =+\infty$. One may similarly get $\Gamma^*$, the subgraph $\Gamma_* \subset \Gamma$ and the limit gradient tree $\Upsilon^*: \Gamma^* \rightarrow M$ as another local limit of $\{w_k\}$. It is easy to see that $C^0_{\rm int} (\Gamma_*) \cap C^0_{\rm int} (\Gamma_\circ) = \emptyset$.\\

Such process can be continued inductively to get a sequence of gradient trees $\Upsilon^i: \Gamma^i \rightarrow M$ for $i = 0,1, \cdots$ as local limits of $\{w_k\}$, with $\Upsilon^0 = \Upsilon^\circ$. If this process terminate at $i=\textsc{k}$, then the collection $\vec{\Upsilon} = \{\Upsilon^0, \cdots, \Upsilon^{\scriptsize{\textsc{k}}}\}$ of gradient trees is called the global limit of $\{w_k\}$.\\

The finiteness of $\textsc{k}$ can be proved by either area method or index method as we explained in the discussion of the global Gromov-Hausdorff limit at the end of section 3. For the area method, we need

\begin{prop}
\label{ml}
There exists $c_0>0$ such that for any pointed limit $w_k: (\Theta_k, t^\circ_k, \epsilon_k^2 g_{\Theta_k}) \rightarrow (M, w_k(t^\circ_k), g)$ to $\Upsilon^\circ: (\Gamma^\circ, \tau^\circ, \gamma^\circ) \rightarrow (M, \textsc{x}^\circ, g)$, we have $\int_{B_R^{\epsilon_k^2 g_{\Theta_k}} (\Theta_k, t_k)} w_k^* \omega \geq \epsilon_k c_0$ for $R>0$ and $k$ large enough.
\end{prop}
{\bf Proof:} If the proposition is not true, then there exists a sequence of pseudo holomorphic $w_k: \Theta_k \rightarrow (M,J)$ with marked points $t^\circ_k \in \Theta_k$ such that $\textsc{x}^\circ = \lim w_k(t^\circ_k)$ is not a critical point for any $f_{ij}$, and $\lim \frac{1}{\epsilon_k} \int_{B_R^{\epsilon_k^2 g_{\Theta_k}} (\Theta_k, t_k)} w_k^* \omega =0$ for any $R>0$. Let $\Upsilon^\circ: (\Gamma^\circ, \tau^\circ, \gamma^\circ) \rightarrow (M, \textsc{x}^\circ, g)$ be the limiting gradient tree. Since $\nabla f_{ij} (\textsc{x}^\circ) \not=0$, one may find $R>0$ and $c_1>0$ such that $|\nabla f_{ij} (\Upsilon^\circ (\tau))| \geq 2c_1$ for $\tau \in B_R^{\gamma^\circ} (\Gamma^\circ, \tau^\circ)$. By proposition \ref{me}, for $\epsilon_k$ small, $\frac{1}{\epsilon_k} |Dw_k (t)| \geq c_1$ for $t \in B_R^{\epsilon_k^2 g_{\Theta_k}} (\Theta_k, t_k)$. Hence

\[
\frac{1}{\epsilon_k} \int_{B_R^{\epsilon_k^2 g_{\Theta_k}} (\Theta_k, t_k)} w_k^* \omega \geq \epsilon_k c_1^2 {\rm Vol}_{g_{\Theta_k}} (B_{R/\epsilon_k}^{g_{\Theta_k}} (\Theta_k, t_k)) \sim \int_{ B_R^{\gamma^\circ} (\Gamma^\circ, \tau^\circ)} d\gamma^\circ >0,
\]

which is a contradiction.
\hfill$\Box$\\

For a limiting gradient tree $\Upsilon^l$, assume $\nu \in C^1_{\rm ext} (\Gamma_l)$ with $I(\nu) =\{i,j\}$, then $\Upsilon^l (\nu)$ is a gradient line of $f_{ij}$. We call the corresponding external vertex of $\Upsilon^l$, (which is a critical point of $f_{ij}$,) the beginning vertex $\textsc{x}'_l$ or the ending vertex $\textsc{x}''_l$ according to the direction of $\nu$ determined by $\nabla f_{ij}$. We may choose the cylindrical complex coordinate $t$ on $\Theta_k (\nu)$ so that the increasing direction of $t^\Re$ is in accordance with this direction of $\nu$.

\begin{prop}
\label{mm}
For a limiting gradient tree $\Upsilon^l$ with $\nu \in C^1_{\rm ext} (\Gamma_l)$, if $\textsc{x}'_l$ (resp. $\textsc{x}''_l$) is not one of $\textsc{x}_i$, then there exists another limiting gradient tree $\Upsilon^{l'}$ with $\nu \in C^1_{\rm ext} (\Gamma_{l'})$ such that $\textsc{x}'_l = \textsc{x}''_{l'}$ (resp. $\textsc{x}''_l = \textsc{x}'_{l'}$).
\end{prop}
{\bf Proof:} Let $t^l_k \in \Theta_k (\nu)$ be the marked point to get the limit $\Upsilon^l$. For a small $c_1>0$ and a large $R>0$, there exists $t^{l'}_k \in \Theta_k (\nu)$ such that $\lim \epsilon_k (t^{l'}_k - t^l_k)^\Re = +\infty$ and ${\rm Dist}_g (w_k (t), \textsc{x}'_l) \leq {\rm Dist}_g (w_k (t^{l'}_k), \textsc{x}'_l) =c_1$ for $t \in \Theta_k (\nu)$ satisfying $(t^l_k)^\Re + \frac{R}{\epsilon_k} \leq t^\Re \leq (t^{l'}_k)^\Re$. Take $t^{l'}_k$ as the marked point, since $c_1>0$ is small, we get the limit $\Upsilon^{l'}$ satisfying $\nu \in C^1_{\rm ext} (\Gamma_{l'})$ and $\textsc{x}'_l = \textsc{x}''_{l'}$.
\hfill$\Box$\\
\begin{prop}
\label{mn}
The global limit $\vec{\Upsilon} = \{\Upsilon^0, \cdots, \Upsilon^{\scriptsize{\textsc{k}}}\}$ of $\{w_k\}$ exists with finite $\textsc{k}$, such that the components connect through there external vertices to form a tree. Furthermore, $\textsc{k} =0$ if ${\rm Ind} (\vec{f}, \vec{\textsc{x}}) =0$.
\end{prop}
{\bf Proof:} The finiteness of $\textsc{k}$ can be proved via either the area method or the index method. For the area method, one start with the total area bound from proposition \ref{md}, then proposition \ref{ml} implies that $\textsc{k}$ is finite.\\

For the index method, notice that $C^0_{\rm int} (\Gamma^i) \cap C^0_{\rm int} (\Gamma^j) = \emptyset$ for $i \not= j$. There are only finitely many $\Upsilon^l$ such that $C^0_{\rm int} (\Gamma^l) \not= \emptyset$. When $C^0_{\rm int} (\Gamma^l) = \emptyset$, there exists $\nu \in C^1(\Gamma)$ with $I(\nu) =\{i,j\}$ such that $\Gamma_l = \nu$ and $\Upsilon^l$ is a gradient line of $f_{ij}$. $\Upsilon^l$ has exactly 2 vertices: the beginning vertex $\textsc{x}'_l$ and the ending vertex $\textsc{x}''_l$, both of which are critical points of $f_{ij}$ so that $\mu (\textsc{x}''_l) +1 \leq \mu (\textsc{x}'_l)$. Proposition \ref{mm} then implies that such $\Upsilon^l$ form a finite connected chain. Consequently, $\textsc{k}$ is finite.\\

With $\textsc{k}$ being finite, proposition \ref{mm} implies that the components of $\vec{\Upsilon}$ connect through there external vertices to form a tree. Furthermore, it is straightforward to check that

\[
{\rm Ind} (\vec{f}, \vec{\textsc{x}}) + \textsc{k} = \sum_{l=0}^{\scriptsize{\textsc{k}}} {\rm Ind} (\Upsilon^l).
\]

According to our generic assumption, gradient tree $\Upsilon$ with negative virtual moduli dimension ${\rm Ind} (\Upsilon)$ does not exist. Consequently, $\textsc{k} =0$ if ${\rm Ind} (\vec{f}, \vec{\textsc{x}}) =0$.
\hfill$\Box$\\

\begin{theorem}
\label{mo}
For a sequence $w_k \in {\cal M}_J(M, \vec{\Lambda}^{\epsilon_k}, \vec{\textsc{x}}^{\epsilon_k})$ with ${\rm Ind} (\vec{f}, \vec{\textsc{x}}) =0$ and $\lim \epsilon_k =0$. By possibly taking subsequence, $w_k: (\Theta_k, \epsilon_k^2 g_{\Theta_k}) \rightarrow (M, g)$ globally converge to a rigid gradient tree $\Upsilon: (\Gamma, \gamma) \rightarrow ({\mathscr L},g|_{\mathscr L})$.
\hfill$\Box$
\end{theorem}

\begin{re}
\label{moa}
Theorem \ref{mo} reduces all our discussion to within $M_\Upsilon$ as a neighborhood of a rigid gradient tree $\Upsilon \subset L_\Upsilon \subset M_\Upsilon$. (See page \pageref{bj} for notations.) With the Lagrangian fibration $\pi_\Upsilon: M_\Upsilon \rightarrow L_\Upsilon$, $M_\Upsilon$ is a Darboux-Weinstein neighborhood of $L_\Upsilon$ with Lagrangian fibres being orthogonal to $L_\Upsilon$. This exactly is the situation of \cite{FO}. In principle, we may now apply results in \cite{FO} to prove theorem \ref{mk}. The discussion in the rest of this paper will be carried out for the Hausdorff Lagrangian submanifold $L = L_\Upsilon \subset {\mathscr L}$ in $M = M_\Upsilon$.
\end{re}

\se{The approximate solution}
\stepcounter{subsection}
{\bf \S \thesubsection\ The local model of the approximate solution:} Let $\Theta \cong \mathbb{D}_{[3]}$. By Riemann mapping theorem, there are unique maps $z_i: \Theta \rightarrow \mathbb{C}$ for $i=1,2,3$ such that

\[
z^\Im_1 (\partial_1 \Theta) = z^\Im_1 (\partial_2 \Theta) = 0,\ z^\Im_1 (\partial_3 \Theta) = \pi,\ z_1 (\{\textsc{p}_1, \textsc{p}_2, \textsc{p}_3\})= \{0, +\infty, -\infty\}.
\]
\[
z^\Im_2 (\partial_2 \Theta) = z^\Im_2 (\partial_3 \Theta) = 0,\ z^\Im_2 (\partial_1 \Theta) = \pi,\ z_2 (\{\textsc{p}_1, \textsc{p}_2, \textsc{p}_3\}) = \{-\infty, 0, +\infty\}.
\]
\[
z^\Im_3 (\partial_1 \Theta) = z^\Im_3 (\partial_3 \Theta) = 0,\ z^\Im_3 (\partial_2 \Theta) = \pi,\ z_3 (\{\textsc{p}_1, \textsc{p}_2, \textsc{p}_3\}) = \{+\infty, -\infty, 0\}.
\]

$z_1,z_2,z_3$ satisfy the relations:

\[
z_2 = \log \frac{e^{z_1} - 1}{e^{z_1}},\ z_3 = \log \frac{e^{z_2} - 1}{e^{z_2}},\ z_1 = \log \frac{e^{z_3} - 1}{e^{z_3}}.
\]
\[
z_3 = -\log (1 - e^{z_1}),\ z_1 = -\log (1 - e^{z_2}),\ z_2 = -\log (1 - e^{z_3}).
\]

Consequently, $z_1 + z_2 + z_3 =\pi i$.\\

Let $q + ip \cong (q,p)$ be the coordinate of $\mathbb{C}^n \cong \mathbb{R}^n + i\mathbb{R}^n$.

\begin{equation}
\label{ona}
u = \frac{z_2p_1 + z_3p_2 +z_1p_3}{\pi}: \Theta \rightarrow \mathbb{C}^n
\end{equation}

satisfies $u (\partial_i \Theta) \subset \Lambda_i =\{p=p_i\}$.\\

Using $t = \frac{z_{i-1}}{\pi}$ as cylindrical coordinate on $\Theta_i$, we have

\[
u(t) = ip_i + t(p_{i+1} - p_i) + \frac{\log (1 - e^{-\pi t})}{\pi}(p_{i-1} - p_i)
\]

On $\Theta_i$, let $\tilde{u}_\epsilon (t) = \tilde{u}^i_\epsilon (t) + \check{u}^i_\epsilon (t)$, where

\[
\tilde{u}^i_\epsilon (t) = ip_i + t(p_{i+1} - p_i),\ \check{u}^i_\epsilon (t) =  \eta(t^\Re) \frac{\log (1 - e^{-\pi t})}{\pi}(p_{i-1} - p_i),
\]

$\eta(t^\Re) = 0$ for $t^\Re \geq 2\epsilon^{-\alpha}$ and $\eta(t^\Re) = 1$ for $t^\Re \leq \epsilon^{-\alpha}$, $0 < \alpha <1$. $\tilde{u}^i_\epsilon (t)$ is the unique holomorphic linear strip between $\Lambda_i$ and $\Lambda_{i+1}$ passing through $ip_i$.\\
\begin{prop}
\label{on}
Up to addition of a real constant, there exists a unique holomorphic map $u: \Theta \rightarrow \mathbb{C}^n$ with bounded derivative that satisfies $u (\partial_i \Theta) \subset \Lambda_i$.
\end{prop}
{\bf Proof:} The existence is implied by the explicit formula (\ref{ona}). For the uniqueness, assume $u_1,u_2$ both satisfy the proposition. Then we have $\delta u = u_1-u_2: \Theta \rightarrow \mathbb{C}^n$ with bounded derivative satisfying $\delta u (\partial \Theta) \subset \mathbb{R}^n$. According to proposition \ref{ed}, $\delta u$ is a real constant map.
\hfill$\Box$\\

\stepcounter{subsection}
{\bf \S \thesubsection\ The non-exceptional case:} The approximate solution can be constructed from a gradient tree $\Upsilon: \Gamma \rightarrow M$ in 2 stages: first for the smooth part of $\Upsilon$ including the ends, second for the 3-valent vertices of $\Upsilon$. This construction is essentially the same as the one in \cite{FO}, with certain differences in detail. (Special construction is needed near an exceptional vertex of $\Upsilon$.)\\

For any $\textsc{x} \in L$, it is convenient to choose complex coordinate $z = z^\Re + iz^\Im$ in a neighborhood $U_{\scriptsize{\textsc{x}}}$ of $\textsc{x}$ in $M$ such that $z(\textsc{x})=0$, $L = \{z^\Im=0\}$, $\omega = dz^\Re \wedge dz^\Im$, $z^\Re = \pi(z)$ is constant on fibres of the transverse Lagrangian fibration $\pi$ near $L$ and the almost complex structure $|J - J_0| = O(|z|)$, where $J_0$ is the complex structure determined by $z$. Under such coordinate, $z^\Re$ can be regarded as coordinate of $L$ and $z \rightarrow z^\Im d z^\Re$ induce the identification of $U_{\scriptsize{\textsc{x}}}$ with a neighborhood of $0_{\scriptsize{\textsc{x}}}$ in $T^* L$. Under such coordinate, $\Lambda^{\epsilon}_{f_i}$ can be locally described as the graph of $z^\Im = \epsilon \frac{\partial f_i}{\partial z^\Re} (z^\Re)$ that corresponds to the 1-form $\epsilon df_i (z^\Re)$ on $L$.\\

For $\nu \in C^1(\Gamma)$ such that $I(\nu) = \{i,j\}$ and $t\in \Theta(\nu)$, define $\tilde{w}^{\nu} (t) = (z^\Re (t), z^\Im (t))$, where

\[
z^\Re (t) = \pi \circ \tilde{w}^{\nu} (t) = \Upsilon (\epsilon t^\Re),\ z^\Im (t) = \epsilon \left( \left(1 - t^\Im \right) \frac{\partial f_i}{\partial z^\Re} + t^\Im \frac{\partial f_j}{\partial z^\Re}\right) (z^\Re (t)).
\]

$\tilde{w}^{\nu}: \Theta(\nu) \rightarrow M$ satisfies $\tilde{w}^{\nu} (\partial_i \Theta(\nu)) \subset \Lambda^{\epsilon}_{f_i}$ and $\tilde{w}^{\nu} (\partial_j \Theta(\nu)) \subset \Lambda^{\epsilon}_{f_j}$. Through 1-form on $L$, $\tilde{w}^{\nu} (t)$ can also be expressed as $\epsilon (1 -t^\Im) df_i(z^\Re (t)) + \epsilon t^\Im df_j(z^\Re (t))$.\\
\begin{prop}
\label{oj}
For $t \in \Theta(\nu)$, we have
\[
|\bar{\partial} \tilde{w}^{\nu}|_{C^1(B_1(t))} \leq C\epsilon^2 |\nabla f_{ij} (\pi \circ \tilde{w}^{\nu} (t))|.
\]
\end{prop}
{\bf Proof:} Through straightforward computation, we have

\begin{equation}
\label{ok}
\frac{\partial z^\Re}{\partial t^\Im} = 0,\ \frac{\partial z^\Re}{\partial t^\Re} = \epsilon \nabla f_{ij} (z^\Re (t)).
\end{equation}
\[
\frac{\partial z^\Im}{\partial t^\Im} =  \epsilon \frac{\partial f_{ij}}{\partial z^\Re} (\Upsilon (\epsilon t^\Re)) = \epsilon \nabla_0 f_{ij} (z^\Re (t)).
\]
\[
\frac{\partial z^\Im}{\partial t^\Re} = \epsilon^2 \left( \left(1 - t^\Im \right) \frac{\partial^2 f_i}{\partial (z^\Re)^2} + t^\Im \frac{\partial^2 f_j}{\partial (z^\Re)^2}\right) (z^\Re (t)) \nabla f_{ij}  (z^\Re (t)).
\]

Here $\nabla$ (resp. $\nabla_0$) respects metric $g$ (resp. $g_0$). Since $|g-g_0| \leq C|J-J_0| = O(\epsilon)$, derivatives of $f_i$ are bounded, we have

\[
\left|\frac{\partial z^\Re}{\partial t^\Re} - \frac{\partial z^\Im}{\partial t^\Im}\right| \leq C\epsilon^2 |\nabla f_{ij} (z^\Re (t))|,\ \left|\frac{\partial z^\Im}{\partial t^\Re}\right| \leq C\epsilon^2 |\nabla f_{ij} (z^\Re (t))|.
\]

These estimates together with the first formula in (\ref{ok}) imply the $C^0$ part of the proposition through the following formula

\[
2\frac{\partial z}{\partial \bar{t}} = \left(\frac{\partial z^\Re}{\partial t^\Re} - \frac{\partial z^\Im}{\partial t^\Im}\right) + i\left(\frac{\partial z^\Im}{\partial t^\Re} + \frac{\partial z^\Re}{\partial t^\Im}\right).
\]

According to the observation that $\frac{\partial z}{\partial \bar{t}} (t)$ is linear on $t^\Im$ and is otherwise a function of $t^\Re$ through $z^\Re (t)$, the $C^1$ part of the proposition (or for even higher derivatives) is an easy consequence of the second formula of (\ref{ok}).
\hfill$\Box$\\

The second stage of the construction concerns $\nu \in C^0_{\rm int} (\Gamma)$. For $\textsc{x} = \textsc{x}_\nu$, we may further assume that $df_{ij} (z^\Re)$ is constant under coordinate $z$ in a small neighborhood of $\textsc{x}_\nu$ for $i,j \in I(\nu)$.\\

Assume that each interior vertex of $\Gamma$ is 3-valent. For $\nu \in C^0_{\rm int} (\Gamma)$, without loss of generality, we may assume $I(\nu) = \{1,2,3\}$. Then there are 3 legs $\nu_i \in C^1 (\Gamma)$ for $i = 1,2,3$ connecting to $\nu$ such that $I(\nu_i) = \{i,i+1\}$. let $\Gamma_\nu$ be the subgraph of $\Gamma$ consists of $\nu$ and $\nu_i \in C^1 (\Gamma)$  for $i = 1,2,3$. $\Upsilon (\Gamma_\nu)$ is a ``Y" shaped gradient tree. Let $\Upsilon_0: \Gamma_\nu \rightarrow \mathbb{C}^n$ be the standard gradient tree in our local model (\ref{ona}) with $(0,p_i) = df_i (\textsc{x}_\nu)$. Choose a constant $c>0$ and take neighborhood $U^{\mathbb{C}}_c = U^{\mathbb{C}}_c(\Upsilon_0 (\Gamma_\nu)) \subset \mathbb{C}^n$ (resp. $U^{\mathbb{R}}_c = U^{\mathbb{R}}_c(\Upsilon_0 (\Gamma_\nu)) \subset \mathbb{R}^n$) of $\Upsilon_0 (\Gamma_\nu)$ with fibration $\pi_0: U^{\mathbb{C}}_c \rightarrow U^{\mathbb{R}}_c$ such that $U^{\mathbb{C}}_c$ contains our local model holomorphic disk solution. For the convenience of discussion in the following proposition, define $\Upsilon_\epsilon (\tau) = \Upsilon (\epsilon\tau)$ for $\tau \in \Gamma_\nu$ (not to be confused with $\Upsilon_0$).\\
\begin{prop}
\label{oi}
There exist an open embedding $\phi_{\mathbb{C}}: U^{\mathbb{C}}_c \rightarrow M$ (resp. $\phi_{\mathbb{R}}: U^{\mathbb{R}}_c \rightarrow L$) such that $\phi_{\mathbb{C}}$ preserves Lagrangian fibres (namely, $\pi \circ \phi_{\mathbb{C}} = \phi_{\mathbb{R}} \circ \pi_0$) and is affine on fibres. Furthermore,
\[
\phi_{\mathbb{C}}: (U^{\mathbb{C}}_c, J_0, g_0) \rightarrow (M, J, g/\epsilon^2)
\]
is pseudo holomorphic up to $O(\epsilon)$-perturbation and is quasi-isometric, $\Upsilon_\epsilon = \phi_{\mathbb{R}} \circ \Upsilon_0$ and $\phi_{\mathbb{C}} (\Lambda_i \cap U^{\mathbb{C}}_c) \subset \Lambda^{\epsilon}_{f_i}$.
\end{prop}
{\bf Proof:} With slight abuse of notation, let $f_4, \cdots, f_{n+1}$ be functions on a small neighborhood of $\textsc{x}_\nu$ so that $df_1, \cdots, df_{n+1}$ forms a non-degenerate simplex in each fibre of $T^* L$ (trivialized via coordinate $z^\Re$) that differ by shift in different fibres. The local coordinate establishes an (isometric) identification

\[
(\mathbb{C}^n = \mathbb{R}^n + i\mathbb{R}^n, ip_i, J_0, g_0) \cong (T_{\scriptsize{\textsc{x}}_\nu} M \cong T_{\scriptsize{\textsc{x}}_\nu} L \oplus T^*_{\scriptsize{\textsc{x}}_\nu} L, \epsilon df_i (\textsc{x}_\nu), J, g /\epsilon^2).
\]

For $x\in L$, we have the affine identification $\varphi_x: \mathbb{R}^n \cong T^*_x L$ identifying $p_i$ and $\epsilon df_i (x)$. In particular, for any $\tau \in \Gamma_\nu$, we have the affine identification $\varphi_{\Upsilon_\epsilon (\tau)}: \mathbb{R}^n \cong T^*_{\Upsilon_\epsilon (\tau)} L$.\\

The tangent map $T_0\varphi_{\Upsilon_\epsilon (\tau)}: \mathbb{R}^n \cong T^*_{\Upsilon_\epsilon (\tau)} L$ is a linear identification. Using the metric $g$, we get the linear identification $[\hat{\varphi}^\Upsilon_{\mathbb{R}}]_\tau: \mathbb{R}^n \cong T_{\Upsilon_\epsilon (\tau)} L$. In such way, we define the identification $\hat{\varphi}^\Upsilon_{\mathbb{R}}: \Upsilon_0^* T\mathbb{R}^n \rightarrow \Upsilon_\epsilon^* TL$. It is straightforward to check that $\hat{\varphi}^\Upsilon_{\mathbb{R}}$ extends the identification $T\varphi^\Upsilon: T\Upsilon_0 (\Gamma_\nu) \rightarrow T\Upsilon_\epsilon (\Gamma_\nu)$, where $\varphi^\Upsilon = \Upsilon_\epsilon \circ \Upsilon_0^{-1}: \Upsilon_0 (\Gamma_\nu) \rightarrow \Upsilon_\epsilon (\Gamma_\nu)$. Since $df_{ij}$ are constants, the bundle map $\hat{\varphi}^\Upsilon_{\mathbb{R}}$ is constant under the trivialization of $TL$ according to $z$. Namely, $D\hat{\varphi}^\Upsilon_{\mathbb{R}} =0$, where $D$ is the derivative according to coordinate $z$.\\

It is straightforward to construct a quasi-isometric $C^{1,1}$-diffeomorphism $\phi_{\mathbb{R}}: (\mathbb{R}^n, g_0) \rightarrow (L, g/\epsilon^2)$ satisfying $|D^2 \phi_{\mathbb{R}}|_{g/\epsilon^2} = O(\epsilon)$ that induces $\hat{\varphi}^\Upsilon_{\mathbb{R}}$ on the germ near $\Upsilon_0 (\Gamma_\nu)$, namely, $T\phi_{\mathbb{R}} |_{\Upsilon_0} = \hat{\varphi}^\Upsilon_{\mathbb{R}}$. Then we may define $\varphi: \mathbb{C}^n \rightarrow T^*L$ as $\varphi (z_0) = \varphi_{\phi_{\mathbb{R}} (z_0^\Re)} (z_0^\Im)$. The desired $\phi_{\mathbb{C}}$ can then be defined by restricting $\varphi$ to $U^{\mathbb{C}}_c \subset \mathbb{C}^n$ together with the identification $T^* L \cong M$ near $L \subset M$.\\

By the construction of $\phi_{\mathbb{C}}$ and $\phi_{\mathbb{R}}$, we clearly have $\Upsilon_\epsilon = \phi_{\mathbb{R}} \circ \Upsilon_0$ and $\phi_{\mathbb{C}} (\Lambda_i \cap U^{\mathbb{C}}_c) \subset \Lambda^{\epsilon}_{f_i}$. To prove that $\phi_{\mathbb{C}}: (U^{\mathbb{C}}_c, J_0, g_0) \rightarrow (M, J, g/\epsilon^2)$ is pseudo holomorphic up to $O(\epsilon)$-perturbation and is quasi-isometric, we first extend $\hat{\varphi}_{\mathbb{R}} = T\phi_{\mathbb{R}}$ to $\hat{\varphi}_{\mathbb{C}}: (T\mathbb{C}^n|_{\mathbb{R}^n}, J_0) \rightarrow (TM|_L, J)$ uniquely as a pseudo holomorphic bundle map. (Here we are using the fact that the Lagrangian fibres are orthogonal to $L$, hence $J$ on $TM|_L \cong TL \oplus T^*L$ coincides with the identification of $TL$ and $T^*L$ under the metric $g|_L$.) Since $df_{ij}$ are constants, we have ${\cal H}f_i = {\cal H}f_1$ for all $i$. From the expression of $\phi_{\mathbb{C}}$, for $z_0 \in \mathbb{C}^n$, $z = \phi_\mathbb{C} (z_0)$ and $\dot{z}_0 \in T_{z_0} \mathbb{C}^n$,

\[
T_{z_0} \phi_{\mathbb{C}} (\dot{z}_0) = \hat{\varphi}_{\mathbb{C}}|_{z_0^\Re} (\dot{z}_0) + i\epsilon {\cal H}f_1|_{z^\Re} \circ \hat{\varphi}_{\mathbb{R}}|_{z_0^\Re} (\dot{z}_0^\Re).
\]

Under metric $g/\epsilon^2$, $\hat{\varphi}_{\mathbb{R}} = T\phi_{\mathbb{R}}$ being quasi-isometric implies that $\hat{\varphi}_{\mathbb{C}}$ is quasi-isometric. Hence $T\phi_{\mathbb{C}}$ is quasi-isometric and as an $O(\epsilon)$-perturbation of $\hat{\varphi}_{\mathbb{C}}$ is pseudo holomorphic up to $O(\epsilon)$-perturbation. (Here we also used the fact that $J(z) - J(z^\Re) = O(\epsilon)$ under coordinate $z$ for $z = \phi_\mathbb{C} (z_0)$ with bounded $z_0^\Im$.)
\hfill$\Box$\\

Since $\phi_{\mathbb{C}}$ is affine on transverse Lagrangian fibres according to proposition \ref{oi}, we have $\tilde{w}^{\nu_i} = \phi_{\mathbb{C}} \circ \tilde{u}^i_{\epsilon}$. Consequently, $\tilde{w} = \phi_{\mathbb{C}} \circ \tilde{u}_{\epsilon}$ defined near $\textsc{x}_\nu$ naturally coincides with $\tilde{w}^{\nu_i} (t)$ on $\Theta(\nu_i)$ for $t^\Re \geq 2\epsilon^{-\alpha}$. In such way, we can piece together the local conformal models ($\Theta(\nu_i)$ etc.) to form the global conformal model $\Theta$ (depending on the gradient tree $\Upsilon$ and $\epsilon$) and extend $\tilde{w}$ to be defined on $\Theta$. It is straightforward to see from the explicit form of $\Theta = \Theta_\epsilon$ that $\displaystyle \lim_{\epsilon \rightarrow 0} (\Theta_\epsilon, \epsilon^{-1} g_{\Theta_\epsilon}) = (\Gamma, \gamma)$ in the sense of Cheeger-Gromov. In particular, the length of $\Theta_\epsilon (\nu)$ is of order $O(\epsilon^{-1})$ for $\nu \in C^1_{\rm int} (\Gamma)$.\\
\begin{prop}
\label{om}
For $\nu \in C^0_{\rm int} (\Gamma)$ and $t\in \Theta (\nu, \frac{1}{\epsilon^\alpha})$, we have $|\bar{\partial} \tilde{w}|_{C^1(B_1(t))} \leq C\epsilon^2$.
\end{prop}
{\bf Proof:} Using the notation $z = \tilde{w} (t) = \phi_{\mathbb{C}} (z_0)$, $z_0 = \tilde{u}_{\epsilon}(t)$, we have

\[
\frac{\partial \tilde{w}}{\partial \bar{t}} = \frac{\partial \phi_{\mathbb{C}}}{\partial z_0} \frac{\partial \tilde{u}_{\epsilon}}{\partial \bar{t}} + \frac{\partial \phi_{\mathbb{C}}}{\partial \bar{z}_0} \overline{\frac{\partial \tilde{u}_{\epsilon}}{\partial t}}.
\]

Proposition \ref{oi} implies that

\[
\left| D_{z_0} \phi_{\mathbb{C}} \right| = O(\epsilon),\ \left| \frac{\partial \phi_{\mathbb{C}}}{\partial \bar{z}_0} \right| \leq C |(\phi_{\mathbb{C}})^* J - J_0| | D_{z_0} \phi_{\mathbb{C}}| \leq C \epsilon^2.
\]

The explicit expression of $\tilde{u}_{\epsilon}$ implies that

\[
|D_t \tilde{u}_{\epsilon}| = O(1),\ \left| \frac{\partial \tilde{u}_{\epsilon}}{\partial \bar{t}} \right| \leq C \epsilon^\alpha e^{-\frac{1}{\epsilon^\alpha}} \leq C\epsilon.
\]

These estimates imply the $C^0$ part of the proposition. It is easy to see that additional derivative will not change the order of each of the estimates, consequently, the $C^1$ part of the proposition.
\hfill$\Box$\\

Propositions \ref{oj} and \ref{om} together imply the following.\\
\begin{prop}
\label{ol}
Assume that the gradient tree $\Upsilon$ is non-exceptional. For $t\in \Theta$, we have $|\bar{\partial} \tilde{w}|_{C^1(B_1(t))} \leq C\epsilon^2$. For $\nu \in C^1(\Gamma)$ such that $I(\nu) = \{i,j\}$ and $t \in \Theta(\nu)$, more precisely, we have\\

\hspace{1in} $|\bar{\partial} \tilde{w}|_{C^1(B_1(t))} \leq C\epsilon^2 |\nabla f_{ij} (\pi \circ \tilde{w} (t))|.$
\hfill$\Box$\\
\end{prop}

\stepcounter{subsection}
{\bf \S \thesubsection\ The exceptional case:} Assume that $\Upsilon: \Gamma \rightarrow L$ is an exceptional gradient tree with exceptional vertex $\textsc{x}_i = \Upsilon (\textsc{p}_i)$. Let $\textsc{x}_i^\epsilon$ be the corresponding Lagrangian intersection point. It is more convenient to use the conformal model $\hat{\Theta}$ associated with the reduced tree $\hat{\Gamma}$. Let $\hat{\nu}_i \in C^1 (\hat{\Gamma})$ be the leg such that $\textsc{x}_i \in \hat{\Upsilon} (\hat{\nu}_i)$. Under $\hat{\pi}: \hat{\Gamma} \rightarrow \Gamma$, $\hat{\pi}^{-1} (\hat{\nu}_i)$ can be separated into $\nu_i, \nu_-, \nu_+ \in C^1(\Gamma)$ with $I(\nu_-) = \{i,j\}$ and $I(\nu_+) = \{i+1,j\}$. It is straightforward to see that $\hat{\Upsilon}: \hat{\nu}_i \rightarrow L$ is a $C^{1,1}$-map, where $\Upsilon = \hat{\Upsilon} \circ \hat{\pi}$.\\

In the construction of the conformal models $\Theta$ and $\hat{\Theta}$, it is convenient to construct $\Theta$ from $\hat{\Theta}$. For any $\hat{\nu} (\not= \hat{\nu}_i) \in C^*(\hat{\Gamma})$, there exists a unique $\nu \in C^*(\Gamma)$ such that $\hat{\nu} = \hat{\pi} (\nu)$. We identify $\Theta(\nu) \cong \hat{\Theta} (\hat{\nu})$. On the other hand, $\hat{\Theta} (\hat{\nu}_i)$ is a strip bounded by $\partial_i \hat{\Theta} (\hat{\nu}_i) \cup \textsc{p}_i \cup \partial_{i+1} \hat{\Theta} (\hat{\nu}_i) = \{\hat{t}^\Im =0\}$ and $\partial_j \hat{\Theta} (\hat{\nu}_i) = \{\hat{t}^\Im =1\}$, where $\hat{t}$ is the cylindrical coordinate of $\hat{\Theta} (\hat{\nu}_i)$ satisfying $\hat{t}(\textsc{p}_i)=0$. It is natural to decompose $\hat{\Theta} (\hat{\nu}_i) = \Theta(\nu_-) \cup \Theta(\nu_+)$ so that $\partial_i \hat{\Theta} (\hat{\nu}_i) = \partial_i \Theta(\nu_-)$ and $\partial_{i+1} \hat{\Theta} (\hat{\nu}_i) = \partial_{i+1} \Theta(\nu_+)$. We may then take $\Theta (\nu_i)$ to be a circular neighborhood of $\textsc{p}_i$ in $\hat{\Theta} (\hat{\nu}_i)$ with the cylindrical coordinate $t$ satisfying $\hat{t} = e^{-\pi t}$. In such way, we get the natural conformal map $\hat{\varpi}: \Theta \rightarrow \hat{\Theta}$ explicitly.\\

\pspicture(-6,-2.2)(6,1.5)
\rput(.5,0){
\pscustom[linecolor=lightgray,fillstyle=solid,fillcolor=lightgray]{
\psline(-1.8,-.4)(-1.8,.4)
\psline[liftpen=1,linearc=.2](-1.8,.4)(-1,.4)(-1,1.2)
\psline[liftpen=1](-1,1.2)(-.2,1.2)
\psline[liftpen=1,linearc=.2](-.2,1.2)(-.2,.4)(.6,.4)
\psline[liftpen=1](.6,.4)(.6,-.4)
\psline[liftpen=1](.6,-.4)(-1.8,-.4)}
\psline[linearc=.2](-1.8,.4)(-1,.4)(-1,1.2)
\psline[linearc=.2](-.2,1.2)(-.2,.4)(.6,.4)
\psline (-1.8,-.4)(.6,-.4)
\psline[linewidth=.5pt,linecolor=gray,linestyle=dashed](-1,.6)(-.2,.6)
\psline[linewidth=.5pt,linecolor=gray,linestyle=dashed](-.6,-.4)(-.6,.6)
\rput(-1.5,.65){{\footnotesize $\partial_{i+1} \Theta$}}
\rput(.2,.65){{\footnotesize $\partial_i \Theta$}}
\rput(-1.4,-.7){{\footnotesize $\partial_j \Theta$}}
\rput(-.3,-1){$\Theta$}}

\rput(-3.5,0){
\psline[linecolor=lightgray,fillstyle=solid,fillcolor=lightgray] (-1.8,-.4)(-1.8,.4) (.6,.4)(.6,-.4)(-1.8,-.4)
\psline (-1.8,.4) (.6,.4)
\psline (-1.8,-.4)(.6,-.4)
\qdisk(-.6,.4){1.5pt}
\psarc[linewidth=.5pt,linecolor=gray,linestyle=dashed](-.6,.4){.4}{180}{0}
\psline[linewidth=.5pt,linecolor=gray,linestyle=dashed](-.6,-.4)(-.6,0)
\rput(-.6,.7){$\textsc{p}_i$}
\rput(-1.4,.65){{\footnotesize $\partial_{i+1} \hat{\Theta}$}}
\rput(.2,.65){{\footnotesize $\partial_i \hat{\Theta}$}}
\rput(-1.4,-.7){{\footnotesize $\partial_j \hat{\Theta}$}}
\rput(-.6,.2){{\footnotesize $U_0$}}
\rput(-1.3,0){{\footnotesize $U_2$}}
\rput(.1,0){{\footnotesize $U_1$}}
\rput(-.3,-1){$\hat{\Theta}$}}

\rput(4.5,0){
\psline[linecolor=lightgray,fillstyle=solid,fillcolor=lightgray] (-1.8,-.4)(-1.8,.3)(-.6,.5) (.6,.3) (.6,-.4)(-1.8,-.4)
\psline (-1.8,.3)(-.6,.5) (.6,.3)
\psline (-1.8,-.4)(.6,-.4)
\qdisk(-.6,.5){1.5pt}
\psarc[linewidth=.5pt,linecolor=gray,linestyle=dashed](-.6,.5){.4}{190}{350}
\psline[linewidth=.5pt,linecolor=gray,linestyle=dashed](-.6,-.4)(-.6,.1)
\rput(-.6,.8){$\textsc{p}_i$}
\rput(-1.5,.65){{\footnotesize $\partial_{i+1} \tilde{\Theta}$}}
\rput(.2,.65){{\footnotesize $\partial_i \tilde{\Theta}$}}
\rput(-1.4,-.7){{\footnotesize $\partial_j \tilde{\Theta}$}}
\rput(-.3,-1){$\tilde{\Theta}$}}

\stepcounter{figure}
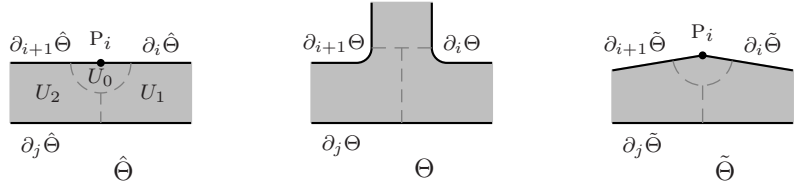
\uput{1.7}[d](0,0){Figure \thefigure: Conformal models in the exceptional case}
\endpspicture

Since $\mu (\textsc{x}_i) =0$, $\textsc{x}_i$ is a non-degenerate local minimal of $f_{i,i+1}$, ${\cal H}_{g|_L} f_{i,i+1}$ can be made proportional to $g|_L$ at $\textsc{x}_i$ by modifying $J$ near $\textsc{x}_i$. Recall that $\tilde{w}^{\nu_-}: \Theta(\nu_-) \rightarrow M$ and $\tilde{w}^{\nu_+}: \Theta(\nu_+) \rightarrow M$ satisfy $\tilde{w}^{\nu_-} (\partial_i \Theta(\nu_-)) \subset \Lambda^{\epsilon}_{f_i}$, $\tilde{w}^{\nu_+} (\partial_{i+1} \Theta(\nu_+)) \subset \Lambda^{\epsilon}_{f_{i+1}}$ and $\tilde{w}^{\nu_-} (\partial_j \Theta(\nu_-)) \cup \tilde{w}^{\nu_+} (\partial_j \Theta(\nu_+)) \subset \Lambda^{\epsilon}_{f_j}$.

\begin{prop}
\label{ou}
If ${\cal H}_{g|_L} f_{i,i+1} = c_i g|_L$ at $\textsc{x}_i$, then $\tilde{w}^{\nu_-} (\Theta(\nu_-))$ and $\tilde{w}^{\nu_+} (\Theta(\nu_+))$ pieced together form a $C^{1,1}$-surface.
\end{prop}
{\bf Proof:} It is obvious that $\tilde{w}^{\nu_-} (\Theta(\nu_-))$ and $\tilde{w}^{\nu_+} (\Theta(\nu_+))$ pieced together form a $C^{0,1}$-surface. In fact, the intersection of the two pieces is the line in $T^*_{\scriptsize{\textsc{x}}_i} L$ connecting $\textsc{x}_i^\epsilon = \epsilon df_i(\textsc{x}_i) = \epsilon df_{i+1}(\textsc{x}_i)$ and $\hat{\textsc{x}}_i^\epsilon = \epsilon df_j(\textsc{x}_i)$. More precisely, the line $\overline{\textsc{x}_i^\epsilon \hat{\textsc{x}}_i^\epsilon} \subset T^*_{\scriptsize{\textsc{x}}_i} L$ is along the direction

\[
\left. \frac{\partial \tilde{w}^{\nu_-}}{\partial t^\Im} \right|_{\scriptsize{\textsc{x}}_i} = \epsilon df_{ij} (\textsc{x}_i) = \epsilon df_{i+1,j} (\textsc{x}_i) = \left.\frac{\partial \tilde{w}^{\nu_+}}{\partial t^\Im} \right|_{\scriptsize{\textsc{x}}_i}.
\]

On the other hand,

\[
\left. \frac{\partial \tilde{w}^{\nu_+}}{\partial t^\Re} \right|_{\scriptsize{\textsc{x}}_i} - \left. \frac{\partial \tilde{w}^{\nu_-}}{\partial t^\Re} \right|_{\scriptsize{\textsc{x}}_i} = \left. \epsilon^2 (1-t^\Im) {\cal H}_{g|_L} f_{i,i+1} \nabla f_{ij} \right|_{\scriptsize{\textsc{x}}_i} = c_i \epsilon^2 (1-t^\Im) df_{ij} (\textsc{x}_i).
\]

Consequently, the tangent spaces of $\tilde{w}^{\nu_-} (\Theta(\nu_-))$ and $\tilde{w}^{\nu_+} (\Theta(\nu_+))$ coincide at their intersection, which implies the proposition.
\hfill$\Box$\\

Assume $t_1$ (resp. $t_2$) to be coordinate on $\Theta(\nu_-)$ (resp. $\Theta(\nu_+)$) with $t^\Re_1\leq 0$ (resp. $t^\Re_2\geq 0$. We may now define a different conformal model $\tilde{\Theta} (\hat{\nu}_i) = \Theta(\nu_-) \cup \Theta(\nu_+)$ with the smooth structure determined by the coordinate transformation

\[
t^\Re_2 = t^\Re_1,\ t^\Im_2 = t^\Im_1 - c_i\epsilon (1-t^\Im_1) t^\Re_1.
\]

Under the affine structure defined by $t_2$, boundaries of $\tilde{\Theta} (\hat{\nu}_i)$ are straight lines and the pullback of metric $g_{\scriptsize{\textsc{x}}^\epsilon_i}/\epsilon^2$ defined on $T_0 \tilde{\Theta} (\hat{\nu}_i)$ can be extended to a flat metric $g_{\tilde{\Theta} (\hat{\nu}_i)}$ on $\tilde{\Theta} (\hat{\nu}_i)$. Choose a complex coordinate $\tilde{t}$ on $\tilde{\Theta} (\hat{\nu}_i)$ compatible with the flat metric and coincide with $t_2$ when $t^\Im_2=0$ and $t^\Re_2\geq 0$. $\tilde{t}$ determines an embedding $\tilde{\Theta} (\hat{\nu}_i) \rightarrow \mathbb{C}$ (see figure 2).

\begin{prop}
\label{os}
$\tilde{w}^{\nu_-}$ and $\tilde{w}^{\nu_+}$ can be pieced together to form a $C^{1,1}$-map $\tilde{w}^{\hat{\nu}_i}: (\tilde{\Theta} (\hat{\nu}_i), g_{\tilde{\Theta} (\hat{\nu}_i)}) \rightarrow (M, J)$ that is pseudo holomorphic up to $O(\epsilon)$-perturbation in $C^{1,1}$.
\end{prop}
{\bf Proof:} It is straightforward to check that $\tilde{w}^{\hat{\nu}_i}$ is a $C^{1,1}$-map using proposition \ref{ou} and computations in its proof. Since $\tilde{w}^{\nu_-}$ and $\tilde{w}^{\nu_+}$ are pseudo holomorphic up to $O(\epsilon)$-perturbation in $C^{1,1}$, and the transitions among $t_1$, $t_2$ and $\tilde{t}$ are pseudo holomorphic (identity maps) up to $O(\epsilon)$-perturbation in $C^{1,1}$, we have that $\tilde{w}^{\hat{\nu}_i}$ is pseudo holomorphic up to $O(\epsilon)$-perturbation in $C^{1,1}$.
\hfill$\Box$\\

Choose $U_0$ (resp. $U_1$, resp. $U_2$) to be a small neighborhood of $B_{\frac{1}{2}}(0)$ (resp. $\Theta (\nu_-) \setminus B_{\frac{1}{2}}(0)$, resp. $\Theta (\nu_+) \setminus B_{\frac{1}{2}}(0)$) in $\hat{\Theta} (\hat{\nu}_i)$ (see figure 2). Define $\psi_0: U_0 \rightarrow \tilde{\Theta} (\hat{\nu}_i)$ as $\tilde{t} = \psi_0(\hat{t}) = \hat{t}^a$. The natural inclusion $\hat{\Theta} (\hat{\nu}_i) \rightarrow \tilde{\Theta} (\hat{\nu}_i) \rightarrow \mathbb{C}$ restricts to $\psi_1: U_1 \rightarrow \mathbb{C}$ (resp. $\psi_2: U_2 \rightarrow \mathbb{C}$). Let $\{\rho_i\}$ be a partition of unity of $\{U_i\}$.

\begin{prop}
\label{ot}
$\psi = \rho_0 \psi_0 + \rho_1 \psi_1 + \rho_2 \psi_2$ is a pseudo holomorphic map from $(\hat{\Theta} (\hat{\nu}_i), g_{\hat{\Theta} (\hat{\nu}_i)})$ to $(\tilde{\Theta} (\hat{\nu}_i), g_{\tilde{\Theta} (\hat{\nu}_i)})$ up to $O(\epsilon)$-perturbation in $C^{1,1}$.
\end{prop}
{\bf Proof:} The proposition is a simple consequence of the facts that each $\psi_i$ is pseudo holomorphic on $U_i$ up to $O(\epsilon)$-perturbation in $C^{1,1}$, and $d\psi_i = {\rm id}$ up to $O(\epsilon)$-perturbation in $C^{1,1}$ on the overlaps of different $U_i$'s.
\hfill$\Box$\\

$\tilde{w}^{\nu_-}$ and $\tilde{w}^{\nu_+}$ coincide with $\tilde{w}^{\hat{\nu}_i} \circ \psi$ near both ends of $\hat{\Theta} (\hat{\nu}_i)$. Piece all these local constructions together, we get $\tilde{w}: \hat{\Theta} \rightarrow M$. Propositions \ref{ol}, \ref{os} and \ref{ot} implies the following:

\begin{prop}
\label{or}
When the gradient tree $\Upsilon$ is exceptional, other than the estimates in proposition \ref{ol} for $\tilde{w}$ away from the exceptional vertices, we have that for $\hat{\nu}_i \in C^1(\hat{\Gamma})$ and $t \in \hat{\Theta} (\hat{\nu}_i)$ under cylindrical coordinate of $\Theta$\\

\hspace{1.5in} $|\bar{\partial} \tilde{w}|_{C^1} (t) \leq C\epsilon^2 |D\psi (t)|.$
\hfill$\Box$
\end{prop}

\se{Estimates of the linearization of $\bar{\partial}$}
Let $\mathbb{D}'$ be a smaller disk in a disk $\mathbb{D}$, $f: \mathbb{D} \rightarrow \mathbb{C}^n$ be a smooth map. The following estimate is standard. (For example, see theorem 4.15 in \cite{GT}.)

\begin{equation}
\label{oc}
|f|_{C^{1,\alpha}(\mathbb{D}')} \leq C |\bar{\partial} f|_{C^{\alpha}(\mathbb{D})} +C' |f|_{C^0(\partial \mathbb{D})}.
\end{equation}

For our application of estimating the linearization of $\bar{\partial}$, we will need to generalize this estimate to the case of almost complex structure and the corresponding free boundary problem with Lagrangian boundary condition.\\

Let $(\Sigma, \jmath)$ be a Riemann surface and ${\rm Map} (\Sigma, M)$ be the set of smooth map $w: \Sigma \rightarrow M$. The tangent space $T_w {\rm Map} (\Sigma, M)$ of ${\rm Map} (\Sigma, M)$ at $w \in {\rm Map} (\Sigma, M)$ can be identified with $\Omega^0 (\Sigma, w^* TM)$. The tangent map $Tw$ can be think of as an element of $\Omega^1 (\Sigma, w^* TM)$. In such way, $Tw$ and $Fw = J\circ Tw - Tw \circ \jmath$ define sections of the infinite dimensional vector bundle ${\cal E}$ over ${\rm Map} (\Sigma, M)$ with fibre $\Omega^1 (\Sigma, w^* TM)$ over $w \in {\rm Map} (\Sigma, M)$. $w$ is pseudo holomorphic if $Fw =0$. \\

Since our discussion will only concern a subset of ${\rm Map} (\Sigma, M)$ containing maps from $\Sigma$ to a small neighborhood $U_{\tilde{w}(\Sigma)}$ of $\tilde{w}(\Sigma)$ in $M$, where $\tilde{w}$ is the approximate solution and $\Sigma$ is a disk, without loss of generality, we may assume that $(T_M, J)$ is trivialized as a $\mathbb{C}^n$ bundle over $U_{\tilde{w}(\Sigma)}$, then the vector bundle ${\cal E}$ is trivialized with fibre $\Omega^1 (\Sigma, \mathbb{C}^n)$. The section $Fw$ can be understood as a map $F: {\rm Map} (\Sigma, M) \rightarrow \Omega^1 (\Sigma, \mathbb{C}^n)$. Similarly $T_w {\rm Map} (\Sigma, M) \cong \Omega^0 (\Sigma, \mathbb{C}^n)$. The tangent map $T_wF: \Omega^0 (\Sigma, \mathbb{C}^n) \rightarrow \Omega^1 (\Sigma, \mathbb{C}^n)$ can be computed as

\begin{equation}
\label{ox}
T_wF \dot{w} = 2i\bar{\partial} \dot{w} + B_w(Tw, \dot{w}), \mbox{ where } \dot{w} \in \Omega^0 (\Sigma, \mathbb{C}^n)
\end{equation}

and $B_{w(z)}: (T^*_z\Sigma \otimes_{\mathbb{R}} \mathbb{C}^n) \times \mathbb{C}^n \rightarrow T''_z \Sigma \otimes_{\mathbb{C}} \mathbb{C}^n$ is an $\mathbb{R}$-bilinear map smoothly depending on $w$. (Here we are using the convention $T^*_z\Sigma \otimes_{\mathbb{R}} \mathbb{C} = T'_z\Sigma \oplus T''_z \Sigma$.)\\

In our situation, since $w(\Sigma)$ is near $L$, The trivialization of $(TM, J)$ can be constructed from a trivialization of $T^*L$ through the fact that $T_z M \cong T^*_{\pi(z)}L \oplus JT^*_{\pi(z)}L$ for $z$ within the symplectic neighborhood of $L$. More precisely, we have $(TM, TL, J) \cong (\mathbb{C}^n, \mathbb{R}^n, i)$. To make the trivialization of $(T_M, J)$ isometric at $L$, one only need to ensure that the trivialization of $T^*L$ is isometric.\\

To compute (\ref{ox}), let $\{\alpha_i\}$ be the trivialization of $T^*L$, then $\{\alpha_i, J\alpha_i\}$ trivialize $T^*M$. $\dot{w}$ (resp. $Fw$) as element in $\Omega^0 (\Sigma, \mathbb{C}^n)$ (resp. $\Omega^1 (\Sigma, \mathbb{C}^n)$) has the $i$-th component

\[
[\dot{w}]_i = \langle \dot{w}, \alpha_i - iJ\alpha_i \rangle,\ [Fw]_i = (i- \jmath) (w^* J\alpha_i + i w^* \alpha_i).
\]
\[
[T_wF \dot{w}]_i = (i - \jmath) d[\dot{w}]_i + (i-\jmath) \imath_{\dot{w}} d(\alpha_i - iJ\alpha_i) = 2i\bar{\partial} [\dot{w}]_i + [B_w(Tw, \dot{w})]_i.
\]

For a disk $\mathbb{D}$ in the interior of $\Sigma$ and a smaller disk $\mathbb{D}' \subset \mathbb{D}$, we have the following interior estimate. (Here the restriction of $w$ and $\dot{w}$ to $\mathbb{D}$ is still denoted by $w$ and $\dot{w}$ for simplicity of notation.)\\
\begin{prop}
Assume that $w \in C^{1,\alpha} (\mathbb{D})$, then for $\dot{w} \in \Omega^0 (\mathbb{D}, w^* TM)$,
\begin{equation}
\label{od}
|\dot{w}|_{C^{1,\alpha} (\mathbb{D}')} \leq C_1 |\dot{w}|_{C^1(\mathbb{D})} + C_2|T_wF \dot{w}|_{C^\alpha (\mathbb{D})}.
\end{equation}
\end{prop}
{\bf Proof:} (\ref{ox}) implies that

\[
|\bar{\partial} \dot{w}|_{C^\alpha (\mathbb{D})} \leq C_1 |\dot{w}|_{C^1(\mathbb{D})} + C_2 |TF (\dot{w})|_{C^\alpha (\mathbb{D})}.
\]

This together with estimate (\ref{oc}) imply the desired estimate.
\hfill$\Box$\\

Consider the boundary version $w: (\Sigma, \partial \Sigma, j) \rightarrow (M, \Lambda, J)$ with half disks $(\mathbb{D}'_+, \partial_0 \mathbb{D}'_+) \subset (\mathbb{D}_+, \partial_0 \mathbb{D}_+) \subset (\Sigma, \partial \Sigma)$, where $\Lambda$ is a Lagrangian submanifold in $(M, \omega)$. Let $\Omega^0 (\Sigma, \vec{E}) = \{\dot{w} \in \Omega^0 (\Sigma, E): \dot{w}|_{\partial\Sigma} \in \Omega^0 (\partial\Sigma, E_\Lambda)\}$, where $\vec{E} = (E, E_\Lambda)$, $E = (w^* TM, J)$ and $E_\Lambda = w^* T \Lambda$ is a real subbundle of $E|_{\partial \Sigma}$. $\Omega^0 (\mathbb{D}_+, \vec{E})$ is the restriction of $\Omega^0 (\Sigma, \vec{E})$ to $\mathbb{D}_+$.\\
\begin{prop}
Assume that $w \in C^{1,\alpha} (\mathbb{D}_+)$, then for $\dot{w} \in \Omega^0 (\mathbb{D}_+, \vec{E})$,
\begin{equation}
\label{oe}
|D\dot{w}|_{C^{\alpha} (\mathbb{D}'_+)} \leq C_1 |\dot{w}|_{C^1(\mathbb{D}_+)} + C_2|TF (\dot{w})|_{C^\alpha (\mathbb{D}_+)}.
\end{equation}
\end{prop}
{\bf Proof:} We first assume that $E \cong \mathbb{C}^n$ and $E_\Lambda \cong \mathbb{R}^n$ under the trivialization. Then according to the proof of (\ref{od}), one only need to prove the following boundary version of (\ref{oc}) for $f: (\mathbb{D}_+, \partial_0 \mathbb{D}_+) \rightarrow (\mathbb{C}^n, \mathbb{R}^n)$,

\begin{equation}
\label{oea}
|f|_{C^{1,\alpha}(\mathbb{D}'_+)} \leq C |\bar{\partial} f|_{C^{\alpha}(\mathbb{D}_+)} +C' |f|_{C^0(\partial_1 \mathbb{D}_+)},
\end{equation}

which is a consequence of (\ref{oc}) if the extension $\tilde{f}$ of $f$ to $\mathbb{D}$ through reflection $\tilde{f}(z) = \overline{f(\bar{z})}$ for $z\in \mathbb{D}_-$ is in $C^1(\mathbb{D})$ (consequently, in $C^{1,\alpha}(\mathbb{D})$). Recall that $\partial_0 \mathbb{D}_+ = \{z\in \mathbb{D}: z^\Im=0\}$. Since $f^\Im (z^\Re)=0$, $\tilde{f} \in C^1(\mathbb{D})$ if and only if
\[
\frac{\partial f^\Re}{\partial z^\Im} (z^\Re) = \left(\frac{\partial f}{\partial \bar{z}} (z^\Re)\right)^\Im =0.
\]

In general, using a smooth bump function $\rho(s) \geq 0$ supported in $[-1,1]$ satisfying $\int_\mathbb{R} \rho (s) ds=1$, we may define

\[
\hat{f} (z) = z^\Im \int_\mathbb{R} \rho (s) \frac{\partial f^\Re}{\partial z^\Im} (z^\Re - z^\Im s) ds \mbox{ for } z \in \mathbb{D}_+.
\]

It is straightforward to check that

\[
\hat{f} (z^\Re) = 0,\ \frac{\partial \hat{f}}{\partial z^\Im} (z^\Re) = \frac{\partial f^\Re}{\partial z^\Im} (z^\Re), \mbox{ and } |\hat{f}|_{C^{1,\alpha}(\mathbb{D}_+)} \leq C |\bar{\partial} f|_{C^{\alpha}(\mathbb{D}_+)}.
\]

Consequently, $f - \hat{f}$ can be extended through reflection to a function in $C^1(\mathbb{D})$. Then (\ref{oc}) implies that

\[
|f - \hat{f}|_{C^{1,\alpha}(\mathbb{D}'_+)} \leq C |\bar{\partial} (f - \hat{f})|_{C^{\alpha}(\mathbb{D}_+)} +C' |f - \hat{f}|_{C^0(\partial_1 \mathbb{D}_+)},
\]

which implies the estimate (\ref{oea}).\\

In general, $T_wF \dot{w} \in \Omega^1 (\Sigma, \mathbb{C}^n)$ depends on the trivialization. Assume that $[T_wF \dot{w}], [T_wF \dot{w}]' \in \Omega^1 (\Sigma, \mathbb{C}^n)$ are determined by 2 different trivializations with transition matrix function $G$. Then (\ref{ox}) implies that

\[
[T_wF \dot{w}]' = [T_wF \dot{w}] + B'_w(Tw, \dot{w}), \mbox{ where } B'_w(Tw, \dot{w}) = 2i (w^*G^{-1}) \bar{\partial} (w^*G).
\]

Consequently, the estimate (\ref{oe}) for different trivializations are equivalent.
\hfill$\Box$\\

\se{Basic setting for the existence}
Start with a generic gradient tree $\Upsilon: \Gamma \rightarrow L$ with the conformal model $\Theta$, we may define a weight function $\beta = \beta_\epsilon: \Theta \rightarrow \mathbb{R}^+$. For $\nu_i \in C^1_{\rm ext} (\Gamma)$ connecting to a non-exceptional (resp. exceptional) vertex $\textsc{p}_i$, let $t$ be the cylindrical coordinate on $\Theta (\nu_i)$ such that $t^\Re (\textsc{p}_i) =+\infty$, $\beta_\epsilon$ can be defined as $\beta_\epsilon (t) = \epsilon e^{-b\epsilon t^\Re}$ (resp. $\beta_\epsilon (t) = \epsilon e^{-(1-b)\pi t^\Re}$) near $\textsc{p}_i$ for a fixed small constant $b>0$, and $\beta_\epsilon =\epsilon$ on the rest of $\Theta$ away from all the vertices. $\beta_\epsilon$ is defined so that the approximate solution $\tilde{w}$ constructed in section 7 satisfies $|D\tilde{w} (t)| \leq C\beta_\epsilon (t)$ for $t\in \Theta$ and $\epsilon$ small. Then propositions \ref{ol} and \ref{or} imply that

\begin{equation}
\label{oz}
|\bar{\partial} \tilde{w}|_{C^1} (t) \leq C\epsilon \beta_\epsilon (t), \mbox{ for } t\in \Theta.
\end{equation}

Recall that the conformal model $\Theta = \Theta_\epsilon$ satisfies $\displaystyle \lim_{\epsilon \rightarrow 0} (\Theta_\epsilon, \epsilon^{-1} g_{\Theta_\epsilon}) = (\Gamma, \gamma)$ in the sense of Cheeger-Gromov. In particular, the length of $\Theta_\epsilon (\nu)$ is $\sim \epsilon^{-1}$ for $\nu \in C^1_{\rm int} (\Gamma)$. Let ${\cal C} = Gr_{\scriptsize{\textsc{n}}}$ denote the moduli space of complete metric tree $(\Gamma, \gamma)$ that is parameterized by $\gamma_\nu$ for $\nu \in C^1_{\rm int} (\Gamma)$ with $\dim_{\mathbb{R}} {\cal C} = \textsc{n} -3$. ${\cal C}$ can also be viewed as parameterizing the moduli space ${\cal I}_{\scriptsize{\textsc{n}}}$ of the corresponding conformal models $\Theta = \Theta_\epsilon$ for each $\epsilon>0$ small. Infinitesimal deformation $\check{\jmath}$ of $(\Theta, \jmath) \in {\cal C}$ can be characterized by deformation $\dot{t}^\Re_{\check{\jmath}}$ of coordinate $t^\Re$ on $\Theta(\nu)$ for $\nu \in C^1_{\rm int} (\Gamma)$ so that $\dot{t}^\Re_{\check{\jmath}}$ is constant function near either end of $\Theta(\nu)$, and the difference of the 2 constants is $\dot{\gamma}_{\nu}$. More precisely, $\check{\jmath} \sim \frac{\partial \dot{t}^\Re_{\check{\jmath}}}{\partial t^\Re}$ pointwise. For example, one may define $\dot{t}^\Re_{\check{\jmath}} = (\dot{\gamma}_{\nu}/\epsilon) \eta (\epsilon t^\Re/ \gamma_{\nu})$ on $\Theta(\nu)$, for a smooth increasing function $\eta$ on $[0,1]$ that equals to 0 (resp. 1) near 0 (resp. 1). Elements in $T_\jmath {\cal C}$ can be equivalently represented by $\check{\jmath}$, $\dot{t}^\Re_{\check{\jmath}}$ and $\dot{\gamma}$ with equivalent norms $\| \check{\jmath} \|_{\cal C} := |\check{\jmath}|_{C^0(\Theta)} \sim |\frac{\partial \dot{t}^\Re_{\check{\jmath}}}{\partial t^\Re}|_{C^0(\Theta)} \sim |\dot{\gamma}|$. ($\check{\jmath}$ is used here to avoid confusion with $j$.) \\

Let $\tilde{{\cal B}}_1 = \tilde{{\cal B}}^\epsilon_1 (\Theta) = \tilde{{\cal B}}_1 (M, \vec{\Lambda}^\epsilon, \vec{\textsc{x}}^\epsilon, \beta_\epsilon)$ be the Banach manifold of maps $w: \Theta \rightarrow M$ such that $w (\partial_i \Theta) \subset \Lambda^\epsilon_{f_i}$ and $|Dw|_{C^\alpha_\beta (\Theta)} < \infty$, and  ${\cal B}_2 = {\cal B}^\epsilon_2 (\Theta) = (\Omega^1 (\Sigma, \mathbb{C}^n), {\|\cdot\|_{{\cal B}_2}} = |\cdot|_{C^\alpha_\beta (\Theta)})$, where

\[
|Dw|_{C^\alpha_\beta (\Theta)} = |Dw|_{C^0_\beta (\Theta)} + [Dw]_{C^\alpha_\beta (\Theta)},\ |Dw|_{C^0_\beta (\Theta)} = \sup_{t\in \Theta} \frac{|Dw|_{C^0 (B_1 (t))}}{\beta (t)},
\]
\[
[Dw]_{C^\alpha_\beta (\Theta)} = \sup_{t\in \Theta} \frac{[Dw]_{C^\alpha (B_1(t))}}{\beta (t)},\ [Dw]_{C^\alpha (B_1)} = \sup_{t_1,t_2 \in B_1} \frac{|Dw (t_1) - Dw (t_2)|}{|t_1 - t_2|^\alpha}.
\]

(Here we suppress reference to $\epsilon$ in $\beta$, $w$ and $\Theta$ for brevity.)\\

Define $F: {\cal B}_1 (= {\cal C} \times \tilde{{\cal B}}_1) \rightarrow {\cal B}_2$ such that $[F (\jmath, w)]_i = (i- \jmath) w^* (i + J) \alpha_i$. Then $F^{-1} (0) = {\cal M}_J(M, \vec{\Lambda}^\epsilon, \vec{\textsc{x}}^\epsilon)$. For $w\in \tilde{{\cal B}}_1$, $T_w \tilde{{\cal B}}_1$ can be understood as the space of sections of $\vec{E}$ that vanish at ends of $\Theta$ with the norm $\|\dot{w}\|_{{\cal B}_1} = |D\dot{w}|_{C^\alpha_\beta (\Theta)}$ for $\dot{w} \in T_w \tilde{{\cal B}}_1$, where $\vec{E} = (E, E_i)$, $E = (w^* TM, J)$ and $E_i = w^* T \Lambda^\epsilon_{f_i}$ is a real subbundle of $E|_{\partial_i \Theta}$. For $(\dot{w}, \check{\jmath}) \in T_{(w, \jmath)} {\cal B}_1$, define $\|(\dot{w}, \check{\jmath})\|_{{\cal B}_1} = \|\dot{w}\|_{{\cal B}_1} + \|\check{\jmath}\|_{\cal C}$. From the expression (\ref{ox}), it is straightforward to derive the following:\\
\begin{prop}
\label{opb}
There exists $C>0$ (independent of $\epsilon$) such that\\

\hspace{1in} $\|T_{(w,\jmath)} F(\dot{w}, \check{\jmath})\|_{{\cal B}_2} \leq C \|(\dot{w}, \check{\jmath})\|_{{\cal B}_1}.$
\hfill$\Box$
\end{prop}

Let ${\cal B}_0$ be the completion of $T_w \tilde{{\cal B}}_1$ under the weaker norm $\|\dot{w}\|_{{\cal B}_0} = |D\dot{w}|_{C^0_\beta (\Theta)}$. Define $\|(\dot{w}, \check{\jmath})\|_{{\cal B}_0} = \|\dot{w}\|_{{\cal B}_0} + \|\check{\jmath}\|_{\cal C}$. The following estimate is a direct consequence of the interior estimate (\ref{od}) and the boundary estimate (\ref{oe}).\\
\begin{prop}
\label{op}
There exists $C_1,C_2 >0$ (independent of $\epsilon$) such that\\

\hspace{1in} $\|\dot{w}\|_{{\cal B}_1} \leq C_1 \|\dot{w}\|_{{\cal B}_0} + C_2 \|T_wF\dot{w}\|_{{\cal B}_2}.$
\hfill$\Box$
\end{prop}

\begin{re}
Since ${\cal C}$ is finite dimensional, proposition \ref{op} easily implies

\begin{equation}
\label{opa}
\|(\dot{w}, \check{\jmath})\|_{{\cal B}_1} \leq C_1 \|(\dot{w}, \check{\jmath})\|_{{\cal B}_0} + C_2 \|T_{(w,\jmath)} F(\dot{w}, \check{\jmath})\|_{{\cal B}_2}.
\end{equation}
\end{re}

For our choice of $\beta$, the $C^0_\beta (\Theta)$-norm and estimate (\ref{oz}) can be rewritten as

\begin{equation}
\label{og}
|D \dot{w}|_{C^0_\beta (\Theta)} \sim \sup_{t\in \Theta} \frac{|D \dot{w} (t)|}{\beta (t)},\ \|F (\tilde{w})\|_{{\cal B}_2} = O(\epsilon).
\end{equation}

\begin{theorem}[Inverse function theorem]
\label{oh}
Assume that ${\cal B}_1$, ${\cal B}_2$ are Banach spaces, $F: {\cal B}_1 \rightarrow {\cal B}_2$ is a map such that $\frac{\partial F}{\partial h}(0): {\cal B}_1 \rightarrow {\cal B}_2$ is invertible. $\|[\frac{\partial F}{\partial h}(0)]^{-1}\| \leq C$. Assume further that we can find $r_0>r>0$ such that for $h \in U_{{\cal B}_1}(r_0)$,

\[
\left\|\frac{\partial F}{\partial h}(h) - \frac{\partial F}{\partial h}(0)\right\| \leq \frac{1}{2C},\ \ \|F(0)\|_{{\cal B}_2} \leq \frac{r}{2C}.
\]

Then there exists a unique $h_0 \in U_{{\cal B}_1}(r)$ such that $F(h_0)=0$.
\hfill$\Box$
\end{theorem}

To relate the deformation of pseudo holomorphic polygons to deformation of gradient trees in our situation, it is necessary to compute $T_wF \dot{w}$ more precisely. Using the identification $(TM, J) \cong \pi^* TL \otimes_{\mathbb{R}} \mathbb{C}$, one can define a connection ${\cal D}$ on $TM$ based on the Levi-Civita connection on $(TL, g|_L)$. Let $2i\bar{\partial}_{\cal D} = (i-\jmath){\cal D}$, we have

\begin{prop}
Assume that $J$ on $M \cong T^*L$ coincides with the canonical almost complex structure $J_{\rm can}$ determined by the metric on $L$. Then
\begin{equation}
\label{of}
T_wF \dot{w} = 2i\bar{\partial}_{\cal D} \dot{w} + O(|w||Tw||\dot{w}|) + O(|Fw||\dot{w}|).
\end{equation}
\end{prop}
{\bf Proof:} $\{\alpha^i\}$ induce the fibre affine coordinate (trivialization) $y: M (\cong T^*L) \rightarrow \mathbb{R}^n$. We have the coordinate $(x,y)$ on $M$, where $x$ is the composition of a coordinate on $L$ and $\pi: M \rightarrow L$. Under the Levi-Civita connection on $L$, ${\cal D} \alpha^i = -\omega_j^i \alpha^j$ and $d \alpha^i = -\omega_j^i \wedge \alpha^j$. It is easy to check that the horizontal vector field on $M (\cong T^*L)$ over $\frac{\partial}{\partial x^i}$ on $L$ is

\[
v_i = \frac{\partial}{\partial x^i} + y_k \left\langle \omega_j^k, \frac{\partial}{\partial x^i}\right\rangle \frac{\partial}{\partial y_j}.
\]

$\{v_i, \frac{\partial}{\partial y_i}\}$ forms the dual basis of $\{\alpha_i, J\alpha_i\}$. Using the formula

\[
d\alpha (u,v) = u\langle \alpha, v\rangle - v\langle \alpha, u\rangle - \langle \alpha, [u,v] \rangle,
\]

\[
[v_i, v_j] = O(|y|),\ [\frac{\partial}{\partial y_i}, v_j] = \left\langle \omega_k^i, \frac{\partial}{\partial x^j}\right\rangle \frac{\partial}{\partial y_k}.
\]

we have

\[
d J\alpha^i \left(\frac{\partial}{\partial y_j}, \frac{\partial}{\partial y_k} \right) =0,\ d J\alpha^i \left(v_j, v_k\right) = O(|y|).
\]
\[
d J\alpha^i \left(\frac{\partial}{\partial y_j}, v_k \right) = v_k (g^{ij}) - g^{il} \left\langle \omega_l^j, \frac{\partial}{\partial x^k}\right\rangle = g^{jl} \left\langle \omega_l^i, \frac{\partial}{\partial x^k}\right\rangle.
\]

Consequently, $d J\alpha^i = - \omega_j^i \wedge J\alpha^j + O(|y|)$.

\[
[T_wF \dot{w}]_i = (i - \jmath) d[\dot{w}]_i + (i-\jmath) \imath_{\dot{w}} d(\alpha_i - iJ\alpha_i)
\]
\[
= (i - \jmath) (d[\dot{w}]_i + \omega_j^i [\dot{w}]_j) - \left\langle \omega_j^i, \dot{w}\right\rangle (i-\jmath) (\alpha_i - iJ\alpha_i) + O(|w||Tw||\dot{w}|)
\]

\hspace{1in} $= 2i[\bar{\partial}_{\cal D} \dot{w}]_i + O(|w||Tw||\dot{w}|) + O(|Fw||\dot{w}|)$.
\hfill$\Box$\\
\begin{prop}
Assume $J = J_{\rm can}$ on $L \subset M$. Then
\begin{equation}
\label{oy}
T_wF \dot{w} = 2i\bar{\partial}_{\cal D} \dot{w} + O((|w| + \epsilon)|Tw||\dot{w}|) + O(|Fw||\dot{w}|) + O(|Tw||D\dot{w}|).
\end{equation}
\end{prop}
{\bf Proof:} (\ref{of}) was proved through (\ref{ox}) by showing the estimate

\begin{equation}
\label{oyb}
[B_w(Tw, \dot{w})]_i = (i - \jmath) \omega_j^i [\dot{w}]_j + O(|w||Tw||\dot{w}|) + O(|Fw||\dot{w}|)
\end{equation}

for $J = J_{\rm can}$. To prove (\ref{oy}), we only need to show similar estimate as (\ref{oyb}) for the new $J$. The only terms in (\ref{oyb}) that are affected by the change of $J$ are $[\dot{w}]_j$, $Fw$ and the term $\imath_{\dot{w}} d(J\alpha_i)$ in $[B_w(Tw, \dot{w})]_i = (i-\jmath) \imath_{\dot{w}} d(\alpha_i - iJ\alpha_i)$. Our assumption implies that $J - J_{\rm can} = O(y)$. Consequently, the changes caused by $[\dot{w}]_j$ and $Fw$ can be absorbed by the term $O(|w||Tw||\dot{w}|)$. The change to the remaining term $\imath_{\dot{w}} d(J\alpha_i)$ can be controlled by $O(|w||Tw||\dot{w}|) + O([\dot{w}]^y |Tw|)$. Since $\dot{w}$ is tangent to $\Lambda^\epsilon_{f_i}$ at $\partial_i \Theta$, we have $[\dot{w}]^y|_{\partial \Theta} = O(\epsilon |\dot{w}|)$, $[\dot{w}]^y = O(\epsilon |\dot{w}|) + O(|D\dot{w}|)$. We have proved (\ref{oyb}) for the new $J$ with the additional term $O([\dot{w}]^y |Tw|) = O(\epsilon |\dot{w}| |Tw|) + O(|D\dot{w}| |Tw|)$. With such estimate, the proposition is a consequence of (\ref{ox}).
\hfill$\Box$\\
\begin{re}
\label{oyc}
For $\bar{\partial}_{\cal D}$ to coincide with the usual $\bar{\partial}$, the trivialization of $T^*L$ has to be parallel, which can always be done at one point in $L$, and is usually not possible globally.
\end{re}

\begin{co}
When deformation $\check{\jmath}$ in $T_\jmath {\cal C}$ is considered, we have

\begin{equation}
\label{oya}
T_{(w, \jmath)} F (\dot{w}, \check{\jmath}) = 2i\bar{\partial}_{\cal D} \check{w} + O((|w| + \epsilon)|Tw||\check{w}| + |Fw||\check{w}| + |Tw||D\check{w}|),
\end{equation}

where $\check{w} = \dot{w} - \dot{t}^\Re_{\check{\jmath}} \frac{\partial w}{\partial t^\Re}$.
\end{co}
{\bf Proof:} (\ref{oya}) is different from (\ref{oy}) in the support of $\check{\jmath}$, which is in the union of $\Theta(\nu)$ for $\nu \in C^1_{\rm ext} (\Gamma)$. Locally on such  $\Theta(\nu)$, by change of coordinates according to $\dot{t}^\Re_{\check{\jmath}}$, the deformation of $(w, \jmath)$ according to $(\dot{w}, \check{\jmath})$ has the same effect on $F$ as the deformation of $w$ according to $\check{w}$. Namely, $T_{(w, \jmath)} F (\dot{w}, \check{\jmath}) = T_w F \check{w}$. In such way, (\ref{oya}) is a consequence of (\ref{oy}).
\hfill$\Box$

\begin{re}
There is one important subtle point: For the estimate (\ref{oya}), it is necessary to require $w$ to be at least $C^2$-bounded on the support of $\check{\jmath}$. (\ref{oya}) will only be applied when $w$ is the approximate solution, which clearly satisfies this condition.\\
\end{re}

{\samepage
\stepcounter{subsection}
{\bf \S \thesubsection\ Limiting estimates}\\
\nopagebreak

For $\epsilon_k \rightarrow 0$,} let $w_k \in {\cal B}^{\epsilon_k}_1 (\Theta_k)$ with uniform bound and $\jmath_k$ be the complex structure on $\Theta_k$. In this section, for the purpose of applying (\ref{oy}), we need to assume the condition $Fw_k = o(\epsilon_k)$, which is clearly satisfied for $w_k$ pseudo holomorphic ($Fw_k=0$) or for the approximate solution $\tilde{w}_k$ constructed from a generic gradient tree $\Upsilon$ ($F\tilde{w}_k = O(\epsilon_k^2)$).

\begin{prop}
\label{ob}
Assume that $\|(\dot{w}_k, \check{\jmath}_k)\|_{{\cal B}_0} =1$ and $\lim \|T_{(w_k,\jmath_k)} F(\dot{w}_k, \check{\jmath}_k)\|_{{\cal B}_2} = 0$, then for $c_1>0$ there exists $c_2>0$ such that $\frac{\epsilon_k}{{\beta_k}(t)} |\dot{w}_k (t)| \geq c_2$ for $k$ large and $t\in \Theta_k$ satisfying $\frac{1}{{\beta_k}(t)} |D\dot{w}_k (t)| \geq c_1$ and $\check{\jmath}_k (t)=0$.
\end{prop}
{\bf Proof:} If the proposition is not true, then there exists $t_k\in \Theta_k$ such that $\check{\jmath}_k (t_k)=0$, $\lim \frac{1}{{\beta_k}(t_k)} |D\dot{w}_k (t_k)| = c_0 \geq c_1$ and $\lim \frac{\epsilon_k}{{\beta_k}(t)} |\dot{w}_k (t_k)| = 0$.\\

Assume $w_k(t_k) \rightarrow \textsc{x}_o \in L$. We use a trivialization of $TM$ near $\textsc{x}_o$ induced from a trivialization of $T^*L$ that is parallel at $\textsc{x}_o \in L$ with respect to the metric on $L$, so that $T_zM$ is identified with $T_{\scriptsize\textsc{x}_o} M$ for $z$ near $\textsc{x}_o$. Then $\dot{w}_k$ can be understood as a map $\Theta_k \rightarrow T_{\scriptsize\textsc{x}_o} M$. Assume that $(\Theta_o, 0, g_{\Theta_o}) = \lim (\Theta_k, t_k, g_{\Theta_k})$ in the sense of Cheeger-Gromov. We take the following limit

\[
\dot{w}_o (t) = \lim_{k\rightarrow +\infty} Q_k (\dot{w}_k (t)), \mbox{ where } Q_k (\dot{z}) = \frac{1}{{\beta_k}(t_k)} \left[\dot{z} - \pi_*\dot{w}_k (t_k)\right] \mbox{ for } \dot{z} \in T_{\scriptsize\textsc{x}_o} M.
\]

Since $\check{\jmath}_k (t_k)=0$, $[T_{(w_k,\jmath_k)} F(\dot{w}_k, \check{\jmath}_k)] (t_k) = [T_{w_k} F\dot{w}_k] (t_k)$. By (\ref{oy}), remark \ref{oyc} and our choice of trivialization of $TM$, we have $\bar{\partial} Q_k (\dot{w}_k (t)) = O(\epsilon_k + {\rm Dist}(w_k(t),\textsc{x}_o))$ near $t_k$. Consequently, the limit $\dot{w}_o: \Theta_o \rightarrow (T_{\scriptsize\textsc{x}_o} M,J) \cong \mathbb{C}^n$ is holomorphic.\\

For $i$ and $\tilde{t}_k \in \partial_i \Theta_k$ such that $|\tilde{t}_k - t_k|$ is uniformly bounded and $\tilde{t}_o = \lim \tilde{t}_k \in \partial_i \Theta_o$, $T_{w_k (\tilde{t}_k)} \Lambda_{f_i}^{\epsilon_k} \subset T_{w_k (\tilde{t}_k)} M$ is defined by $\dot{z}^\Im = \epsilon_k {\cal H}f_i(w_k (\tilde{t}_k)) \dot{z}^\Re$. Hence

\[
\dot{w}_k (\tilde{t}_k) = \epsilon_k {\cal H}f_i(w_k (\tilde{t}_k)) \pi_* \dot{w}_k (\tilde{t}_k) + \pi_* \dot{w}_k (\tilde{t}_k),
\]
\[
Q_k (\dot{w}_k (\tilde{t}_k)) = \frac{\epsilon_k}{{\beta_k} (t_k)} {\cal H}f_i(w_k (\tilde{t}_k)) \pi_* \dot{w}_k (\tilde{t}_k) + \frac{1}{{\beta_k} (t_k)} \pi_* [\dot{w}_k (\tilde{t}_k) - \dot{w}_k (t_k)].
\]

Since $\lim \frac{\epsilon_k}{{\beta_k}(t)} |\dot{w}_k (t_k)| = 0$ and $\frac{1}{{\beta_k}(t_k)} |\dot{w}_k (\tilde{t}_k) - \dot{w}_k (t_k)| \leq C_{|\tilde{t}_k - t_k|}|\dot{w}_k|_{C_{\beta_k}^1 (\Theta_k)}$, we have

\[
\dot{w}_o (\tilde{t}_o) = \lim_{k\rightarrow +\infty} Q_k (\dot{w}_k (\tilde{t}_k)) = \lim_{k\rightarrow +\infty} \frac{1}{{\beta_k} (t_k)} \pi_* [\dot{w}_k (\tilde{t}_k) - \dot{w}_k (t_k)] \in T_{\scriptsize\textsc{x}_o} L.
\]

Consequently, under the identification $(T_{\scriptsize\textsc{x}_o} M, T_{\scriptsize\textsc{x}_o} L, J) \cong (\mathbb{C}^n, \mathbb{R}^n, i)$, $\dot{w}_o: \Theta_o \rightarrow \mathbb{C}^n$ can be viewed as a holomorphic map such that $|D\dot{w}_o| \leq \beta_o$ and $\dot{w}_o (\partial \Theta_o) \subset \mathbb{R}^n$, where $\beta_o = \lim \beta_k/\beta_k(t_k)$. $\beta_o$ is unbounded on $\Theta_o$ only when $t_k \rightarrow \textsc{p}_i$ corresponds to an exceptional vertex, then $\beta_o (t) = e^{-(1-b)\pi t^\Re}$ for the cylindrical coordinate $t$ on $\Theta_o$. By proposition \ref{ed}, $\dot{w}_o$ has to be a constant map. Consequently $D\dot{w}_o = 0$.\\

On the other hand, $\lim \frac{1}{{\beta_k} (t_k)} D\dot{w}_k (t_k) = D\dot{w}_o (0)$ Consequently,
\[
|D\dot{w}_o (0)| = \lim_{k\rightarrow +\infty} \left|\frac{1}{{\beta_k} (t_k)} D\dot{w}_k (t_k)\right| = c_0>0,
\]

which is a contradiction to $D\dot{w}_o = 0$.
\hfill$\Box$\\
\begin{co}
\label{oo}
Assume that $\|(\dot{w}_k, \check{\jmath}_k)\|_{{\cal B}_0} =1$ and $\lim \|T_{(w_k,\jmath_k)} F(\dot{w}_k, \check{\jmath}_k)\|_{{\cal B}_2} = 0$, then for any $t_k\in \Theta_k$ such that $\frac{1}{{\beta_k}(t_k)} |D\dot{w}_k (t_k)| \geq c_1$ for certain $c_1>0$, $\lim w_k(t_k)$ can not be an exceptional vertex $\textsc{x}_i$.
\end{co}
{\bf Proof:} If the corollary is not true, then $\lim w_k(t_k) = \textsc{x}_i$ is an exceptional vertex. Hence $\check{\jmath}_k (t_k)=0$ and $\epsilon_k \hat{t}_k \rightarrow 0$, where $\hat{t}_k = e^{-\pi t_k}$ is the coordinate of $t_k \in \hat{\Theta}_k (\hat{\nu}_i)$ under the cylindrical coordinate of $\hat{\Theta}_k (\hat{\nu}_i)$ that vanishes at $\textsc{p}_i$, and $\hat{\nu}_i \in C^1(\hat{\Gamma})$ is the leg that contains $\textsc{p}_i$. Proposition \ref{ob} then implies that there exists $c_2>0$ such that $\frac{\epsilon_k}{{\beta_k}(t_k)} |\dot{w}_k (t_k)| \geq c_2$ for large $k$. On the other hand,

\[
|\dot{w}_k (t_k)| = \int_{\scriptsize\textsc{p}_i}^{t_k} |D\dot{w}_k (t)| dt \leq \int_{\scriptsize\textsc{p}_i}^{t^\Re_k} \beta_k (t) dt \leq C \max (\beta_k (t_k), \epsilon_k |\hat{t}_k|)
\]

implies that $\lim \frac{\epsilon_k}{\beta_k (t_k)} |\dot{w}_k (t_k)| = 0$, which is a contradiction.
\hfill$\Box$\\

Let $\Upsilon(\tau)$ be a gradient line of $f_{ij}$. The infinitesimal deformation $\dot{\Upsilon}$ of gradient lines of $f_{ij}$ at $\upsilon$ satisfies the equation

\[
\frac{D \dot{\Upsilon} (\tau)}{d\tau} = [{\cal H}f_{ij} (\Upsilon) \dot{\Upsilon}] (\tau).
\]

If the coordinate $\tau$ is allowed to deviate from the canonical parametrization of gradient line with the infinitesimal deformation $\dot{\tau}$, then the equation is

\begin{equation}
\label{oqc}
\frac{D \check{\Upsilon} (\tau)}{d\tau} = [{\cal H}f_{ij} (\Upsilon) \check{\Upsilon}] (\tau), \mbox{ where } \check{\Upsilon} (\tau) = \left[\dot{\Upsilon} - \dot{\tau} \frac{d\Upsilon}{d\tau} \right] (\tau).
\end{equation}

Recall that each non-degenerate directed leg of $\Gamma$ under a gradient tree $\Upsilon: \Gamma \rightarrow L$ possess a canonical coordinate up to a constant, and when $\Upsilon$ is generic, only possible degenerate legs are those external legs mapped to exceptional vertices. The infinitesimal deformation of a generic gradient tree $\Upsilon$ can be described by a pair $(\dot{\Upsilon}, \dot{\tau})$ satisfying (\ref{oqc}) on each leg $\nu$ of $\Gamma$ with $I(\nu) = \{i,j\}$, where $\dot{\Upsilon}$ is a piecewise smooth continuous section of $\Upsilon^* TL$ and $\dot{\tau}$ is a smooth function up to constant on each internal leg of $\Gamma$, and $\dot{\tau}=0$ on other legs of $\Gamma$. A different description of the same infinitesimal deformation of $\Upsilon$ is completely determined by a different $\dot{\tau}$ with the same variation $\dot{\gamma}_\nu$ on each $\nu \in C^1_{\rm int}$. For this reason, to avoid ambiguity, one may standardize $\dot{\tau} (\tau) = \dot{\gamma}_{\nu} \eta (\tau/ \gamma_{\nu})$ on $\nu \in C^1_{\rm int}$. The description of infinitesimal deformation of a non-generic gradient tree (not needed in this paper) is more complicated because the topology of $\Gamma$ may change under deformation.

\begin{re}
An infinitesimal deformation of a generic gradient tree $\Upsilon$ is trivial if $\dot{\Upsilon}=0$ for the unique standardized representation $(\dot{\Upsilon}, \dot{\tau})$. Consequently, $\dot{\tau} =0$ and $\dot{\gamma}=0$. Therefore, if one can show that $\dot{\Upsilon}\not=0$ or $\dot{\tau}\not=0$ for a standardized representation $(\dot{\Upsilon}, \dot{\tau})$, then the corresponding infinitesimal deformation of the generic gradient tree $\Upsilon$ is non-trivial.
\end{re}

Let $w_k =\tilde{w}_k$ be the approximate solution for a generic gradient tree $\Upsilon$. For the purpose of applying (\ref{oya}), in addition to the condition $Fw_k = o(\epsilon_k)$, we need to assume $D^2w_k = o(\epsilon_k)$, which is clearly satisfied for the approximate solution $\tilde{w}_k$. Assume that $\|(\dot{w}_k, \check{\jmath}_k)\|_{{\cal B}_0} =1$ and $\lim \|T_{(w_k,\jmath_k)} F(\dot{w}_k, \check{\jmath}_k)\|_{{\cal B}_2} = 0$. We may define

\begin{equation}
\label{oqa}
\dot{\Upsilon} (\tau) = \lim_{k\rightarrow +\infty} \dot{w}_k (\tau/\epsilon_k),\ \dot{\tau} (\tau) = \dot{\gamma}_{\nu} \eta (\tau/ \gamma_{\nu}),\ \dot{\gamma}_\nu = \lim_{k\rightarrow +\infty} \epsilon_k \dot{\gamma}^k_\nu.
\end{equation}

\begin{prop}
\label{oq}
$(\dot{\Upsilon}, \dot{\tau})$ represents an infinitesimal deformation of the generic gradient tree $\Upsilon$. The convergence in (\ref{oqa}) is uniformly $C^1$ on any compact subset in the smooth part of $\Gamma$.
\end{prop}
{\bf Proof:} $\|\dot{w}_k\|_{{\cal B}_0} \leq 1$ implies that $\dot{\Upsilon}$ is a continuous section of $\Upsilon^* TL$. For the first part of the proposition, it only remains to be shown that (\ref{oqc}) is satisfied on each leg $\nu$ of $\Gamma$ with $I(\nu) = \{i,j\}$.\\

Assume $t_k \rightarrow \tau_o$ under $(\Theta_k, \epsilon_k g_{\Theta_k}) \rightarrow (\Gamma, \gamma)$. Then

\[
\check{\Upsilon} (\tau_o +x) = \lim_{k\rightarrow +\infty} \check{w}_k \left(t_k + \frac{x}{\epsilon_k}\right), \mbox{ where } \check{w}_k = \dot{w}_k - \dot{t}^\Re_{\check{\jmath}_k} \frac{\partial w_k}{\partial t^\Re}.
\]

According to lemma \ref{mg}, if (\ref{oqc}) is not true at $\tau=\tau_o$, then by possibly passing to subsequence and adjusting $t_k$, there exists $c_1>0$ such that for $k$ large,

\begin{equation}
\label{oqb}
\left|\frac{1}{\epsilon_k} D\check{w}_k (t_k) - [{\cal H}f_{ij} (\Upsilon) \check{\Upsilon}] (\tau_o)\right| \geq c_1.
\end{equation}

We use a trivialization of $TM$ near $\textsc{x}_o$ induced from a trivialization of $T^*L$ that is parallel at $\textsc{x}_o \in L$ with respect to the metric on $L$, so that $T_zM$ is identified with $T_{\scriptsize\textsc{x}_o} M$ for $z$ near $\textsc{x}_o$. Then $\dot{w}_k$ can be understood as a map $\Theta_k \rightarrow T_{\scriptsize{\textsc{x}}_o} M$. We take the following limit

\[
\dot{w}_o (t) = \lim_{k\rightarrow +\infty} Q_k (\check{w}_k (t^\Re_k + t)), \mbox{ where } Q_k (\dot{z}) = \frac{1}{\epsilon_k} \left[\dot{z} - \check{w}_k (t^\Re_k)\right] \mbox{ for } \dot{z} \in T_{\scriptsize\textsc{x}_o} M.
\]

For $\tilde{t}_k$ such that $|\tilde{t}_k - t_k|$ is uniformly bounded, $\lim Q_k (T_{\pi(w_k (\tilde{t}_k))} L) = T_{\scriptsize\textsc{x}_o} L$, $\lim Q_k (T_{w_k (\tilde{t}^\Re_k)} \Lambda_{f_i}^{\epsilon_k}) = p_i + T_{\scriptsize\textsc{x}_o} L$ and $\lim Q_k (T_{w_k (\tilde{t}^\Re_k+i)} \Lambda_{f_j}^{\epsilon_k}) = p_j + T_{\scriptsize\textsc{x}_o} L$, where $p_i = [{\cal H}f_i (\Upsilon) \check{\Upsilon}] (\tau_o)$ and $p_j = [{\cal H}f_j (\Upsilon) \check{\Upsilon}] (\tau_o)$. By our assumption of $F w_k = o(\epsilon_k)$ and $D^2w_k = o(\epsilon_k)$, we have $|\check{w}_k|=O(1)$, $|D\check{w}_k|=o(1)$, $|Tw_k| = O(\epsilon_k)$ and $|w_k| = O(\epsilon_k)$. By (\ref{oya}), remark \ref{oyc} and our choice of trivialization of $TM$, $\bar{\partial} Q_k (\check{w}_k (t^\Re_k + t)) = o(1)$ uniformly for bounded $t$.\\

Consequently, using the identification $(T_{\scriptsize\textsc{x}_o} M, T_{\scriptsize\textsc{x}_o} L, J) \cong (\mathbb{C}^n, \mathbb{R}^n, i)$, the limit $\dot{w}_o: \Theta_o \rightarrow \mathbb{C}^n$ is holomorphic with bounded derivative, where $\Theta_o$ is a strip with 2 boundary components $\partial_i \Theta_o$ and $\partial_j \Theta_o$ such that $\dot{w}_o (\partial_i \Theta_o) \subset p_i + \mathbb{R}^n$ and $\dot{w}_o (\partial_j \Theta_o) \subset p_j + \mathbb{R}^n$. By proposition \ref{ec}, $\dot{w}_o$ has to be linear. Consequently $D\dot{w}_o = [{\cal H}f_{ij} (\Upsilon) \check{\Upsilon}] (\tau_o)$.\\

On the other hand, $\frac{1}{\epsilon_k} D\check{w}_k (t_k) \rightarrow D\dot{w}_o (0)$ Consequently,
\[
\lim_{k\rightarrow +\infty} \left|\frac{1}{\epsilon_k} D\check{w}_k (t_k) - [{\cal H}f_{ij} (\Upsilon) \check{\Upsilon}] (\tau_o)\right| = 0,
\]

which is a contradiction to (\ref{oqb}). According to remark \ref{mga}, when $\tau_o$ varies in $\Gamma$, the contradiction also imply that $\check{w}_k (\tau/\epsilon_k)$ $C^1$-converge to $\check{\Upsilon} (\tau)$ uniformly on any compact subset in the smooth part of $\Gamma$. This imply the second part of the proposition.
\hfill$\Box$\\

\se{Invertibility of the linearized operator}
\stepcounter{subsection}
{\bf \S \thesubsection\ Injectivity:} Since $\Theta$ is non-compact, ${\cal B}_1 \rightarrow {\cal B}_0$ is not a compact embedding. But we still have\\
\begin{prop}
\label{oa}
If $\Upsilon$ is a rigid gradient tree, then there exists $C >0$ (independent of $\epsilon$) such that
\[
\|(\dot{w}, \check{\jmath})\|_{{\cal B}_1} \leq C \|T_{(\tilde{w}, \jmath)}F (\dot{w}, \check{\jmath})\|_{{\cal B}_2}.
\]
\end{prop}
{\bf Proof:} If the assertion of the proposition is not true, then by proposition \ref{op}, there exist $\epsilon_k \rightarrow 0$ and the corresponding approximate solutions $\tilde{w}_k: \Theta_k \rightarrow M$ and $(\dot{w}_k, \check{\jmath}_k) \in T_{(\tilde{w}_k, \jmath_k)} {\cal B}_1$ such that $\|(\dot{w}_k, \check{\jmath}_k)\|_{{\cal B}_0} =1$ and $\|T_{(\tilde{w}_k, \jmath_k)}F (\dot{w}_k, \check{\jmath}_k)\|_{{\cal B}_2} \rightarrow 0$. One can find $t_k \in \Theta_k$, such that $\textsc{x} = \displaystyle \lim w_k(t_k)$ exists and

\begin{equation}
\label{oab}
1 \geq \frac{|D \dot{w}_k(t_k)|}{{\beta_k}(t_k)} + |\check{\jmath}_k (t_k)| \geq 1-\epsilon_k.
\end{equation}

According to the first part of proposition \ref{oq}, the limit of $(\dot{w}_k, \check{\jmath}_k)$ will produce a deformation $(\dot{\Upsilon}, \dot{\tau})$ of the gradient tree $\Upsilon$. Since $\Upsilon$ is rigid, this is a contradiction if $(\dot{\Upsilon}, \dot{\tau})$ is non-trivial.\\

If $\textsc{x} \not= \textsc{x}_i$ for any $i$, then proposition \ref{ob} implies that $(\dot{\Upsilon}, \dot{\tau})$ is non-trivial if $\check{\jmath}_k (t_k)=0$ for $k$ large. If $\check{\jmath}_k (t_k) \not=0$ for $k$ large, then $\beta_k(t_k) = \epsilon_k$ for $k$ large and $(\dot{\Upsilon}, \dot{\tau})$ is non-trivial by (\ref{oab}) and the second part of proposition \ref{oq}.\\

Corollary \ref{oo} implies that the only remaining case is when $\textsc{x} = \textsc{x}_i$ is a non-exceptional vertex. Then we have the pointed Gromov-Hausdorff limit

\[
(\Gamma^\circ, \tau^\circ, \gamma^\circ) = \lim_{k \rightarrow +\infty} (\Theta_k, t_k, \epsilon_k^2g_{\Theta_k}) \cong (\mathbb{R}, 0, g_{\mathbb{R}}).
\]

With $\delta_k = {\rm Dist}_k (\textsc{x}_i, w_k(t_k))$, we have the pointed Gromov-Hausdorff limit

\[
(T_{\scriptsize\textsc{x}_i}M, 0, \textsc{v}_o, g) = \lim_{k \rightarrow +\infty} (M, \textsc{x}_i, w_k(t_k), g/\delta^2_k),
\]

where $|\textsc{v}_o|=1$. $u_k = \frac{\epsilon_k\delta_k}{\beta(t^\circ_k) \|\dot{w}_k\|_{{\cal B}_0}} \dot{w}_k$ can be considered as a map $u_k: (\Theta_k, t_k, \epsilon_kg_{\Theta_k})$ $ \rightarrow (M, w_k(t_k), g/\delta^2_k)$. Then

\[
|Du_k (t)| = \frac{|D \dot{w}_k(t)|}{\beta_k (t_k)\|\dot{w}_k\|_{{\cal B}_0}} \leq \frac{\beta_k(t)}{\beta_k(t_k)} \mbox{ and } |Du_k (t_k)| \geq 1-\epsilon_k.
\]

Let $u^\circ: (\Gamma^\circ, \tau^\circ, \gamma^\circ) \rightarrow (T_{\scriptsize\textsc{x}_i}M, \textsc{v}_o, g)$ (resp. $\beta^\circ (\tau)$) be the limit of $u_k$ (resp. $\frac{1}{\epsilon_k} {\beta_k}(t_k + \tau/\epsilon_k)$) as $k \rightarrow +\infty$. Then $u^\circ: (\mathbb{R}, 0, g_{\mathbb{R}}) \rightarrow (T_{\scriptsize\textsc{x}_i}M, \textsc{v}_o, g)$ satisfies

\begin{equation}
\label{oaa}
\left|\frac{du^\circ (\tau)}{d\tau}\right| \leq \frac{\beta^\circ (\tau)}{\beta^\circ (0)} = e^{b\tau} \mbox{ and } |Du^\circ (0)| = 1.
\end{equation}

Proposition \ref{oq} implies that $u^\circ$ satisfies

\[
\frac{du^\circ (\tau)}{d\tau} = {\cal H}f_{i,i+1} (\textsc{x}_i) u^\circ (\tau).
\]

It is elementary to observe that no solution of this ordinary differential equation satisfy the condition (\ref{oaa}).
\hfill$\Box$\\

\stepcounter{subsection}
{\bf \S \thesubsection\ Surjectivity:} With proposition \ref{oa}, the surjectivity of the linearized operator $T_{(w,\jmath)} F$ can be viewed as a consequence of the assertion that the operator $T_{(w,\jmath)} F$ has index 0. Since the index of $T_{(w,\jmath)} F$ is the virtual dimension of the moduli space ${\cal M}_J(M, \vec{\Lambda}^{\epsilon}, \vec{\textsc{x}}^{\epsilon})$ of the pseudo holomorphic polygon $w$, any definition of Fukaya category should have verified this assertion for the norm of its choice. Then following certain standard procedure in functional analysis, one should be able to show that the operator $T_{(w,\jmath)} F$ also has index 0 under our norm. In the following, we will prove that $T_{(w,\jmath)} F$ is surjective based on the surjectivity proved in \cite{FO}. By taking the weight function $\beta$ to be the constant $\epsilon$, we may define the corresponding space $C^{1,\alpha}_\epsilon (\Theta)$. The corresponding Sobolev space $W_\epsilon^{1,p} (\Theta)$ is used in \cite{FO}. Assume $\alpha \leq \frac{p-n}{p}$. We start with the following generalization of (\ref{oc}) and (\ref{oea}).

\begin{prop}
For $\dot{w} \in \Omega^0 (D, \mathbb{C}^n)$,
\begin{equation}
\label{pa}
|\dot{w}|_{C^{1,\alpha} (D')} \leq C_1 |\dot{w}|_{W^{1,p}(D)} + C_2|\bar{\partial} \dot{w}|_{C^\alpha (D)}.
\end{equation}
For boundary version, assume $\dot{w} \in \Omega^0 (D_+, \mathbb{C}^n)$ and $\dot{w}|_{\partial_0 D_+} \in \mathbb{R}^n$,
\begin{equation}
\label{pb}
|D\dot{w}|_{C^{\alpha} (D'_+)} \leq C_1 |\dot{w}|_{W^{1,p}(D_+)} + C_2|\bar{\partial} \dot{w}|_{C^\alpha (D_+)}.
\end{equation}
\end{prop}
{\bf Proof:} These estimates are standard. We will show that (\ref{pa}) and (\ref{pb}) are direct consequence of (\ref{oc}) and (\ref{oea}) without any need of $L^p$-estimate for $\bar{\partial}$.\\

First assume that $\dot{w}$ is compactly supported in $D$. Let $\dot{w}_\delta = \dot{w} * \rho_\delta$. Then

\[
\bar{\partial} \dot{w}_\delta = (\bar{\partial} \dot{w}) * \rho_\delta.
\]

Assume $|\dot{w}_\delta|_{C^1(D)} = C_\delta \rightarrow \infty$. Let $u_\delta = \dot{w}_\delta /C_\delta$, then $|u_\delta|_{C^1(D)} = 1$ and $|\bar{\partial} u_\delta|_{C^\alpha(D)} \leq |\bar{\partial} \dot{w}|_{C^\alpha (D)} /C_\delta \rightarrow 0$. By (\ref{oc}), $|u_\delta|_{C^{1,\alpha} (D)} \leq C$. By passing to subsequence, $u_\delta$ has a limit $u_0$ that is holomorphic, compactly supported in $D$ and $|u_0|_{C^1(D)} =1$, which is a contradiction. Consequently, $|\dot{w}_\delta|_{C^1(D)} \leq C$. By (\ref{oc}), $|\dot{w}_\delta|_{C^{1,\alpha} (D)} \leq C$. Let $\delta \rightarrow 0$, we have $|\dot{w}|_{C^{1,\alpha} (D)} \leq C|\bar{\partial} \dot{w}|_{C^\alpha (D)}$, which implies (\ref{pa}).\\

In general, let $\rho$ be a bump function supported inside $D$ and $\rho|_{D'}=1$. Then we have

\[
|\dot{w}|_{C^{1,\alpha} (D')} \leq C|\bar{\partial} \rho \dot{w}|_{C^\alpha (D)} \leq C_1 |\dot{w}|_{C^\alpha (D)} + C_2|\bar{\partial} \dot{w}|_{C^\alpha (D)}.
\]

Since $\alpha \leq \frac{p-n}{p}$, by Sobolev inequality, we have $|\dot{w}|_{C^\alpha (D)} \leq C |\dot{w}|_{W^{1,p}(D)}$. Consequently, we have (\ref{pa}).\\

(\ref{pb}) is a direct consequence of (\ref{pa}) if $\bar{\partial} \dot{w}$ is in $C^\alpha (D)$ when $f$ is extended by reflection from $D_+$ to $D$. In general, one need to use $\hat{f}$ defined in the proof of (\ref{oea}) and follow the arguments in the proof of (\ref{oea}).
\hfill$\Box$\\
\begin{co}
\label{pc}
$\dot{w} \in W_\epsilon^{1,p} (\Theta)$ and $TF \dot{w} \in C^{\alpha}_\epsilon (\Theta)$ imply that $\dot{w} \in C^{1,\alpha}_\epsilon (\Theta)$.
\end{co}
{\bf Proof:} This is a consequence of formula (\ref{ox}) and estimates (\ref{pa}) and (\ref{pb}).
\hfill$\Box$\\
\begin{prop}
\label{pe}
For $\epsilon$ small, $T_{(w,\jmath)} F: C^{1,\alpha}_\beta \rightarrow C^\alpha_\beta$ is surjective if $T_{(w,\jmath)} F: W^{1,p}_\epsilon \rightarrow L^p_\epsilon$ is surjective.
\end{prop}
{\bf Proof:} For any $y \in C^\alpha_\beta$, assume that there exists $\dot{w} \in W^{1,p}_\epsilon$ such that $y = T_{(w,\jmath)} F \dot{w}$. (See sections 6 and 7 in \cite{FO}.) Since $\dot{w} \in W_\epsilon^{1,p}$, corollary \ref{pc} implies that $\dot{w} \in C^{1,\alpha}_\epsilon (\Theta)$. We need to show that $\dot{w} \in C^1_\beta (\Theta)$ for $\epsilon$ small.\\

If not so, then there exists $\epsilon_k \rightarrow 0$ and $t_k \in \Theta_k$ such that $w_k(t_k) \rightarrow \textsc{x}_i$ and

\begin{eqnarray*}
\frac{|D w_k (t)|}{\beta_k (t)} \leq \frac{|D w_k (t_k)|}{\beta_k (t_k)} &\mbox{ for }& t^\Re \leq t^\Re_k.\\
|D w_k (t)| \leq (1+\epsilon_k) |D w_k (t_k)| &\mbox{ for }& t^\Re \geq t^\Re_k.
\end{eqnarray*}

When $\textsc{x}_i$ is a non-exceptional vertex, by the same argument as in the last part of the proof of proposition \ref{oa}, we have $u^\circ: (\mathbb{R}, 0, g_{\mathbb{R}}) \rightarrow (T_{\scriptsize\textsc{x}_i}M, \textsc{v}_o, g)$ satisfying the linear ordinary differential equation

\[
\frac{du^\circ (\tau)}{d\tau} = {\cal H}f_{i,i+1} (\textsc{x}_i) u^\circ (\tau),
\]

and the bounds

\begin{equation}
\label{pd}
\left|\frac{du^\circ (\tau)}{d\tau}\right| \leq e^{-b\tau} \mbox{ for } \tau \leq 0 \mbox{ and } |Du^\circ (\tau)| \leq |Du^\circ (0)| = 1 \mbox{ for } \tau \geq 0.
\end{equation}

It is elementary to observe that no solution of this ordinary differential equation satisfies the bounds (\ref{pd}).\\

We omit the case of exceptional vertex because the corresponding result was not proved in \cite{FO}.
\hfill$\Box$\\

This more involved proposition is necessary only because the weight function $\beta$ in this paper is incompatible with the weight function in \cite{FO}. We are now in the position to prove the existence of the pseudo holomorphic polygons.

\begin{theorem}
\label{pf}
For a rigid gradient tree $\Upsilon$ and $\epsilon$ small, let $\tilde{w}: (\tilde{\Theta}, \tilde{\jmath}) \rightarrow M$ (depending on $\epsilon$) be the approximate solution constructed in section 7, then there exists a unique pseudo holomorphic polygon $w: (\Theta, \jmath) \rightarrow (M,J)$ in a small $\|\cdot\|_{{\cal B}_1}$-neighborhood of $(\tilde{w}, \tilde{\jmath})$.
\end{theorem}
{\bf Proof:} By (\ref{oz}) or (\ref{og}), we have $\| F(\tilde{w}, \tilde{\jmath}) \|_{{\cal B}_2} = O(\epsilon)$. By (\ref{ox}), it is easy to check that $\|T_{(\tilde{w}, \tilde{\jmath})} F - T_{(w, \jmath)} F\| \leq C\|(\tilde{w}, \tilde{\jmath}) - (w, \jmath)\|_{{\cal B}_1}$. By propositions \ref{oa} and \ref{pe} (more generally \cite{disk}), we have $\|(T_{(\tilde{w}, \tilde{\jmath})} F)^{-1}\| \leq C$. Then by the inverse function theorem \ref{oh}, there exists a unique $w: (\Theta, \jmath) \rightarrow (M,J)$ in a small $\|\cdot\|_{{\cal B}_1}$-neighborhood of $(\tilde{w}, \tilde{\jmath})$ such that $F(w, \jmath)=0$.
\hfill$\Box$

\begin{re}
Our treatment of surjectivity here is clearly not very satisfying in term of our $C^{1,\alpha}$-method, because the result from sections 6 and 7 in \cite{FO} used a lot of quite non-trivial arguments and results with Sobolev norms. Besides, we also need to handle the case of exceptional gradient tree that was not discussed in \cite{FO}. In a forthcoming paper (\cite{disk}), we will provide a completely different geometric proof of the surjectivity using the $C^{1,\alpha}$-method we developed in this paper.\\

The main difficulty of our index theorem lie in the fact that $\Theta$ is complete instead of compact. If $\Theta$ is a compact disk with smooth boundary, the corresponding index theorem is quite standard. Let $\Lambda = \bigcup \Lambda^\epsilon_{f_i}$. The basic idea in \cite{disk} is to construct a family of smooth Lagrangian submanifolds $\tilde{\Lambda}_{\tilde{\epsilon}}$ as smoothing of the singular Lagrangian $\Lambda$ through so-called Lagrangian surgery so that $\tilde{\Lambda}_{\tilde{\epsilon}} \rightarrow \Lambda$ when $\tilde{\epsilon} \rightarrow 0$. When $\tilde{\epsilon}$ is small, we can establish 1-1 correspondence between pseudo holomorphic disks with boundary in $\tilde{\Lambda}_{\tilde{\epsilon}}$ and pseudo holomorphic polygon with boundary in $\Lambda$. Furthermore, the index is preserved under the limiting process. In this way, we can prove that the index of $T_{(w,\jmath)} F$ is zero, since the corresponding index for holomorphic disk is zero by standard results. Clearly, this idea has other more general implications that fit into a more general theory that we will discuss in \cite{disk}. We should note that Lagrangian surgery originated from \cite{P} was discussed extensively in the last chapter of \cite{FOOO} for apparently somewhat different purpose.
\end{re}

\se{Uniqueness}
Notice that theorem \ref{pf} contains a uniqueness statement. Therefore to show the uniqueness of the pseudo holomorphic polygon near a rigid gradient tree $\Upsilon$, it is only necessary to show that the pseudo holomorphic polygon $(w,\jmath)$ is in a small $\|\cdot\|_{{\cal B}_1}$-neighborhood of the approximate solution $(\tilde{w}, \tilde{\jmath})$ constructed from $\Upsilon$.\\
\begin{prop}
\label{ov}
Assume that $w_k$ is pseudo holomorphic, then $\|Dw_k\|_{{\cal B}_2} =O(1)$.
\end{prop}
{\bf Proof:} It is straightforward to see that the estimate away from vertices is a consequence of proposition \ref{mc} and estimates (\ref{nf}), (\ref{nk}). Hence if the proposition is not true, then by possibly taking subsequence, there exist a vertex $\textsc{x}_i$ and $t_k \in \Theta_k (\nu_i)$ such that $\textsc{x}_i = \lim w_k (t_k)$ and

\begin{eqnarray*}
\frac{|D w_k (t)|}{\beta_k (t)} \leq \frac{|D w_k (t_k)|}{\beta_k (t_k)} &\mbox{ for }& t^\Re \leq t^\Re_k.\\
|D w_k (t)| \leq (1+\epsilon_k) |D w_k (t_k)| &\mbox{ for }& t^\Re \geq t^\Re_k.
\end{eqnarray*}

Let $\epsilon'_k := |D w_k (t_k)|$. Through the same arguments as in the proof of proposition \ref{mb}, applying to the sequence of maps  $w_k: (\Theta_k, t_k, g_{\Theta_k}) \rightarrow (M, w_k (t_k), g/(\epsilon'_k)^2)$, using Cheeger-Gromov convergence, we arrive at $w_o: (\Theta_o, t_o, g_{\Theta_o}) \rightarrow (T_{\scriptsize\textsc{x}_i} M, 0, g_{\scriptsize\textsc{x}_i})$ such that $w_o$ is holomorphic, $w_o (\partial_i \Theta_o) \subset p_i + T_{\scriptsize\textsc{x}_i} L$ (resp. $w_o (\partial_j \Theta_o) \subset p_j + T_{\scriptsize\textsc{x}_i} L$) and $|D w_o (t)| \leq |D w_o (t_o)| =1$ when $\textsc{x}_i$ is a non-exceptional vertex, or when $\textsc{x}_i$ is an exceptional vertex

\[
|D w_o (t_o)| =1,\ |D w_o (t)| \leq e^{(1-b)\pi (t^\Re_o - t^\Re)} \mbox{ for } t^\Re \leq t^\Re_o,\ |D w_o (t)| \leq 1 \mbox{ for } t^\Re \geq t^\Re_o.
\]

Let $p_{ij} := p_j -p_i$ and $\tilde{\epsilon}_k = \epsilon'_k /\epsilon_k$. It is straightforward to see that

\[
|p_{ij}| \leq |Dw_o (t_o)|=1, \mbox{ and } p_{ij} = \lim_{k \rightarrow +\infty}  \frac{1}{\tilde{\epsilon}_k} df_{ij}(w_k (t_k)).
\]

Let $\tilde{w}_o = w_o(t^\Re_o) + t^\Im_o p_{ij} - t^\Re_o Jp_{ij}$. $\check{w}_o = w_o - \tilde{w}_o$ satisfies $\check{w}_o (\partial \Theta_o) \subset T_{\scriptsize\textsc{x}_i} L$ and $|D \check{w}_o (t)| \leq C$ when $\textsc{x}_i$ is a non-exceptional vertex, or when $\textsc{x}_i$ is an exceptional vertex

\[
|D \check{w}_o (t)| \leq Ce^{(1-b)\pi (t^\Re_o - t^\Re)} \mbox{ for } t^\Re \leq t^\Re_o,\ |D \check{w}_o (t)| \leq C \mbox{ for } t^\Re \geq t^\Re_o.
\]

According to proposition \ref{ed}, $\check{w}_o$ is a constant map. Consequently

\[
|p_{ij}| = |D\tilde{w}_o| = |Dw_o (t_o)|=1,\ {\rm Dist}_g (\textsc{x}_i, w_k (t_k)) \sim |df_{ij}(w_k (t_k))| \sim \tilde{\epsilon}_k.
\]

Let $f_{ij}^k (x) = \frac{1}{\tilde{\epsilon}_k^2} f_{ij} (\tilde{\epsilon}_k x)$, $\omega^k = \frac{1}{\tilde{\epsilon}_k^2} \omega$, $g^k = \frac{1}{\tilde{\epsilon}_k^2} g$. $(M, \textsc{x}_i, w_k (t_k), g^k, f_{ij}^k)$ converge to $(M^\circ \cong T_{\scriptsize\textsc{x}_i} M, 0, \textsc{v}, g^\circ, f_{ij}^\circ)$, where $f_{ij}^\circ$ is the quadratic approximation of $f_{ij}$ at $\textsc{x}_i$.\\

When $\textsc{x}_i$ is an exceptional vertex, for any $R>0$,

\[
|D w_k (t)|_{g^k} = \tilde{\epsilon}_k |D w_k (t)|_g \leq C_R\epsilon_k \mbox{ for } t^\Re \geq t^\Re_k -R.
\]

With this estimate, by similar arguments as in section 6, $w_k: (\Theta_k, t_k, \epsilon_k^2 g_{\Theta_k}) \rightarrow (M, w_k (t_k), g^k)$ converge to a gradient line $w^\circ: (\Gamma^\circ, \tau_o, \gamma^\circ) \rightarrow (T_{\scriptsize\textsc{x}_i} M, \textsc{v}, g^\circ)$ of the function $f_{ij}^\circ$ flowing toward $0 \in T_{\scriptsize\textsc{x}_i} M$. Since 0 is a local minimal of $f_{ij}^\circ$, this is a contradiction.\\

When $\textsc{x}_i$ is a non-exceptional vertex, for any $R>0$,

\[
|D w_k (t)|_{g^k} = \tilde{\epsilon}_k |D w_k (t)|_g \leq C_R\epsilon_k \mbox{ for } t^\Re \geq t^\Re_k -R/\epsilon_k.
\]

With this estimate, by similar arguments as in section 6, $w_k: (\Theta_k, t_k, \epsilon_k^2 g_{\Theta_k}) \rightarrow (M, w_k (t_k), g^k)$ converge to a gradient line $w^\circ: (\Gamma^\circ, 0, \gamma^\circ) \rightarrow (T_{\scriptsize\textsc{x}_i} M, \textsc{v}, g^\circ)$ of the function $f_{ij}^\circ$ flowing toward 0. Satisfying the estimate

\[
\left|\frac{d w^\circ}{d\tau} (0)\right| =1,\ \left|\frac{d w^\circ}{d\tau} (\tau)\right| \leq e^{-b\tau} \mbox{ for } \tau \leq 0,\ \left|\frac{d w^\circ}{d\tau} (\tau)\right| \leq 1 \mbox{ for } \tau \geq 0.
\]

Since $f_{ij}^\circ$ is the quadratic approximation of $f_{ij}$ at $\textsc{x}_i$, and $\textsc{x}_i$ is a non-degenerate critical point of $f_{ij}$, there exists $c>0$ such that $|df_{ij}^\circ (x)| = |({\cal H}f_{ij}|_{\scriptsize{\textsc{x}}_i}) x| \geq c|x|$ for $x \in T_{\scriptsize\textsc{x}_i} M$. Then no such gradient line exists for $b<c$.
\hfill$\Box$\\
\begin{prop}
\label{ow}
Assume that $w_k$ is pseudo holomorphic, then $\lim |\check{w}_k|_{{\cal B}_1} =0$, where $\check{w}_k = w_k - \tilde{w_k}$.
\end{prop}
{\bf Proof:} If the proposition is not true, then there exist $c_1>0$ and $t_k \in \Theta_k$ such that

\[
|D\check{w}_k|_{C^\alpha (B_1(t_k))} \geq c_1 \beta_k(t_k).
\]

We first need to show that $\lim \frac{\epsilon_k}{{\beta_k}(t_k)} |\check{w}_k (t_k)| =0$. This is obvious if $\beta_k (t_k) \sim \epsilon_k$ for $k$ large. Otherwise, $t_k \in \Theta_k (\nu_i)$ and $\lim w_k (t_k) = \textsc{x}_i$ for some $i$. If $\textsc{x}_i$ is exceptional, then as in the proof of the corollary \ref{oo}, we have that

\[
|\check{w}_k (t_k)| = \int_{\scriptsize\textsc{p}_i}^{t_k} |D\check{w}_k (t)| dt \leq \int_{\scriptsize\textsc{p}_i}^{t^\Re_k} \beta_k (t) dt \leq C \max (\beta_k (t_k), \epsilon_k |\hat{t}_k|)
\]

implies that $\frac{\epsilon_k}{\beta_k (t_k)} |\check{w}_k (t_k)| \rightarrow 0$.\\

If $\textsc{x}_i$ is not exceptional, by proposition \ref{ov}, $|D \check{w}_k (t_k)| \leq C\epsilon_k e^{-2b\epsilon_k t^\Re_k}$. Hence $|\check{w}_k (t_k)| \leq Ce^{-2b\epsilon_kt^\Re_k}$, $\frac{\epsilon_k}{\beta_k (t_k)} |\check{w}_k (t_k)| \leq Ce^{-b\epsilon_kt^\Re_k} \rightarrow 0$.\\

Under a local coordinate of $M$ near $\textsc{x}_i$, $\check{w}_k$ can be understood as a map $\Theta_k \rightarrow \mathbb{C}^n$. Assume that $(\Theta_o, 0, g_{\Theta_o}) = \lim (\Theta_k, t_k, g_{\Theta_k})$ in the sense of Cheeger-Gromov. We take the following limit

\[
\check{w}_o (t) = \lim_{k\rightarrow +\infty} Q_k (\check{w}_k (t)), \mbox{ where } Q_k (z) = \frac{1}{{\beta_k}(t_k)} \left(z - [\check{w}_k (t_k)]^\Re \right) \mbox{ for } z \in \mathbb{C}^n.
\]

(\ref{oz}) implies that $|\bar{\partial} Q_k (\check{w}_k (t))| = O(\epsilon_k)$, hence the limit $\check{w}_o$ is holomorphic. Apply proposition \ref{ov}, we can show that $|D\check{w}_k|_{C^{\alpha'}_{\beta_k} (\Theta_k)} = O(1)$ for certain $\alpha'>\alpha$. Consequently

\[
|D\check{w}_o|_{C^\alpha (B_1(t_k))} = \lim_{k\rightarrow +\infty} \frac{1}{{\beta_k}(t_k)} |D\check{w}_k|_{C^\alpha (B_1(t_k))} \geq c_1.
\]

For $j$ and $\tilde{t}_k \in \partial_j \Theta_k$ such that $|\tilde{t}_k - t_k|$ is uniformly bounded and $\tilde{t}_o = \lim \tilde{t}_k \in \partial_j \Theta_o$,

\[
\check{w}_k (\tilde{t}_k) = [\check{w}_k (\tilde{t}_k)]^\Re + i\epsilon_k df_j([\check{w}_k (\tilde{t}_k)]^\Re),
\]
\[
Q_k (\check{w}_k (\tilde{t}_k)) =  \frac{1}{{\beta_k} (t_k)} [\check{w}_k (\tilde{t}_k) - \check{w}_k (t_k)]^\Re + i \frac{\epsilon_k}{{\beta_k} (t_k)} df_j([\check{w}_k (\tilde{t}_k)]^\Re).
\]

Since $\frac{\epsilon_k}{{\beta_k}(t)} |\check{w}_k (t_k)| \rightarrow 0$ and $\frac{1}{{\beta_k}(t_k)} |\check{w}_k (\tilde{t}_k) - \check{w}_k (t_k)| \leq C_{|\tilde{t}_k - t_k|}|\check{w}_k|_{C_{\beta_k}^1 (\Theta_k)}$, we have

\[
\check{w}_o (\tilde{t}_o) = \lim_{k\rightarrow +\infty} Q_k (\check{w}_k (\tilde{t}_k)) = \lim_{k\rightarrow +\infty} \frac{1}{{\beta_k} (t_k)} [\check{w}_k (\tilde{t}_k) - \check{w}_k (t_k)]^\Re \in \mathbb{R}^n.
\]

Consequently, when $\beta_o = \lim \beta_k/\beta_k(t_k)$ is bounded on $\Theta_o$, $\check{w}_o: \Theta_o \rightarrow \mathbb{C}^n$ can be viewed as a holomorphic map such that $|D\check{w}_o| \leq \beta_o$ and $\check{w}_o (\partial \Theta_o) \subset \mathbb{R}^n$. $\beta_o$ is unbounded only when $t_k \rightarrow \textsc{p}_i$ corresponds to an exceptional vertex. Then $\beta_o (t) = e^{-(1-b)\pi t^\Re}$ for the cylindrical coordinate $t$ on $\Theta_o$. By proposition \ref{ed}, $\check{w}_o$ has to be a constant map. Consequently $|D\check{w}_o|_{C^\alpha (B_1(t_k))} =0$, which is a contradiction.
\hfill$\Box$\\

Combine theorem \ref{mo} and proposition \ref{ow}, we have

\begin{theorem}
\label{qa}
For any $r_0>0$ there exists $\epsilon_0>0$ such that when $\epsilon\leq \epsilon_0$, for any pseudo holomorphic polygon $w: (\Theta, \jmath) \rightarrow (M, J)$ such that $(w,\jmath) \in {\cal M}_J(M, \vec{\Lambda}^{\epsilon}, \vec{\textsc{x}}^{\epsilon})$ with ${\rm Ind} (\vec{f}, \vec{\textsc{x}}) =0$, there exists a rigid gradient tree $\Upsilon: (\Gamma, \gamma) \rightarrow (L,g|_L)$ such that ${\rm Dist}_{{\cal B}_1} ((w,\jmath), (\tilde{w}, \tilde{\jmath})) \leq r_0$.
\hfill$\Box$
\end{theorem}

\ifx\undefined\bysame
\newcommand{\bysame}{\leavevmode\hbox to3em{\hrulefill}\,}
\fi

\end{document}